\documentclass{article}
\usepackage[numbers,sort&compress]{natbib}
\usepackage{amsfonts}
\usepackage{amsmath}
\usepackage[english]{babel}
\usepackage{bbm}
\usepackage{color}
\usepackage{graphicx}
\usepackage[utf8x]{inputenc}
\usepackage{lmodern}
\usepackage{multicol}
\usepackage{comment}
\usepackage{caption}
\usepackage{subcaption}

\usepackage{hyperref}
\usepackage{authblk}

\newcommand{\rP}{r_\P}

\newcommand{\calE}{\mathcal{E}}

\newcommand{\calI}{\mathcal{I}}

\newcommand{\calL}{\mathcal{L}}

\newcommand{\calP}{\mathcal{P}}

\newcommand{\Ns}{N}

\renewcommand{\P} {E}
\newcommand  {\E} {e} 

\newcommand{\hh}{h}
\newcommand{\Th}{\Omega_{\hh}}

\newcommand{\hP}{\hh_{\P}}

\newcommand{\mP}{\ABS{\P}}

\newcommand{\mE}{\ABS{\E}}

\newcommand{\ABS}    [1]{\vert#1\vert}

\newcommand{\EOD}{\end{document}}

\newcommand{\RSet}[1][]{{\mathbb{R}}^{#1}}
\newcommand{\NSet}[1][]{{\mathbb{N}}^{#1}}
\newcommand{\Network}{\Omega}
\newcommand{\Fracture}{F}
\newcommand{\Trace}{S}
\newcommand{\FractureIndices}{\mathcal{I}}
\newcommand{\TraceIndices}{\mathcal{M}}
\newcommand{\TracePair}{\sigma}
\newcommand{\MeshParameter}{h}
\newcommand{\OriginalMesh}{\widetilde{\Omega}_{\MeshParameter}}
\newcommand{\Mesh}{\Omega_{\MeshParameter}}
\newcommand{\Diameter}{h}
\newcommand{\MeshCellTwoD}{E}
\newcommand{\MeshCellOneD}{e}

\newcommand{\VemOrder}{k}
\newcommand{\Polynomial}{\mathbb{P}}

\newcommand{\scalar}[3]{\left({#2}, {#3}\right)_{#1}}
\newcommand{\BilinearForm}[3]{#1 \left({#2}, {#3}\right)}

\newcommand{\HilbertNorm}[2][]{\left\lvert \left\lvert{#2}\right\rvert \right\rvert_{#1}}
\newcommand{\seminorm}[2][]{\left\lvert{#2}\right\rvert_{#1}}
\newcommand{\SobolevL}[1]{{\rm L}^{#1}}
\newcommand{\SobolevH}[2][]{{\rm H}_{#1}^{#2}}
\newcommand{\SobH}[2]{\mathop{{\SobolevH{#1}}{\rm(}#2{\rm)}}\nolimits}
\newcommand{\SobHo}[2]{\mathop{{\SobolevH[0]{#1}}{\rm(}#2{\rm)}}\nolimits}
\newcommand{\SobL}[2][2]{\mathop{\SobolevL{#1}(#2)}\nolimits}
\newcommand{\Projection}{\Pi}
\newcommand{\Boundary}[1]{\partial {#1}}
\newcommand{\SolutionSpace}{V}
\newcommand{\Basis}{\mathcal{M}}

\newcommand{\Matrixize}[1]{\mathbb{#1}}

\newcommand{\Diffusion}{\Matrixize{K}}
\newcommand{\HydraulicHead}{u}

\newcommand{\Energy}{\mathcal{E}}
\newcommand{\VemConvergence}{\alpha}

\newcommand{\AgglomerationParameter}{\lambda}
\newcommand{\DiscreteError}[1][]{e_{\MeshParameter#1}}

\title{Mesh Quality Agglomeration algorithm for the Virtual Element Method applied to Discrete Fracture Networks}

\author[1]{T.~Sorgente}
\author[2]{F.~Vicini}
\author[2]{S.~Berrone}
\author[1]{S.~Biasotti}
\author[1]{G.~Manzini}
\author[1]{M.~Spagnuolo}

\affil[1]{Istituto di Matematica Applicata e Tecnologie Informatiche `E.~Magenes', Consiglio Nazionale delle Ricerche, Italy}
\affil[2]{ Dipartimento di Scienze Matematiche `G.L.~Lagrange', Politecnico di Torino, Italy,}

\begin{document}
\maketitle


\begin{abstract}
We propose a quality-based optimization strategy to reduce the total number of degrees of freedom associated to a discrete problem defined over a polygonal tessellation with the Virtual Element Method.
The presented \emph{Quality Agglomeration} algorithm relies only on the geometrical properties of the problem polygonal mesh, agglomerating groups of neighboring elements.
We test this approach in the context of fractured porous media, in which the generation of a global conforming mesh on a Discrete Fracture Network leads to a considerable number of unknowns, due to the presence of highly complex geometries and the significant size of the computational domains.
We show the efficiency and the robustness of our approach, applied independently on each fracture for different network configurations, exploiting the flexibility of the Virtual Element Method in handling general polygonal elements.
\end{abstract}


\section{Introduction}
\label{sec:introduction}
Over the last fifty years, computer simulations of Partial Differential Equations (PDEs) have dramatically increased their impact on research, design, and production and are now an indispensable tool for modeling and analyzing phenomena from physics, engineering, biology, and medicine. 
The most popular techniques for such computer-based simulations (e.g., the Finite Element Method) rely on suitable descriptions of geometrical entities, such as the computational domain and its properties, which are generally encoded by a mesh, or tessellation.
Despite decades of research and significant results, techniques for generating meshes with good geometrical properties are still studied and developed.
Current meshing algorithms typically produce an initial mesh, which we expect to be dominated by well-shaped elements, and optionally perform optimization steps to maximize the geometrical quality of the elements.
We can broadly group such optimization strategies into two categories \cite{lo2014finite}: geometrical methods (also called \textit{smoothing} or \textit{untangling}) that keep the mesh connectivity fixed while changing only the locations of the mesh vertices \cite{erten2009mesh, knupp2012introducing, vartziotis2008mesh}, and topological methods (also \textit{re-meshing}) that keep the vertices fixed while changing only the connectivity \cite{alliez2008recent, misztal2009tetrahedral}.
In both cases, the critical point is the definition of the concept of \textit{quality} for the mesh elements.
This concept is still a discussed topic in the numerical analysis community, particularly when associated with generic polygonal meshes, as proved by the vast number of very different attempts to define a universal quality indicator
\cite{knupp2001algebraic, stimpson2007verdict, chalmeta2013measuring, zunic2004new, huang2020anisotropic, sorgente2021role, sorgente2021indicator, BERRONE2022103770}.
As a result, it is easier to define mesh optimization strategies that are only related to the optimization of the size of the elements, or the number of incident edges, instead of a generic concept of quality.

In this work, we present an algorithm called \textit{Quality Agglomeration} that automatically optimizes the number of degrees of freedom (DOFs) in the discrete problem defined over a polygonal tessellation.
The algorithm acts as a topological optimization method, agglomerating groups of neighboring elements.
In particular, it minimizes an energy functional based on the mesh quality indicator introduced in \cite{sorgente2021role}, which involves the use of the Virtual Element Method (VEM)~\cite{BPVEM} in the numerical discretization.
The algorithm performs energy minimization via the graph-cut technique (also known as maximum-flow or minimum-cut algorithm \cite{goldberg1988new}), an efficient graph-based technique aimed at splitting a graph into two or more parts, minimizing a certain energy.

We also present an application of our algorithm in the context of fractured porous media.
The flow simulation in underground fractured porous media is a fundamental task for a huge number of applications, such as carbon capture processes, the monitoring of aquifers, or the contaminant transport in the subsoil.
Among the existing approaches for the simulation of underground phenomena in fractured media \cite{fumagalli2019numerical}, the Discrete Fracture Network (DFN) model is particularly suited for scenarios in which the permeability of the surrounding media is negligible with respect to the fracture transmissivities.
Indeed, in DFN simulations the porous rock matrix is neglected and the fluid is confined to move on the rock discontinuities, called \textit{fractures}, which act as channels for the flow.
The DFN model is characterized by a dimensional reduction of each fracture to a planar polygon immersed in a three-dimensional space, and the coupling conditions are imposed on fracture segment intersections \cite{martin2005modeling}.
Due to the deficiency of direct measurements of the subsoil properties, the position, orientation, size, and aperture of each polygonal fracture are generated according to different probability distributions.
Furthermore, real fractured media are characterized by a huge computational domain, that involves a large number of unknowns, and often present an intricate system of intersections at different scales, e.g., large faults are crossed by small fractures.
Smart approaches to tackle DFN simulations involve the use of non-conforming meshes at fracture interfaces, combined with constrained optimization problems or high-performance domain decomposition strategies \cite{BPST, berrone2021three, BURMAN2019726, KMJR, BSV, BERRONE2019904}.
Alternative methods rely on a global conforming mesh on the whole network: this requires a numerical method able to deal with generic polygonal tessellations, as the generation of a standard FEM mesh over a complex domain would be computationally demanding.
For this reason, polygonal methods such as the Hybrid High Order \cite{chave2018hybrid, PFA}, the Mimetic Finite Difference \cite{BeiraodaVeiga-Lipnikov-Manzini:2014,Lipnikov-Manzini-Shashkov:2014,antonietti2016mimetic} or the Virtual Element Method (VEM) \cite{fumagalli2018dual, BENEDETTO2014135, benedetto2016globally, BERRONE2021103502} have recently been investigated.
This work focuses on the second mesh approach combined with the VEM.

The generation of a polygonal tessellation of a DFN conforming to the fracture interfaces is not a trivial task, due to the large number of fractures typically found in real applications, the different sizes of neighboring domains, and the high complexity of their intersections.
In particular, it frequently leads to a high number of degrees of freedom (DOFs) in the discrete problem, and to the generation of numerous elongated elements with aligned small edges.
Such ``bad-shaped'' cells impact on the accuracy of the simulation, producing high condition numbers and numerical errors, while the huge number of DOFs translates into a significant computational cost.

For these reasons, we propose the application of the Quality Agglomeration algorithm locally on each fracture of the network, which is able to reduce the total unknowns of the discrete problem with a control on the geometric quality of the elements.
To evaluate the effects of the Quality Agglomeration, we run the algorithm over two DFNs of increasing complexity and solve a numerical problem over the original and the agglomerated versions.
Then, we compare the results produced by the VEM over the different networks in terms of computational cost, approximation error, and convergence rate.

The paper is organized as follows.
We devote Section~\ref{sec:agglomeration} to the description of the mesh quality agglomeration algorithm.
In Section~\ref{sec:model} we introduce the DFN model, the numerical problem to be solved and the VEM discretization used.
In Section~\ref{sec:results} we report the numerical tests performed on two DFNs of increasing complexity.
Finally, in Section~\ref{sec:conclusions} we collect conclusions and possible future research directions.

\section{Mesh Quality Agglomeration}
\label{sec:agglomeration}
Consider the two-dimensional domain $\Omega$, tessellated with an initial polygonal mesh $\Mesh$.
The Quality Agglomeration algorithm aims at reducing the total number of DOFs in the virtual element approximation, keeping control over the mesh geometric quality.
By Quality Agglomeration, we mean the process in which we glue some mesh elements to form bigger elements.
The gluing strategy is driven by the minimization of a functional related to a notion of mesh “quality” suited explicitly to the VEM.
Furthermore, we admit a finite number of constraints, so that the gluing process preserve somes mesh nodes and edges.
In particular, we are interested in operating over the two following types of elements:
\begin{enumerate}
  \item Elements located in areas of the domain that are poor in details or features. In such a case, the algorithm can increase the size of the elements without a significant impact on the accuracy of the simulation;
  \item ``Badly-shaped'' or pathological elements according to the chosen notion of quality. The gluing algorithm can often improve the quality of such elements by merging them with neighboring elements.
\end{enumerate}

Introducing the notation, we consider the mesh $\Mesh$ formed by polygonal elements $\P$ with diameter $\Diameter_{\MeshCellTwoD}$, centroid $(x_{\MeshCellTwoD}, y_{\MeshCellTwoD})$ and area $\seminorm{\MeshCellTwoD}$.
The boundary $\Boundary{\MeshCellTwoD}$ of the elements is subdivided into edges $\MeshCellOneD$ with length $\seminorm{\MeshCellOneD}$.
The positive number $\MeshParameter$ indicates the mesh size, i.e., the finite constant that bounds all the element diameters $\Diameter_{\MeshCellTwoD}$.

\subsection{Mesh quality}
\label{subsec:indicator}
The first crucial component of the Quality Agglomeration algorithm is the notion of \textit{mesh quality}.
The standard approach in the VEM literature to ensure a good quality mesh is to impose some \textit{geometrical} (or \textit{regularity}) \textit{assumptions} that a mesh must respect to guarantee the optimal behavior of the method.
In \cite{sorgente2021role}, the authors isolate the four principal assumptions typically required (even if not all simultaneously) for the convergence of the VEM on a given mesh family: 
\begin{itemize}
    \item[] \textbf{G1}:\quad there exists a real number $\rho\in(0,1)$, independent of $h$, such that every polygon $\P\in\Th$ is star-shaped with respect to a disc with radius $\rP\ge\rho\hP$;
    \item[] \textbf{G2}:\quad there exists a real number $\rho\in(0,1)$, independent of $h$, such that for every polygon $\P\in\Th$, each edge $\MeshCellOneD \in \Boundary{\MeshCellTwoD}$ satisfies $\seminorm{\MeshCellOneD}\ge\rho\hP$;
    \item[] \textbf{G3}:\quad there exists a positive integer $\Ns$, independent of $\hh$, such that the number of edges of every polygon $\P\in\Th$ is (uniformly) bounded by $\Ns$;
    \item[] \textbf{G4}:\quad for every polygon $\P\in\Th$, the $1$-dimensional mesh $\calI_\P$ induced by $\Boundary{\MeshCellTwoD}$ is piecewise quasi-uniform.
\end{itemize}
Using these assumptions as absolute conditions that a mesh can only satisfy or violate has been shown to be particularly restrictive \cite{sorgente2021role}.
Instead, we define the quality of a mesh as a measure of \textit{how much} it satisfies the above conditions.
This approach is more accurate as it captures small quality differences between meshes and does not exclude a priori all the particular cases of meshes that are only slightly outside the geometrical assumptions.

In \cite{sorgente2021role}, the authors derive four scalar functions $\varrho_s:\Th\to[0,1]$, which measure how a mesh element $\P\in\Th$ meets the requirements of assumption \textbf{Gs}, for $\textbf{s}=1,\ldots,4$.
In particular, $\varrho_s=0$ if $\P$ does not respect \textbf{Gs}, and the higher $\varrho_s$ the better $\P$ is shaped with respect to \textbf{Gs}.
We briefly report the indicator definitions and refer the reader to the original paper and to \cite{sorgente2021vem} for a more complete discussion:
\begin{align*}
\varrho_1(\P) &= \frac{\lvert kernel(\P)\rvert}{\mP} = 
    \begin{cases} 
    1 & \mbox{if $\P$ is convex} \\
    \in (0,1) & \mbox{if $\P$ is concave and star-shaped} \\
    0 & \mbox{if $\P$ is not star-shaped} \\
    \end{cases}\\
\varrho_2(\P) &= \frac{\min(\sqrt{\mP},\ \min_{\E\in\partial\P}\mE)} {\hP} \\
\varrho_3(\P) &= \frac{3}{\#\left\{\E\in\partial\P \right\}} = 
    \begin{cases} 
    1 & \mbox{if $\P$ is a triangle} \\
    \in (0,1) & \mbox{otherwise} \\
    \end{cases}\\
\varrho_4(\P) &= \min_j \frac{\min_{\E\in\calI_\P^j}\mE}{\max_{\E\in\calI_\P^j}\mE}
\end{align*}
The $kernel$ operator computes the kernel of a polygon, intended as the set of points from which the whole polygon is visible, and $\calI_\P^j$ are all the 1-dimensional disjoint sub-meshes corresponding to the edges of $\P$ (we consider each $\calI_\P^j$ as a mesh as it may contain more than one edge, see \cite{beirao2022sharper}) such that $\calI_\P=\cup_j\calI_\P^j$, where $\calI_\P$ is the 1-dimensional mesh induced by $\partial\P$ introduced in assumption \textbf{G4}.

We combine together the four indicators $\varrho_1, \varrho_2$, $\varrho_3$ and $\varrho_4$ into a \textit{quality indicator} $\varrho:\Th\to [0,1]$, which measures the overall quality of an element $\P\in\Th$:
\begin{equation}
\varrho(\P) := \sqrt{\frac{\varrho_1(\P)\varrho_2(\P) + \varrho_1(\P)\varrho_3(\P)+ \varrho_1(\P)\varrho_4(\P)}{3}}.
\label{eq:indicator}
\end{equation}
We have $\varrho(\P)\to1$ if $\P$ is a perfectly-shaped element, e.g. an equilateral triangle or a square, $\varrho(\P)=0$ if and only if $\P$ is not star-shaped, and $0<\varrho(\P)<1$ otherwise.
We remark that all indicators $\{\varrho_s\}_{s \in \{1,\ldots,4\}}$, and consequently $\varrho$, only depend on the geometrical properties of the mesh elements; therefore their values can be computed before applying the VEM, or any other numerical scheme.

We point out that the local quality indicator \eqref{eq:indicator} is the restriction to a single element of the mesh quality indicator introduced in \cite{sorgente2021role}.
A strict correspondence has been proven between the values measured by the mesh quality indicator and the performance of the VEM in terms of approximation errors and convergence rates.
In particular, the VEM is likely to produce small approximation errors over meshes with a high quality value in the sense of \eqref{eq:indicator}.
Moreover, given a collection of refined meshes with decreasing mesh size, the VEM is expected to converge rapidly at the optimal rate if the quality of the meshes remains constant throughout the refinement process (note that $\varrho$ is scale-independent).

\subsection{Mesh agglomeration}
\label{subsec:graphcut}
The second component of the Quality Agglomeration algorithm is the procedure that agglomerates groups of elements in the mesh, following the information provided by the quality indicator.
We perform the element agglomeration by solving an optimization problem that balances the mesh elements' quality and number. 
The energy functional $\calE:\Mesh\to\RSet$ that we want to minimize is:
\begin{equation}
\label{eq:graphcut:energy}
    \calE:=\sum_{\P,\P'\in\Mesh} dc(\P,\P') + \lambda\sum_{\P,\P'\in\Mesh} sc(\P,\P').
\end{equation}
The factor $\lambda\in[0,1]$ of \eqref{eq:graphcut:energy} is the \textit{agglomeration parameter} that balances the importance of the \textit{cost functions} $dc$ and $sc$.
Empirically, high values of $\lambda$ lead to more aggressive agglomerations.
We define the \textit{cost functions} $dc$ and $sc$ from the product space $\Mesh\times\Mesh$ to the unit interval $[0,1]$ as follows:
\begin{itemize}
    \item the \textit{data cost} ($dc$) represents the cost of agglomerating two elements $\P,\P'\in\Mesh$, and measures the potential quality of the element $\P\cup\P'$.
    We define it as:
    \begin{align}
    \label{eq:graphcut:dc}
        dc(\P,\P'):=\begin{cases}
        0 &\mbox{ if $\P=\P'$}\\
        1-\varrho(\P\cup\P') &\mbox{ if $\P$ and $\P'$ are adjacent};\\
        1 &\mbox{ otherwise}.
        \end{cases}
    \end{align}
    Here, $\P\cup\P'$ is the boolean union of the two neighboring elements, and $\varrho$ is the elemental quality indicator \eqref{eq:indicator} of the union;
    \item the \textit{smoothness cost} ($sc$) encodes information on the structure of the mesh, setting a zero weight to the non-adjacent elements:
    \begin{align}
    \label{eq:graphcut:sc}
        sc(\P,\P'):=\begin{cases}
        1 &\mbox{ if $\P$ and $\P'$ are adjacent and $\P\neq\P'$},\\
        0 &\mbox{ otherwise}.
        \end{cases}
    \end{align}
\end{itemize}

To minimize the energy functional defined in \eqref{eq:graphcut:energy}, we interpret the mesh $\Mesh$ as a graph in which a node represents an element $\P\in\Mesh$. The graph connects two nodes if their corresponding elements are adjacent, i.e., if they share at least one mesh edge.
Moreover, we label each node with an integer number that provides information about the node.
Let $L$ denote the set of all possible labels and $\calL:\Mesh\to L$ the map that assigns a label $l\in L$ to each node $\P\in\Mesh$.
For implementation reasons, we multiply $dc$ and $sc$ by the number of elements in $\Mesh$ and round them off to obtain an integer value.

We use the graph-cut technique \cite{boykov2001fast,kolmogorov2004energy,boykov2004experimental} to solve the energy optimization problem, in the implementation provided by the \textit{Multi-label Optimization} method \cite{gco}.
Graph-cut iterates and opportunely re-labels the nodes until it satisfies a given convergence criterion.
At the end of the iterations, the algorithm agglomerates the nodes sharing the same label.
In this technique, the cost functions are defined in the form $\widetilde{dc}:\Mesh\times L\to [0,1]$ and $\widetilde{sc}:L\times L\to [0,1]$.
To align to this notation, we set $L:=\{1,\ldots,\#\Mesh\}$, and define the trivial labeling $\widetilde{\calL}:\Mesh\to L$ that bijectively maps each element of the mesh to its \textit{index}, i.e., its position in the array of the mesh data structure containing all the elements.
Then, we opportunely compose the cost functions \eqref{eq:graphcut:dc} and \eqref{eq:graphcut:sc} with the inverse map $\widetilde{\calL}^{-1}:L\to\Mesh$, which, therefore, connects a label $l\in L$ to the (unique) element $\P\in\Mesh$ with index $l$:
\begin{align*}
\widetilde{dc}(\P,l) &:= dc(\P,\widetilde{\calL}^{-1}(l)) = dc(\P,\P'),\\
\widetilde{sc}(l_1,l_2) &:= sc(\widetilde{\calL}^{-1}(l_1),\widetilde{\calL}^{-1}(l_2)) = sc(\P_1,\P_2),
\end{align*}
where $\P'$, $\P_1$, and $\P_2$ are the graph nodes whose corresponding elements have indices $l$, $l_1$, and $l_2$, respectively.
Substituting $\widetilde{dc}$ and $\widetilde{sc}$ into \eqref{eq:graphcut:energy}, we obtain a new energy functional, which depends on the particular labeling:
\begin{align}
    \widetilde{\calE}(\calL):=\sum_{\P\in\Mesh}\widetilde{dc}(\P,l_\P) + \lambda\sum_{\P,\P'\in\Mesh} \widetilde{sc}(l_\P,l_{\P'}).
    \label{eq:energy2}
\end{align}
The minimization problem reads as $\min_{\calL\in \calP}\widetilde{\calE}(\calL)$, where $\calP$ is the set of all the possible labelings.
We solve this problem by using the \textit{alpha-beta swap} algorithm \cite{boykov2001fast}.
This technique iteratively segments the graph nodes labeled with a given label $\alpha$ to those with another label $\beta$.
These two labels change after each iteration, scouting all the possible combinations.
Other algorithms exist in the literature, for example, the so-called \textit{alpha-expansion} algorithm, which requires the function $sc$ to be a metric (i.e., to respect the triangular inequality).
The results obtained with this latter algorithm are less useful in our context because it generally leads to uneven label distributions.
Indeed, the alpha-expansion algorithm tries to expand each label as much as possible, thus assigning the same label to large part of the domain and different labels only to small isolated areas.
This strategy leads to meshes with unbalanced elements.

\subsection{The algorithm}
\label{subsec:algorithm}
\begin{figure}[htbp]
    \centering
    \subfloat[\label{fig:graphcut_1}]{\includegraphics[width=0.23\textwidth]{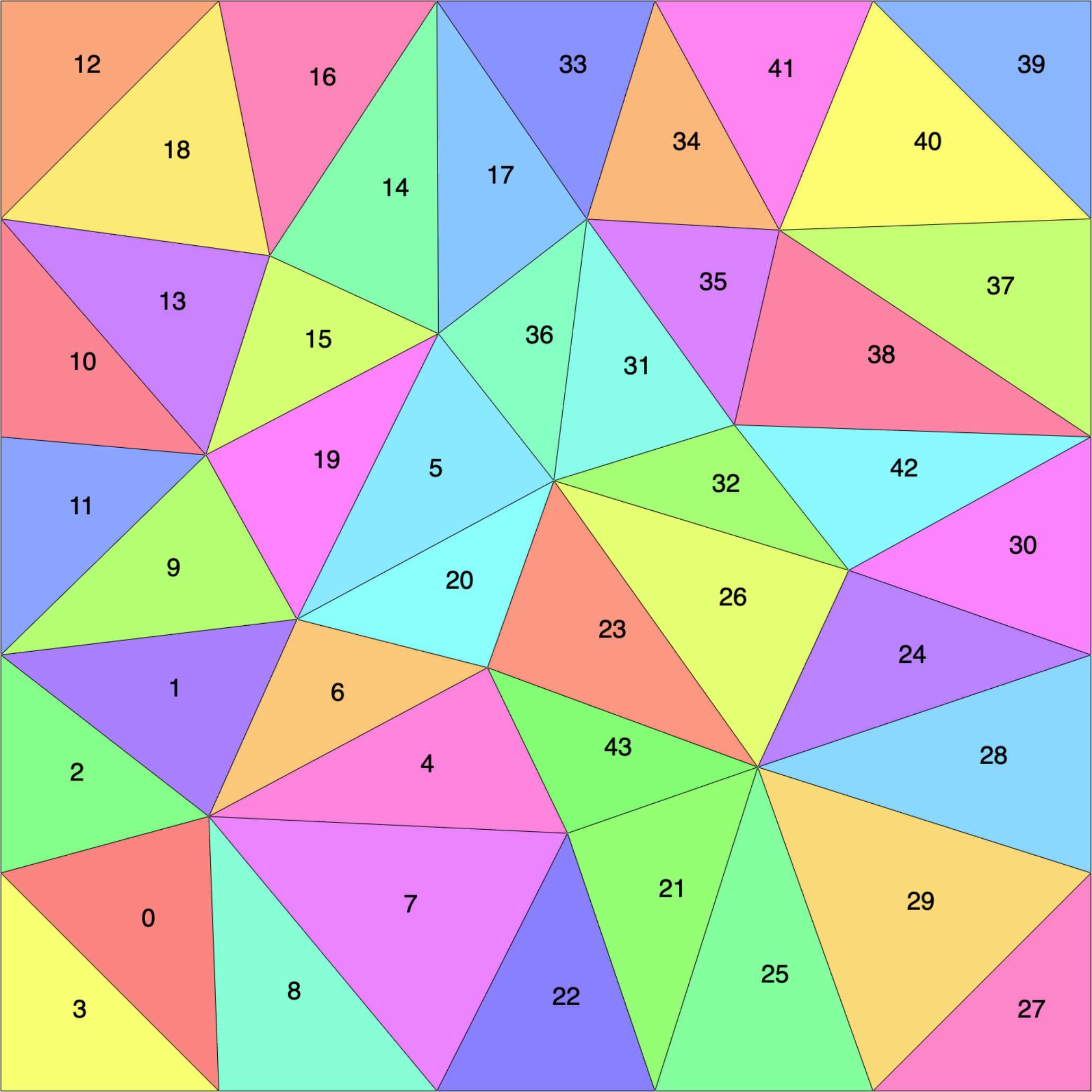}}\
    \subfloat[\label{fig:graphcut_2}]{\includegraphics[width=0.23\textwidth]{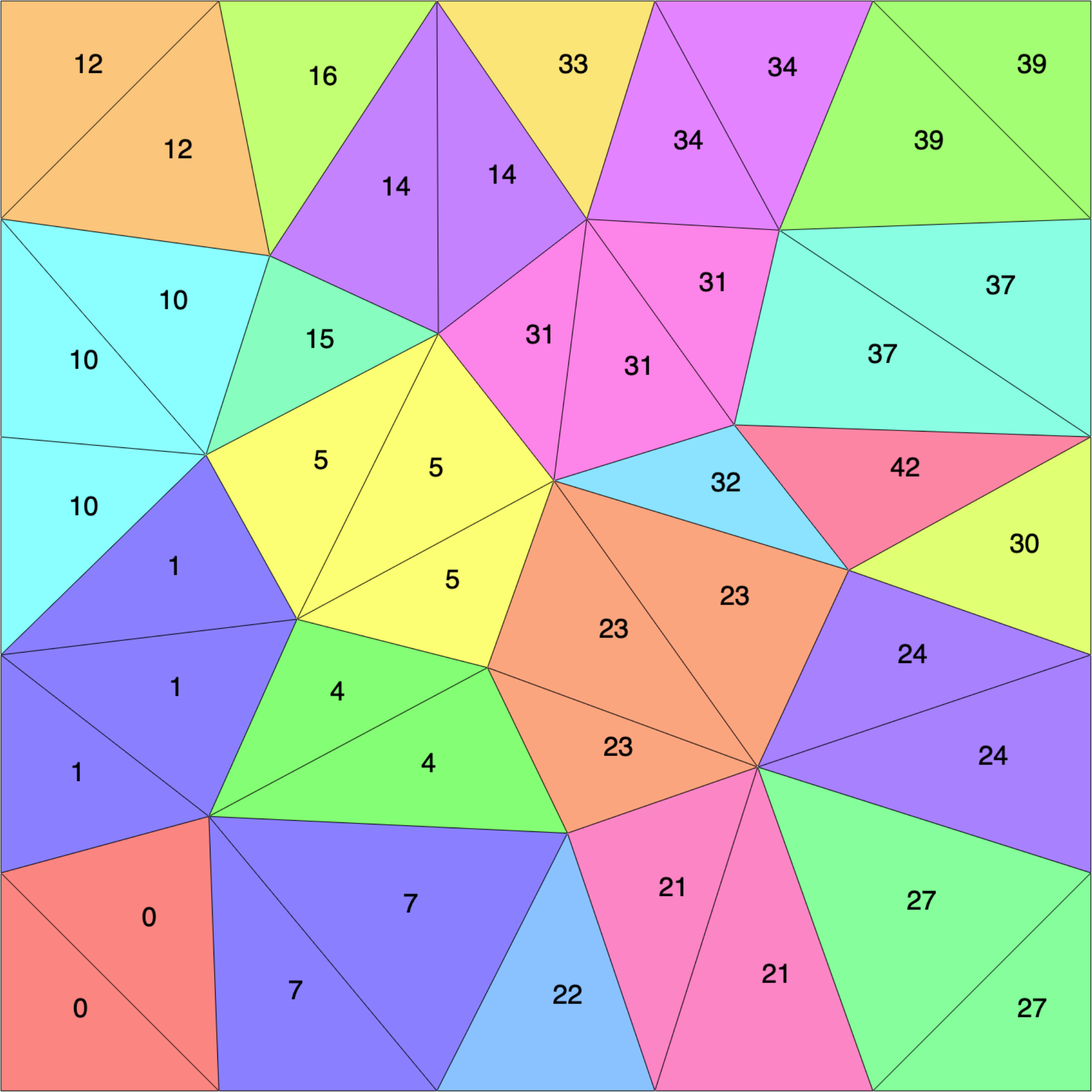}}\
    \subfloat[\label{fig:graphcut_3}]{\includegraphics[width=0.23\textwidth]{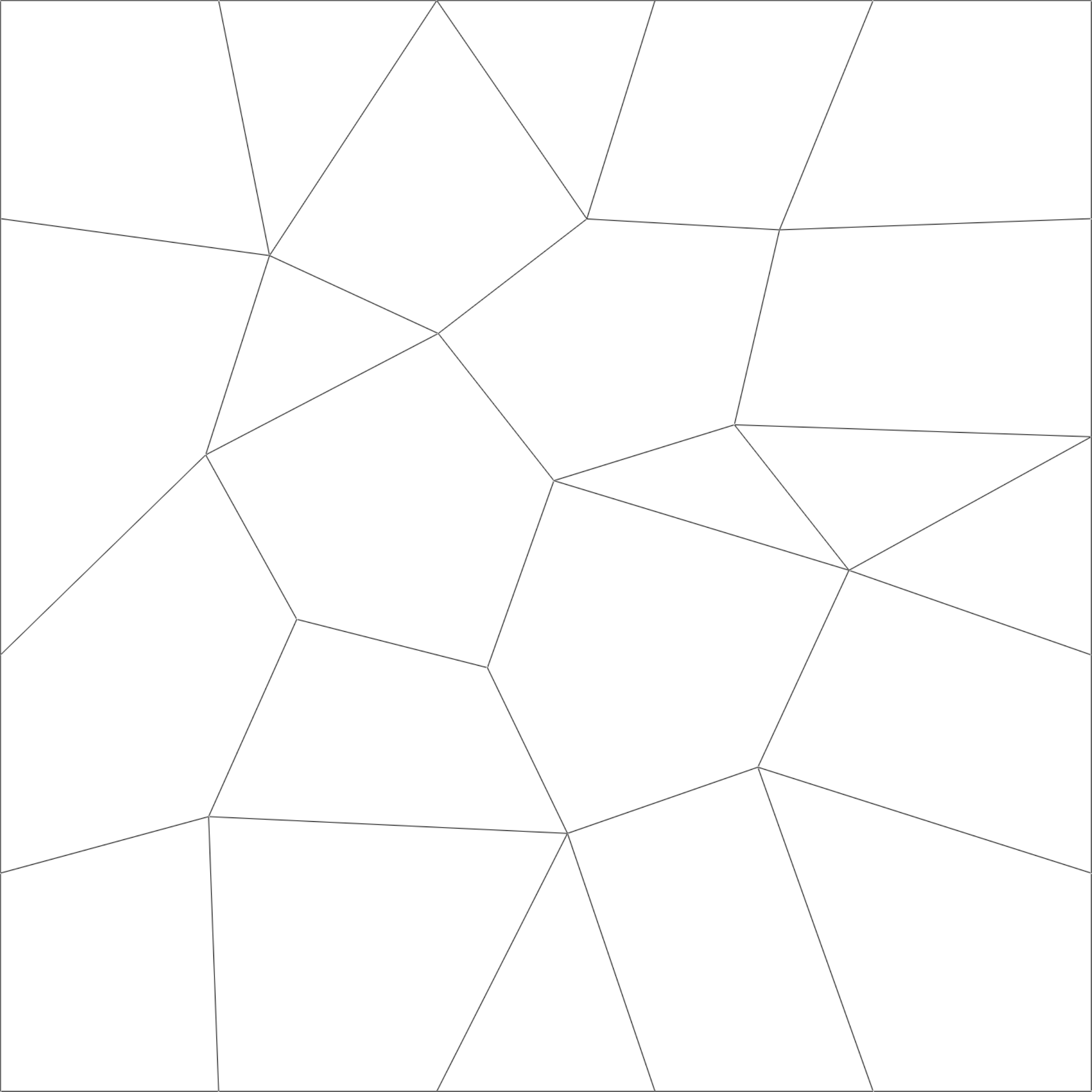}}\
    \subfloat[\label{fig:graphcut_4}]{\includegraphics[width=0.23\textwidth]{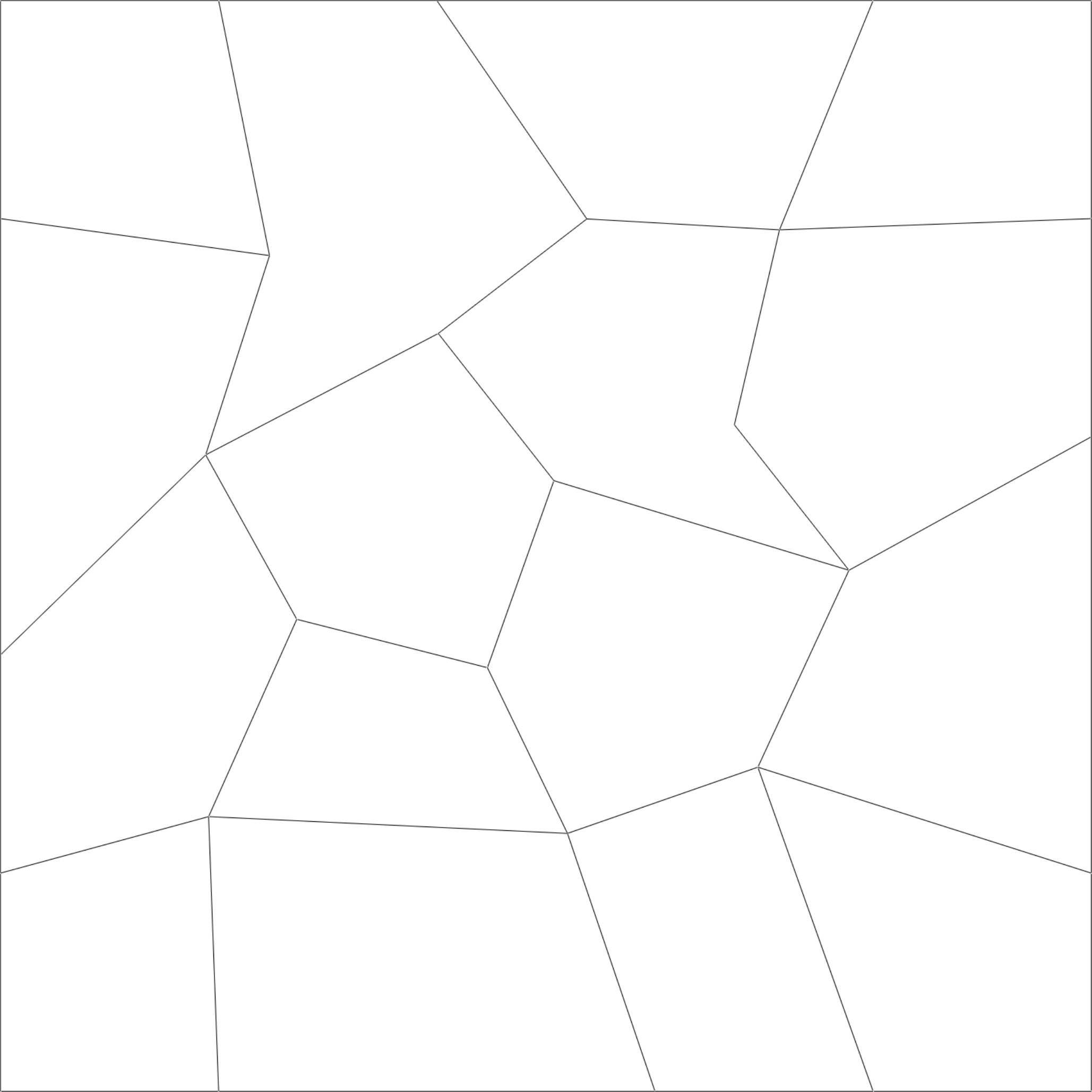}}
    \caption{(a) Initial mesh $\OriginalMesh$, colored w.r.t. element labels (in black), (b) mesh after graph-cut with $\lambda=0.25$, colored w.r.t. element labels (in black), (c) agglomerated mesh $\Mesh$ with $\lambda=0.25$, (d) agglomerated mesh $\Mesh$ with $\lambda=1.0$.}
    \label{fig:graphcut}
\end{figure}

We now detail the Quality Agglomeration algorithm for the VEM.
Given an input mesh $\OriginalMesh$, we build the corresponding graph by generating a node for each element in the mesh and connecting the nodes corresponding to adjacent elements.
In the case of a constrained mesh edge that has to be preserved, we consider the elements sharing that edge as non-adjacent.
We initialize the node labels with the trivial labeling $\widetilde{\calL}$, see Figure~\ref{fig:graphcut_1}.
Then, we run the graph-cut technique with a selected value of the parameter $\lambda$, and find new labels for the elements (Figure~\ref{fig:graphcut_2}).
Note that graph-cut only operates over the labels without actually modifying the mesh.
Graph-cut reaches convergence when the energy term \eqref{eq:energy2} attains a (local) minimum.
This task is typically accomplished in a few iterations.
Last, a post-processing step is required to merge all elements sharing the same label (Figure~\ref{fig:graphcut_3}) and create the agglomerated mesh $\Mesh$.
To further improve the mesh quality, the algorithm merges aligned edges in the newly generated elements, while preserving possible constraints on nodes or edges of the initial mesh $\OriginalMesh$.

Since the agglomerated mesh $\Mesh$ is the solution to a minimization problem, it is also optimal in the number of elements.
However, we can compute different optimal versions of the same mesh $\OriginalMesh$ by choosing different values of the parameter $\lambda$ in \eqref{eq:energy2}, as shown in Figures~\ref{fig:graphcut_3}-\ref{fig:graphcut_4}.
Higher values of $\lambda$ provide agglomerated meshes containing a smaller number of elements, edges, nodes, and therefore degrees of freedom of the VEM space, thus requiring a smaller computational cost for the VEM.
At the same time, we expect the errors produced by the agglomerated meshes to be slightly higher, as for every removed DOF we have less information on the numerical problem.
We are interested in understanding how this reduction of degrees of freedom impacts the accuracy of the VEM, and if it affects the convergence rate of the method.

\section{DFN Problem formulation}
\label{sec:model}

To test the Quality Agglomeration algorithm, we introduce the steady-state flow problem in DFNs as a possible application.
We consider a network domain $\Network \subset \RSet[3]$ made of a finite number $\bar{i}$ of fractures $\Fracture_i$, $i\in\FractureIndices:=\{1,\ldots,\bar{i}\}$; thus, $\Network := \bigcup_{i \in \FractureIndices} \Fracture_i$, such that each $\Fracture_i$ has at least one intersection with at least one $\Fracture_j$, for $j \in \FractureIndices \setminus \{i\}$.

In the DFN model, each $\Fracture_i$ is represented by a
two-dimensional polygonal tessellation oriented in $\RSet[3]$ that
approximates a geological fracture immersed in the impervious rock
matrix, as done in
\cite{AHMED201749, benedetto2016globally, BERRONE2021103502, martin2005modeling, DF}.
With this assumption, the network results in a collection of two-dimensional elements.
Therefore, we refer to it with the same symbol $\Network$ of Section~\ref{sec:agglomeration} even if it is not properly contained in $\RSet[2]$. 

Let $\Boundary{\Network} := \bigcup_{i \in \FractureIndices} \Boundary{\Fracture_i}$ be the domain boundary.
We split it into a Dirichlet part $\Gamma_D$ and a Neumann part $\Gamma_N$, such that $\Boundary{\Network} = \Gamma_D \cup \Gamma_N$, $\Gamma_D \cap \Gamma_N = \emptyset$, and $\seminorm{\Gamma_{D}} \neq 0$.

The fractures intersect along a finite number $\bar{m}$ of segments, denoted as $\Trace_m$, with $m\in\TraceIndices:=\{1,\ldots,\bar{m}\}$.
For the sake of simplicity, we assume that each intersection is formed by exactly two fractures, i.e., $\Trace_m := \Fracture_i \cap \Fracture_j$, fixing a unique pair of fracture indices $\TracePair_m = \{i, j\}$ for each $m \in \TraceIndices$.
We point out that this is a simplifying assumption that does not alter the results.
Finally, we denote by $\TraceIndices_i$ the set of the indices of the fracture intersections $\Trace_m$ which lie on $\Fracture_i$, i.e., $\TraceIndices_i:=\{m\in\TraceIndices:\ \Trace_m \cap \Fracture_i \neq \emptyset\}$.

In what follows, we use the standard definition and notation of Sobolev spaces, norms, and seminorms, cf.~\cite{Adams-Fournier:2003}; thus, given a bounded, connected subset $\omega$ of $\RSet[2]$, we denote the dot product and the norm in $\SobL[2]{\omega}$-space with $\scalar{\mathcal{\omega}}{\cdot}{\cdot}$ and $\HilbertNorm[\mathcal{\omega}]{\,\cdot\,}$ respectively. 
In addition, $\seminorm[\omega]{\,\cdot\,}$ is the seminorm of $\SobH{1}{\omega}$ and $\nabla_i$ represents the tangential component of the tridimensional gradient operator $\nabla$ on the plane of $\Fracture_i$.

We seek the distribution of the hydraulic head $\HydraulicHead$ in the whole network $\Network$; we prescribe the Dirichlet boundary condition on $\Gamma_D$ through function $g \in \SobH{\frac{1}{2}}{\Gamma_D}$ and set the homogeneous Neumann boundary condition on $\Gamma_N$.
We define the functional spaces:
\begin{align*}
    \SolutionSpace^D := &\left\{ v : v_{\mid{\Gamma_D}} = g,\ v_{\mid\Fracture_i} \in \SobH{1}{\Fracture_i}\ \forall i \in \FractureIndices,\ v_{i\mid{\Trace_m}} = v_{j\mid{\Trace_m}}\ \forall m \in \TraceIndices \right\},\\
    \SolutionSpace := &\left\{ v : v_{\mid\Fracture_i} \in \SobHo{1}{\Fracture_i}\ \forall i \in \FractureIndices,\ v_{i\mid{\Trace_m}} = v_{j\mid{\Trace_m}}\ \forall m \in \TraceIndices \right\}.
\end{align*}
The weak formulation reads: find $\HydraulicHead \in \SolutionSpace^D$ such that $\HydraulicHead - \HydraulicHead^D \in \SolutionSpace$ with $\HydraulicHead^D \in \SolutionSpace^D$ such that $\HydraulicHead^D_{\mid\Gamma_D} = g$ and, for all $v \in \SolutionSpace$,
\begin{equation}
    \sum_{i \in \FractureIndices} \scalar{\Fracture_i}{\Diffusion_i \nabla_i \HydraulicHead_i}{\nabla_i v_{\mid\Fracture_i}} = \sum_{i \in \FractureIndices} \scalar{\Fracture_i}{f_{\mid\Fracture_i}}{v_{\mid\Fracture_i}}, \label{eq:weak}
\end{equation}
where $\HydraulicHead_i$ is the restriction of the hydraulic head to $\Fracture_i$ and $\Diffusion_i : \RSet[2] \to \RSet[2 \times 2]$ is the in-plane transmissivity on the fracture.

According to~\cite[Theorem~2.7.7]{Brenner-Scott:2008}, problem \eqref{eq:weak} is well posed since the symmetric bilinear form
\begin{equation}
    \BilinearForm{a}{v}{w} := \sum_{i \in \FractureIndices} \BilinearForm{a^{\Fracture_i}}{v}{w} = \sum_{i \in \FractureIndices} \scalar{\Fracture_i}{\Diffusion_i \nabla_i w_{\mid\Fracture_i}}{\nabla_i v_{\mid\Fracture_i}} \label{eq:bilinear}
\end{equation}
is coercive and continuous on $\SolutionSpace$.

\subsection{Network discretization}
\label{subsec:problem}

We describe the approach used to construct an optimal polygonal tessellation $\Mesh$ on $\Network$, globally conforming at the fracture intersections, exploiting the Quality Agglomeration algorithm of Section~\ref{subsec:algorithm}.
The approach is composed of three steps.

First, we independently discretize each fracture domain $\Fracture_i \in \Network$ by a classical triangular planar mesh $\mathcal{T}_{\MeshParameter}^i$ of given size $\MeshParameter$.
On each $\mathcal{T}_{\MeshParameter}^i$, given the set of local fracture segments $\Trace_m$ with $m \in \TraceIndices_i$, we split all the cells $T \in \mathcal{T}_{\MeshParameter}^i$ that intersect the segments $\Trace_m$ by the direction of the segments, \cite{benedetto2016globally}.
Thus, on each $\Fracture_i$ we obtain  a local polygonal mesh $\mathcal{P}_{\MeshParameter}^i$, which is conforming to all $\Trace_m$ with $m \in \TraceIndices_i$.
Since we do not implement any particular geometrical smoothing technique on this cutting phase, the meshes $\mathcal{P}_{\MeshParameter}^i$ may present small edges and elongated ``bad-shaped'' cells, in the sense of the geometrical assumptions from Section~\ref{subsec:indicator}.

Second, the Quality Agglomeration algorithm described in Section~\ref{subsec:algorithm} processes all the meshes $\mathcal{P}_{\MeshParameter}^i$ independently and imposes the constraints on the mesh vertices that are the segment endpoints of each fracture intersections $\Trace_m$, $m \in \TraceIndices_i$.
Thus, on each domain $\Fracture_i$ we obtain a locally optimal quality discretization $\Mesh^i$ that still respects the conformity constraints along boundaries and interfaces $\Trace_m$ with $m \in \TraceIndices_i$.

Finally, we collect the optimal meshes $\Mesh^i$ to generate the global polygonal conforming mesh $\Mesh$ on the whole network.
We perform the conforming operation working only on the intersection segments $\Trace_m$, $m \in \TraceIndices$.
For each $m \in \TraceIndices$, where $\TracePair_m = \{i, j\}$, we unify the two sets of mesh nodes in $\RSet[3]$ of $\Mesh^i$ and $\Mesh^j$ lying on the segment $\Trace_m$, and create new mesh elements $\MeshCellTwoD\in\Mesh$ that have aligned edges in correspondence of the new nodes; we recall that the VEM allows elements with such an elaborated geometry.

With this approach, we obtain that each element $\MeshCellTwoD \in \Mesh$ lies in only one $\Fracture_i \in \Network$.
Therefore, quantities introduced locally on $\MeshCellTwoD$ are considered with respect to the two-dimensional tangential reference system of $\Fracture_i$, and we omit in what follows the suffix $_i$ on all the operators (e.g., $\nabla_i$ becomes $\nabla$).

\subsection{Virtual Element approximation}
\label{subsec:vem}
We introduce $\Polynomial_{n}(\MeshCellTwoD)$ and $\Polynomial_{n}(\MeshCellOneD)$, the spaces of polynomials of degree less or equal to $n \in \NSet[+]$ defined on each $\MeshCellTwoD \in \Mesh$ and each $\MeshCellOneD \in \Boundary{\MeshCellTwoD}$, respectively.
We remark that $\Polynomial_{n}(\MeshCellTwoD)$ is built on the reference system that is tangential to the fracture $\Fracture_i$ to which the mesh element $\MeshCellTwoD$ belongs.
We define $\Projection^{\nabla}_{n, \MeshCellTwoD} : \SobH{1}{\MeshCellTwoD} \to \Polynomial_{n}(\MeshCellTwoD)$ as the $\SobH{1}{\MeshCellTwoD}$-orthogonal projection operator, computed for any $p \in \Polynomial_{n}(\MeshCellTwoD)$ and $v \in \SobH{1}{\MeshCellTwoD}$ with the conditions:
\begin{equation*}
    \begin{cases}
        \scalar{\MeshCellTwoD}{\nabla \Projection^{\nabla}_{n, \MeshCellTwoD} v}{\nabla p} = \scalar{\MeshCellTwoD}{\nabla v}{\nabla p} &                   \\
        \scalar{\Boundary{\MeshCellTwoD}}{\Projection^{\nabla}_{n, \MeshCellTwoD} v}{1} = \scalar{\Boundary{\MeshCellTwoD}}{v}{1}        & \text{if } n = 1 \\
        \scalar{\MeshCellTwoD}{\Projection^{\nabla}_{n, \MeshCellTwoD} v}{1} = \scalar{\MeshCellTwoD}{v}{1}                              & \text{otherwise}
    \end{cases}
\end{equation*}
Similarly, we let $\Projection^{0}_{n, \MeshCellTwoD}$ and $\Projection^{0}_{n, \MeshCellTwoD} \nabla$ denote the $\SobL{E}$-orthogonal projection operators on $\Polynomial_{n}(\MeshCellTwoD)$ for functions $v$ in $\SobH{1}{\MeshCellTwoD}$ and on $\Polynomial_{n}(\MeshCellTwoD) \times \Polynomial_{n}(\MeshCellTwoD)$, respectively.

We locally define the Virtual Element space of order $\VemOrder\geq1$ on $\MeshCellTwoD$ as
\begin{eqnarray*}
    \SolutionSpace^{\VemOrder, \MeshCellTwoD}_{\MeshParameter} := \Bigl\{ v \in \SobH{1}{\MeshCellTwoD} : & \Delta v \in \Polynomial_{\VemOrder}(\MeshCellTwoD),\ v_{\mid \Boundary{\MeshCellTwoD}} \in C^{0}(\Boundary{\MeshCellTwoD}),\ v_{\mid\MeshCellOneD} \in \Polynomial_{\VemOrder}(\MeshCellOneD)\ \forall \MeshCellOneD \in \MeshCellTwoD, \\
    & \scalar{\MeshCellTwoD}{v}{p} = \scalar{\MeshCellTwoD}{\Projection^{\nabla}_{\VemOrder, \MeshCellTwoD} v}{p} \ \forall p \in \Polynomial_{\VemOrder}(\MeshCellTwoD) \setminus \Polynomial_{\VemOrder - 2}(\MeshCellTwoD) \Bigl\}.
\end{eqnarray*}
Given a function $v \in \SolutionSpace^{\VemOrder, \MeshCellTwoD}_{\MeshParameter}$, we select the following standard degrees of freedom on the element $\MeshCellTwoD$:
\begin{itemize}
    \item the values of $v$ in the vertices of $\MeshCellTwoD$;
    \item if $\VemOrder > 1$, the values of $v$ on the $\VemOrder - 1$ internal Gauss-Lobatto quadrature points on every edge $\MeshCellOneD \in \MeshCellTwoD$;
    \item if $\VemOrder > 1$, the $\frac{\VemOrder (\VemOrder - 2)}{2}$ moments of $v$ on $E$: $\frac{1}{\seminorm{\MeshCellTwoD}}\scalar{\MeshCellTwoD}{v}{m_j}$, $\forall m_j \in \Basis_{\VemOrder-2}(\MeshCellTwoD)$,
\end{itemize}
where $\Basis_{n}(\MeshCellTwoD)$ is the selected basis for the space $\Polynomial_{n}(\MeshCellTwoD)$.
The chosen DOFs are unisolvent for $\SolutionSpace^{\VemOrder, \MeshCellTwoD}_{\MeshParameter}$, and all the projection operators $\Projection^{\nabla}_{\VemOrder, \MeshCellTwoD}$, $\Projection^{0}_{\VemOrder, \MeshCellTwoD}$, and $\Projection^{0}_{\VemOrder - 1, \MeshCellTwoD}\nabla$ are computable \cite{BPVEM}.
Similarly to \cite{CMSO}, we use the set of two-dimensional scaled monomials
\begin{equation}
    \Basis_{n}(\MeshCellTwoD) := \left\{ m \in\Polynomial_{n}(\MeshCellTwoD) : m(x, y) = \frac{(x-x_{\MeshCellTwoD})^{\alpha_x}(y-y_{\MeshCellTwoD})^{\alpha_y}}{\Diameter_{\MeshCellTwoD}^{(\alpha_x+\alpha_y)}}, \ 0 \leq \alpha_x+\alpha_y \leq n\right\} \label{eq:monomial}
\end{equation}
as a basis of $\Polynomial_{n}(\MeshCellTwoD)$.
Finally, using the Lagrangian basis functions with respect to the DOFs as a basis for $\SolutionSpace^{\VemOrder, \MeshCellTwoD}_{\MeshParameter}$, we define the global discrete space $\SolutionSpace^{\VemOrder}_{\MeshParameter}$ on the whole domain as
\begin{equation*}
    \SolutionSpace^{\VemOrder}_{\MeshParameter} := \Big\{v \in \SolutionSpace : v_{\mid \MeshCellTwoD} \in \SolutionSpace^{\VemOrder, \MeshCellTwoD}_{\MeshParameter} \ \forall \MeshCellTwoD \in \Mesh \Big\}.
\end{equation*}
We can now define the discrete counterpart of the bilinear form introduced in \eqref{eq:bilinear} as $a_{\MeshParameter} : \SolutionSpace^{\VemOrder}_{\MeshParameter} \times \SolutionSpace^{\VemOrder}_{\MeshParameter} \to \RSet$, such that for all $v, w \in \SolutionSpace^{\VemOrder}_{\MeshParameter}$:
\begin{equation*}
    \BilinearForm{a_{\MeshParameter}}{v}{w} := \sum_{\MeshCellTwoD \in \Mesh} \BilinearForm{a_{\MeshParameter}^{\MeshCellTwoD}}{v}{w},
\end{equation*}
with
\begin{equation*}
    \BilinearForm{a^{\MeshCellTwoD}_{\MeshParameter}}{v}{w} := \scalar{\Fracture_i}{\Diffusion_i \Projection^{0}_{\VemOrder - 1, \MeshCellTwoD} \nabla w_{\mid \Fracture_i}}{\Projection^{0}_{\VemOrder - 1, \MeshCellTwoD} \nabla v_{\mid\Fracture_i}} + S^{\Fracture_i}_{\MeshCellTwoD}(v-\Projection^{\nabla}_{\VemOrder, \MeshCellTwoD}v, w - \Projection^{\nabla}_{\VemOrder, \MeshCellTwoD}w),
\end{equation*}
where $\Fracture_i$ is the fracture to which the mesh element $\MeshCellTwoD$ belongs.

The bilinear form $S^{\Fracture_i}_{\MeshCellTwoD}$ is selected such that
\begin{equation*}
    \exists\ \alpha, \beta > 0 :\ \alpha\ \BilinearForm{a^{\Fracture_i}_{\MeshCellTwoD}}{u}{u} \leq \BilinearForm{a_{\MeshParameter}^{\MeshCellTwoD}}{u}{u} \leq \beta\ \BilinearForm{a^{\Fracture_i}_{\MeshCellTwoD}}{u}{u} \ \forall u \in \SolutionSpace^{\VemOrder, \MeshCellTwoD}_{\MeshParameter},
\end{equation*}
where $a^{\Fracture_i}_{\MeshCellTwoD}$ is the restriction to the mesh element $\MeshCellTwoD$ of the bilinear form $a^{\Fracture_i}$ introduced in \eqref{eq:bilinear}.
Different choices for $S^{\Fracture_i}_{\MeshCellTwoD}$ are used in the literature, see for example \cite{BLR}.
In this work we select the typical form
\begin{equation}
    S^{\Fracture_i}_{\MeshCellTwoD}(v,w) := \HilbertNorm[{\SobL[\infty]{\MeshCellTwoD}}]{\Diffusion_i} \sum_{\ell = 1}^{\# \SolutionSpace^{\VemOrder, \MeshCellTwoD}_{\MeshParameter}} \text{dof}^{\MeshCellTwoD}_{\ell} (v)\ \text{dof}^{\MeshCellTwoD}_{\ell} (w),
    \label{eq:stabilization}
\end{equation}
where $\text{dof}^{\MeshCellTwoD}_{\ell} (\cdot)$ is the operator that selects the $\ell$-th degree of freedom of $\SolutionSpace^{\VemOrder, \MeshCellTwoD}_{\MeshParameter}$.

The discrete formulation of problem \eqref{eq:weak} reads as: Find $\HydraulicHead_{\MeshParameter} \in \SolutionSpace^{\VemOrder}_{\MeshParameter}$ such that, for all $v_{\MeshParameter} \in \SolutionSpace^{\VemOrder}_{\MeshParameter}$, it holds
\begin{equation}
    \BilinearForm{a_{\MeshParameter}}{\HydraulicHead_{\MeshParameter}}{v_{\MeshParameter}} = \sum_{\MeshCellTwoD \in \Mesh} \scalar{\MeshCellTwoD}{f}{\Projection^{0}_{\VemOrder - 1, \MeshCellTwoD} v_{\MeshParameter}}. \label{eq:discrete}
\end{equation}
For a detailed discussion about the approximation of the source term we refer to \cite{BPVEM}.
Finally, we can prove the well-posedness of problem \eqref{eq:discrete} and the optimal convergence rates measured in the numerical results of Section~\ref{sec:results} by using standard arguments for Virtual Element methods, see \cite{VEMGSO}.

\section{Numerical results}
\label{sec:results}
We present some numerical results for two increasingly complex networks: the first deals with a known hydraulic head in a simple setting, whereas the second simulates a quite realistic DFN configuration.

In both cases, we apply the Quality Agglomeration algorithm proposed in Section~\ref{sec:agglomeration} on each $\Fracture_i$ with $i \in \FractureIndices$, taking the parameter $\AgglomerationParameter \in [0, 1]$ constant for all the fractures of the network $\Network$.
We associate the value $\lambda = 0.0$ to the non-agglomerated strategy, and compare two different $\AgglomerationParameter$ values, namely $0.25$ and $1.0$, to test a moderate and an aggressive agglomeration strategy, respectively.
We stress that the mesh optimization approach is performed independently on each domain $\Fracture_i$.
Thus, we can tackle large networks coming from real applications with a naturally distributed parallel strategy.

We compare the values of the energy functional \eqref{eq:energy2} measured in the global polygonal conforming mesh $\Mesh$ of the network before and after the optimization process.
We test the effectiveness of the algorithm by comparing the performance of the VEM from Section~\ref{subsec:vem} with $\VemOrder = \{1, 2, 3\}$ for each $\AgglomerationParameter$ value.
We set a constant transmissivity $\Diffusion_i = \Matrixize{I}$ for all $i \in \FractureIndices$, and solve the linear system generated by problem \eqref{eq:discrete} directly, exploiting the Cholesky factorization of the C++ Eigen library, see \cite{eigenweb}.

\subsection{Network 1 - Simple DFN with known solution}
\label{sec:results:1}
\begin{figure}[htbp!]
    \centering
    \subfloat[\label{plot:RUN_T5_DMN_1}]{\includegraphics[width=0.4\textwidth]{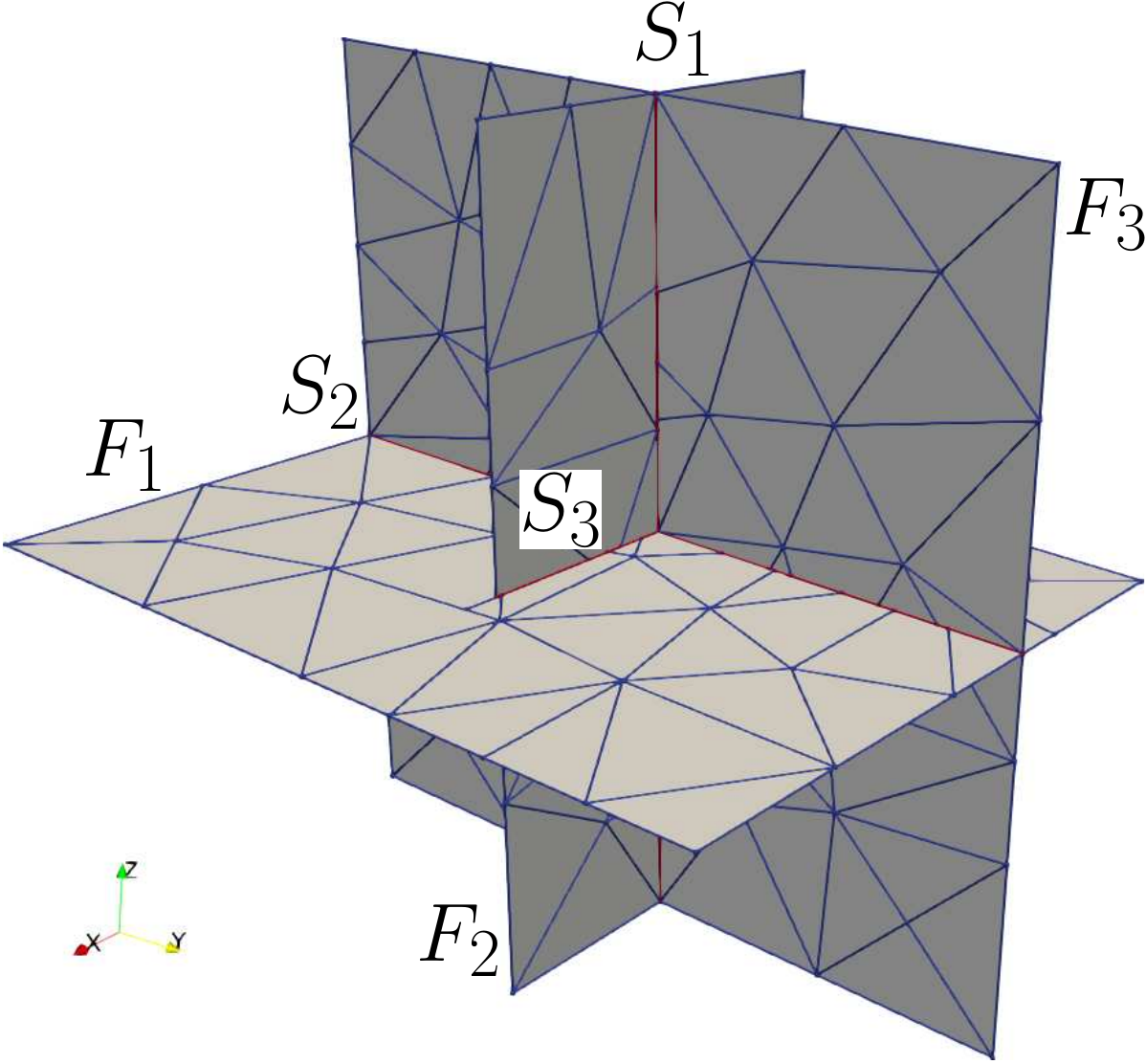}}
    \subfloat[\label{plot:RUN_T5_DMN_2}]{\includegraphics[width=0.4\textwidth]{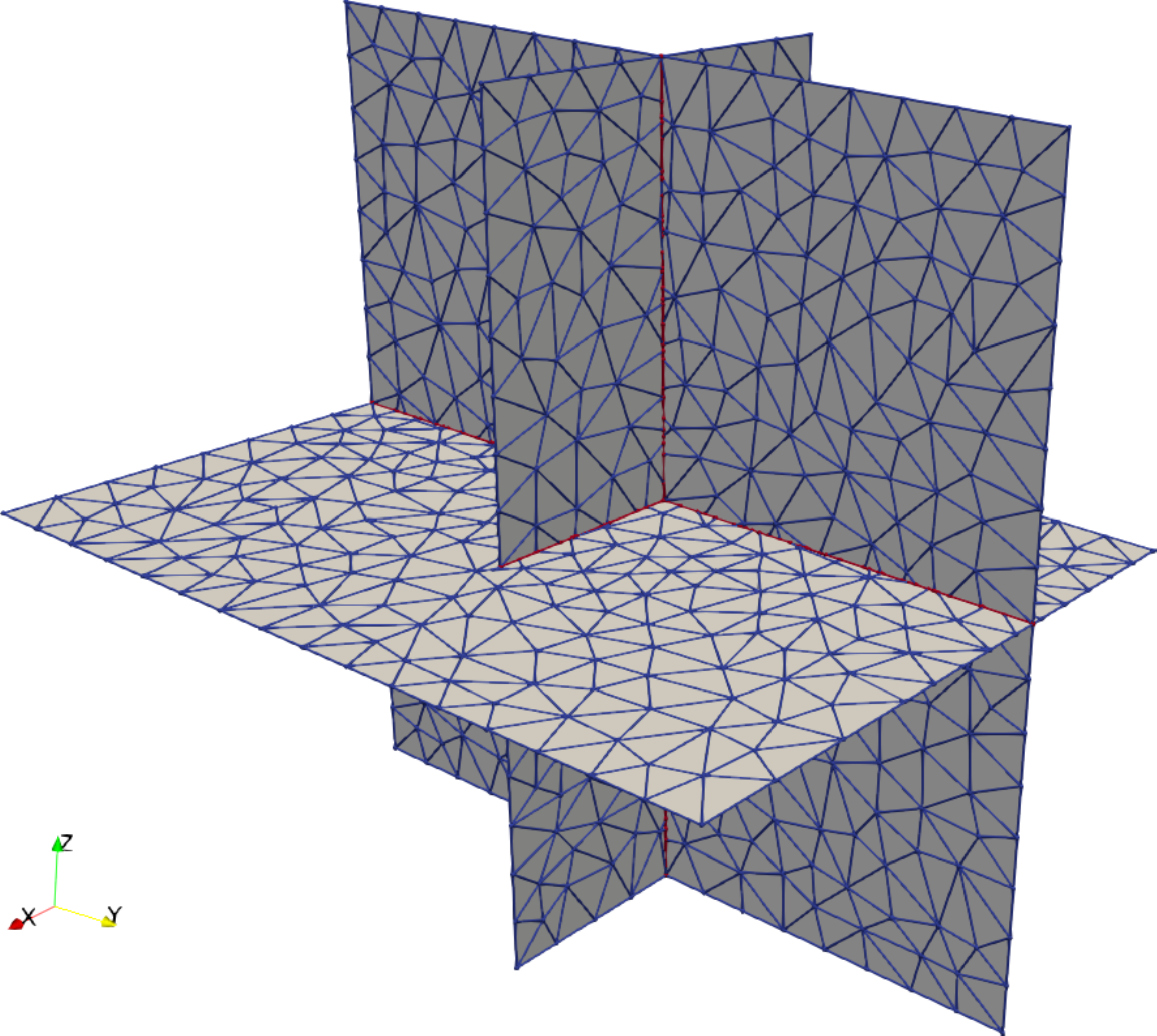}}
    \caption{Network 1; discretization with mesh M1 (a) and M2 (b), with constrained edges along the interfaces marked in red.}
    \label{plot:RUN_T5_DMN}
\end{figure}

In the first test, we analyze the convergence properties of the proposed method on a simple network that is composed by three fractures aligned with the reference system:
\begin{eqnarray*}
    \Fracture_1 &= \{(x, y, z) \in \RSet[3] : -1.0 \leq x \leq 0.5,\ -1.0 \leq y \leq 1.0,\ z = 0.0 \},\\
    \Fracture_2 &= \{(x, y, z) \in \RSet[3] : -1.0 \leq x \leq 0.0,\ y = 0.0,\ -1.0 \leq z \leq 1.0 \},\\
    \Fracture_3 &= \{(x, y, z) \in \RSet[3] : x  =0.5,\ -1.0 \leq y \leq 1.0,\ -1.0 \leq z \leq 1.0 \}.\\
\end{eqnarray*}
These fractures intersect along three interfaces, $\Trace_1, \Trace_2$, and $\Trace_3$, see Figure~\ref{plot:RUN_T5_DMN}.

\begin{figure}[htbp!]
    \centering
    \subfloat[\label{plot:RUN_T5_DMN_3_C2_ZOOM}]{\includegraphics[width=0.48\textwidth]{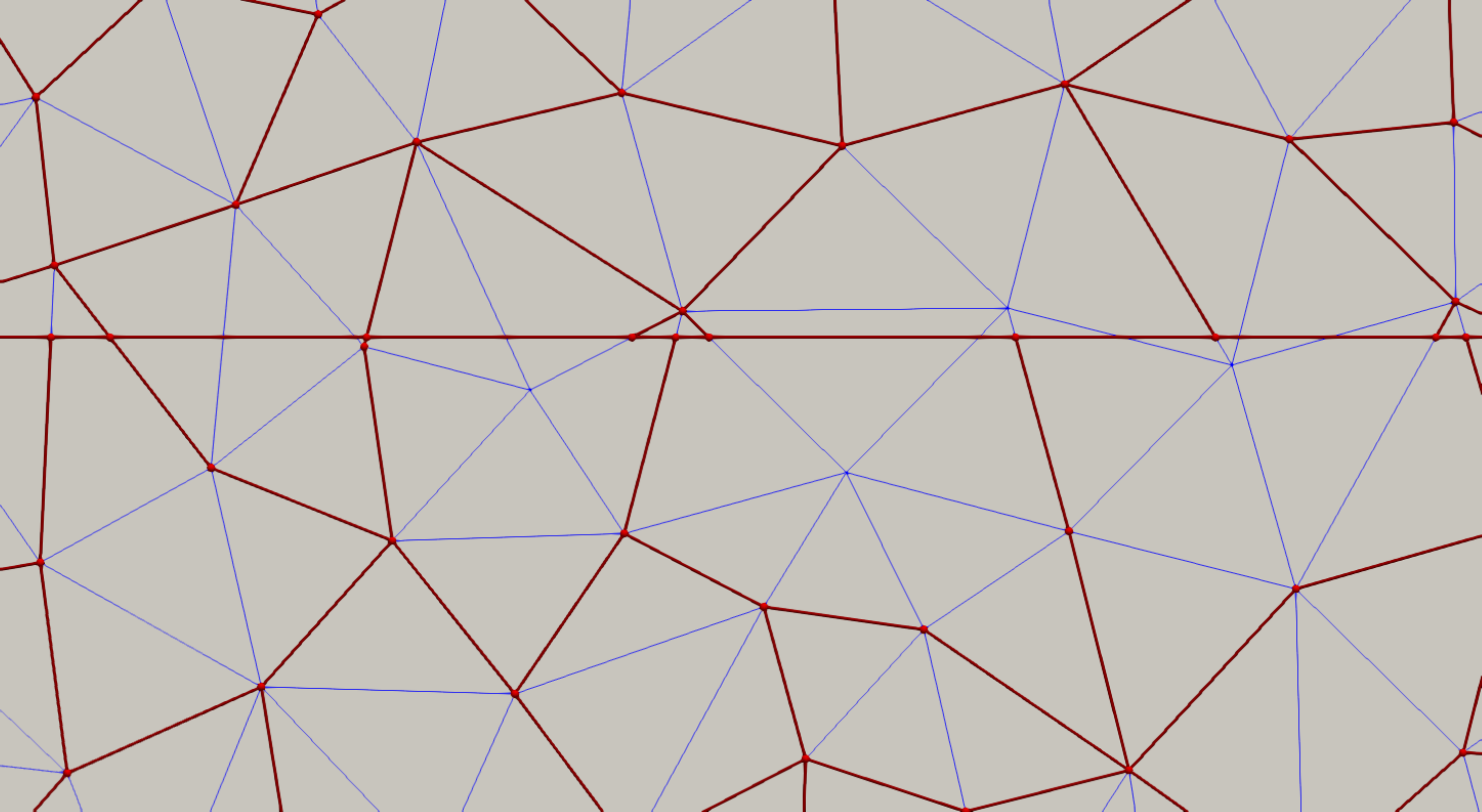}}\
    \subfloat[\label{plot:RUN_T5_DMN_3_C2}]{\includegraphics[width=0.38\textwidth]{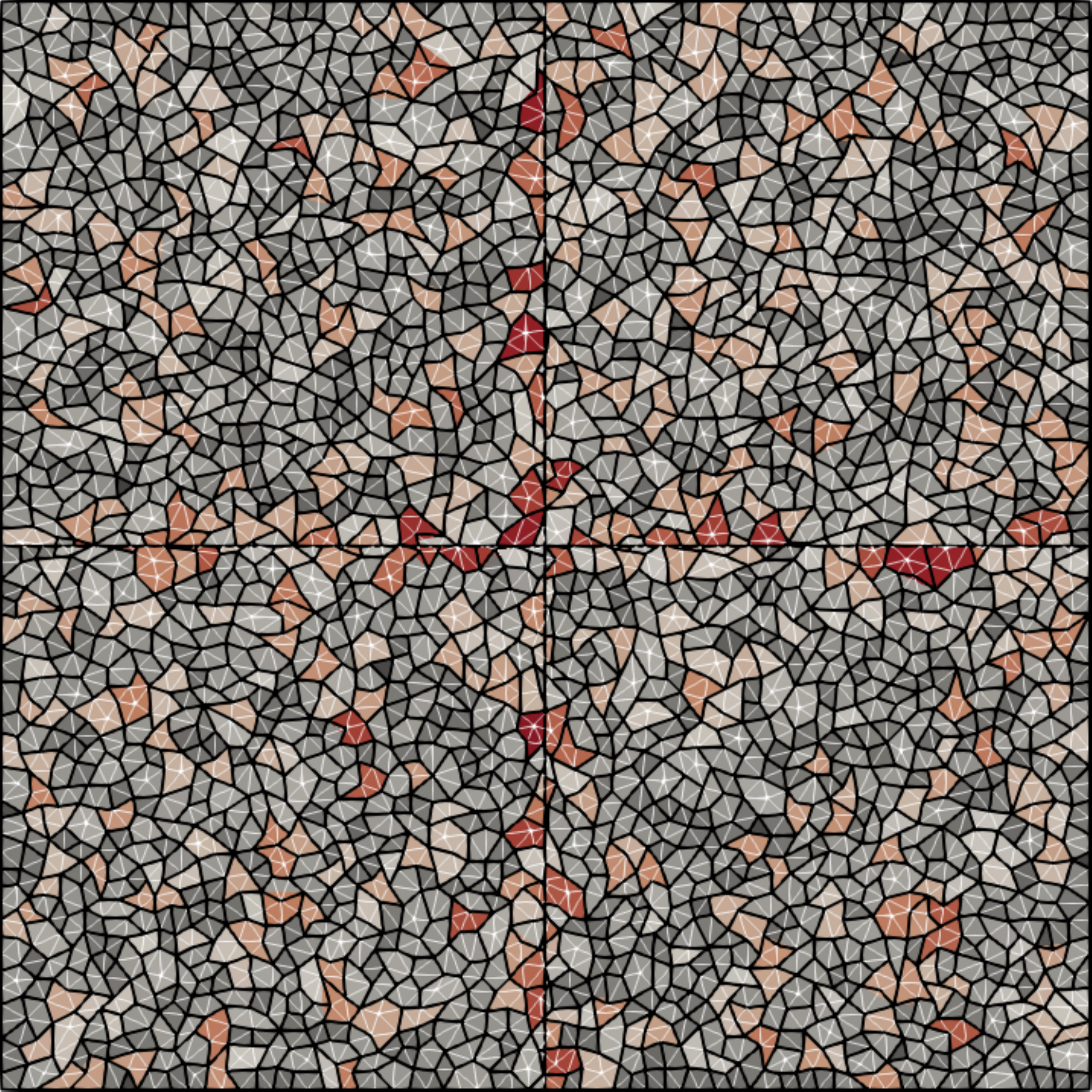}}
    \caption{Network 1, mesh M3; (a) visual comparison between fracture $F_3$ with $\AgglomerationParameter = 0.0$ (blue lines) and $\AgglomerationParameter = 1.0$ (red lines) (b) fracture $F_3$ with $\AgglomerationParameter = 1.0$ colored w.r.t. the $\varrho$ value on each cell, from red ($\varrho \approx 0$) to white ($\varrho \approx 1$).}
    \label{plot:RUN_T5_DMN_3}
\end{figure}

We test three tessellations of decreasing size, labeled M1, M2, and M3.
We first discretize the network with three triangular meshes fixing the maximum area of the triangular elements to $10^{-1}, 10^{-2}$, and $10^{-3}$.
Then, we cut the element along the fracture intersection segments as described in Section~\ref{subsec:problem}.
Finally, we apply the Quality Agglomeration algorithm described in Section~\ref{sec:agglomeration} to M1, M2, and M3, with $\AgglomerationParameter=0.25$ and $\AgglomerationParameter=1.0$.
In Figure~\ref{plot:RUN_T5_DMN_3_C2_ZOOM} we present a visualization of the effects of the agglomeration: the algorithm deletes several edges from the original mesh with $\AgglomerationParameter = 0.0$ (in blue).
By observing the localization of the quality indicator $\varrho$ on a single fracture (Figure~\ref{plot:RUN_T5_DMN_3_C2}), we notice how the most pathological elements (in red) are around the interfaces.

Table~\ref{tab:RUN_T5_ENERGY} summarizes the results of the optimization process applied to meshes M1, M2, and M3 on fracture $\Fracture_3$ (we do not report $\Fracture_1$ and $\Fracture_2$, as the measured values were comparable).
We analyze the energy functional \eqref{eq:graphcut:energy} and compare its value $\calE_1$ over the original meshes against its value $\calE_2$ over the meshes agglomerated with the two $\AgglomerationParameter$ values.
In columns $(\Energy_1 - \Energy_2)$ we report the percentage of energy saved by the optimization process, compared to the total original energy of the non-agglomerated mesh.
In columns $\Energy_{2, \text{dc}}$ and $\Energy_{2, \text{sc}}$, we report the contribution of the data cost and the smoothness cost to the total energy $\Energy_2$.
We remark that, while $\Energy_{1, \text{dc}}$ and $\Energy_{1, \text{sc}}$ are equal over the same mesh, their combination $\calE_1$ also depends on the value of $\AgglomerationParameter$.
We can observe how the process always leads to a reduction of the total energy within a small number of iterations (this latter is reported in column $It$).
However, using $\AgglomerationParameter = 0.25$ leads to a more conservative strategy (around $2 \%$ of energy saving), counter to the aggressive approach of using $\AgglomerationParameter = 1.0$ (more than $30 \%$ of energy saving).

\begin{table}[htbp!]
    \caption{Network 1; energy functional over $\Fracture_3$ before ($\calE_1$) and after ($\calE_2$) the optimization, with the detail of the data cost ($\calE_{dc}$) and smoothness cost ($\calE_{sc}$) contribution. 
    Column \textit{It} shows the number of iterations needed to converge.}
    \centering
    \includegraphics[width=0.65\textwidth]{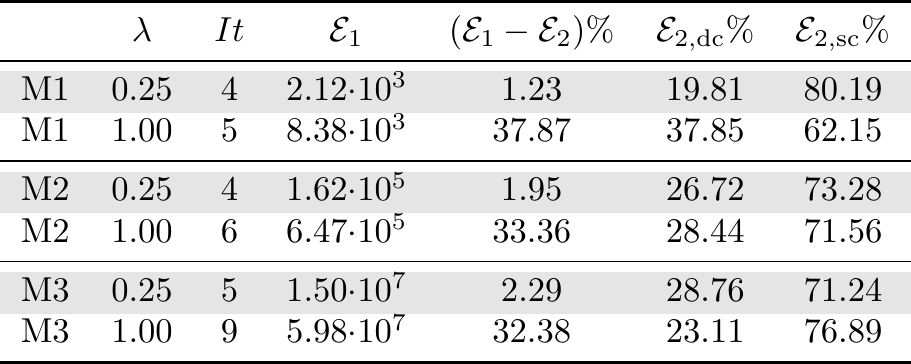}
    \label{tab:RUN_T5_ENERGY}
\end{table}

In Table~\ref{tab:RUN_T5}, we report the number of elements $\seminorm{\Mesh}$ of the meshes for the different $\AgglomerationParameter$ values, and the relative number of DOFs (which varies according to the VEM order $k$).
We observe that $\seminorm{\Mesh}$ is roughly $30\%$ smaller for $\AgglomerationParameter=0.25$ and $70\%$ smaller for $\AgglomerationParameter=1.0$ than the value for the original mesh for $\AgglomerationParameter=0$.
As far as DOFs are concerned, we recall from Section~\ref{subsec:vem} that, for $k=1$, the degrees of freedom correspond to the mesh vertices; thus, since the agglomeration algorithm mainly removes edges and elements, the number of DOFs remains almost untouched in the linear formulation, especially with $\AgglomerationParameter=0.25$.
As soon as we raise the order, however, the difference becomes significant: around $50\%$ reduction in the case of $k=2$ and $\AgglomerationParameter=1.0$, and slightly more for $k=3$.
\begin{table}[htbp!]
    \caption{Network 1; analysis of the numerical errors for the different meshes and VEM orders $\VemOrder$. $\AgglomerationParameter$ is the agglomeration parameter, $\seminorm{\Mesh}$ the number of elements in the mesh, $\DiscreteError$ the discrete error, $\Matrixize{A}$ the stiffness matrix and NNZ its non-zero elements.}
    \centering
    \begin{subtable}{.8\textwidth}
        \centering
        \caption{$\VemOrder = 1$}
        \includegraphics[width=\textwidth]{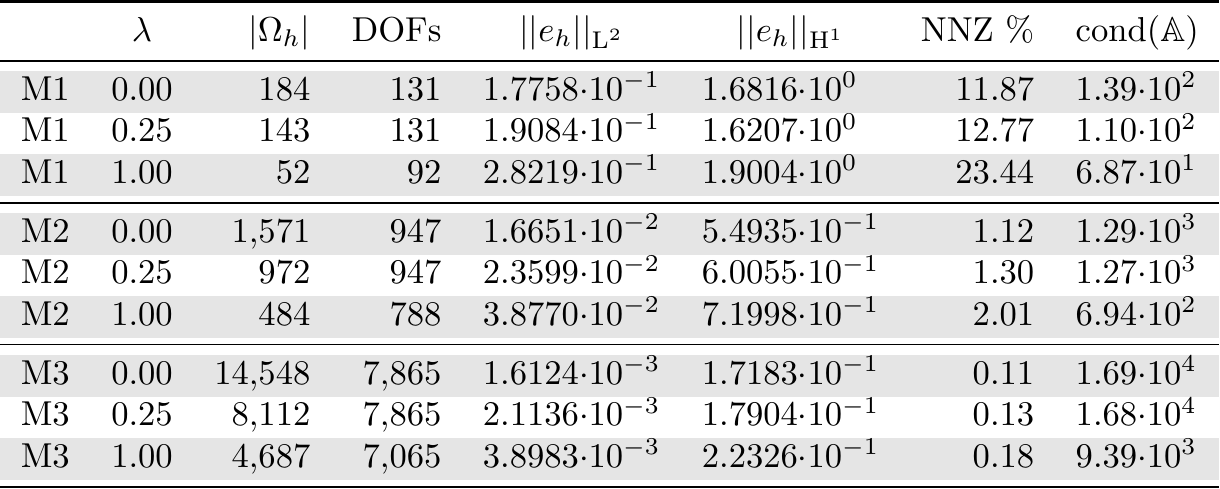}
        \label{tab:RUN_T5_O1}
    \end{subtable}
    \begin{subtable}{.8\textwidth}
        \centering
        \caption{$\VemOrder = 2$}
        \includegraphics[width=\textwidth]{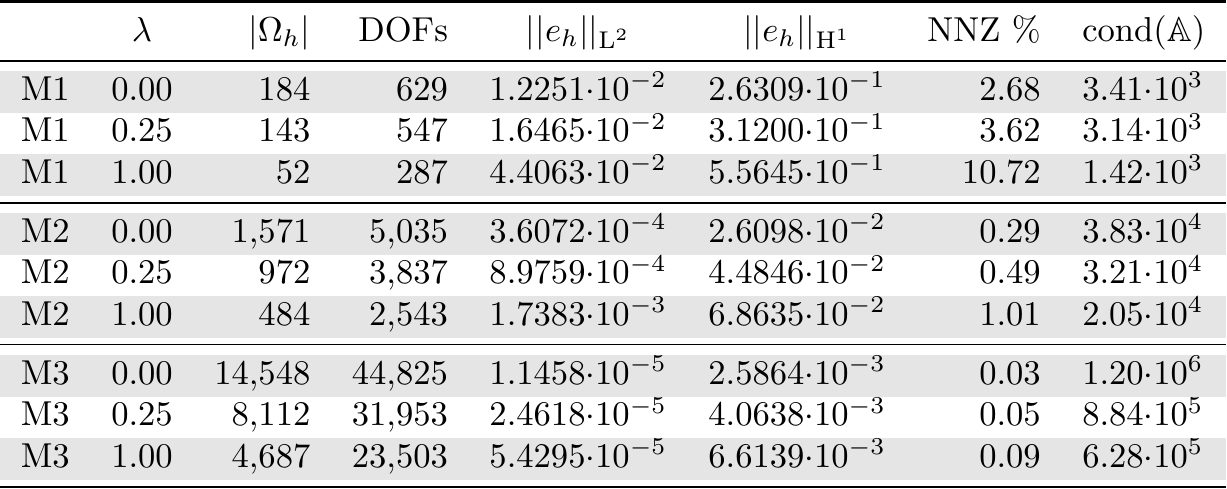}
        \label{tab:RUN_T5_O2}
    \end{subtable}
    \begin{subtable}{.8\textwidth}
        \centering
        \caption{$\VemOrder = 3$}
        \includegraphics[width=\textwidth]{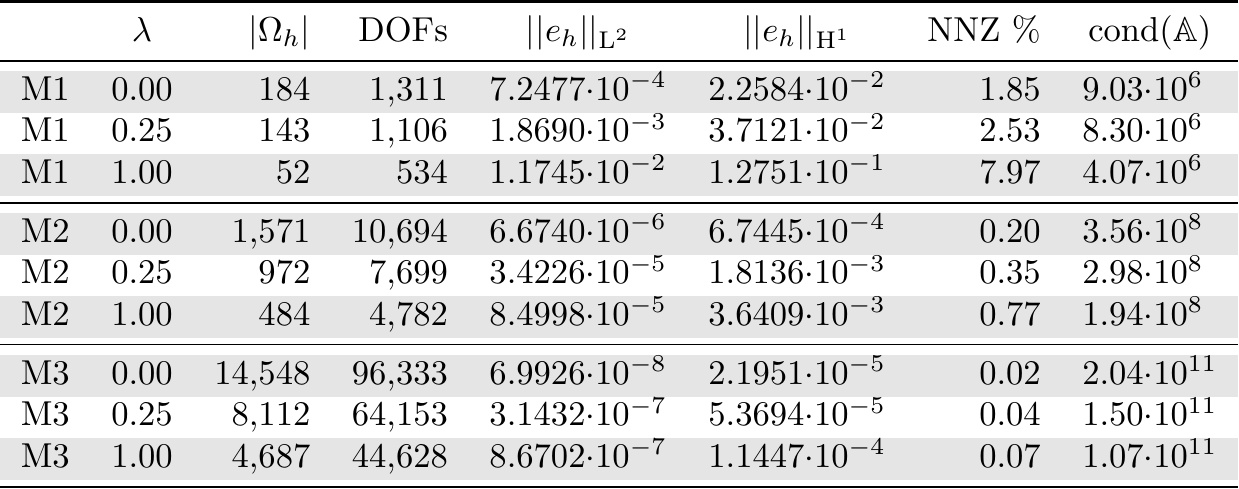}
        \label{tab:RUN_T5_O3}
    \end{subtable}
    \label{tab:RUN_T5}
\end{table}

\begin{figure}[htbp!]
    \subfloat[$\AgglomerationParameter = 0.0$\label{plot:RUN_T5_SOL_C0}]{\includegraphics[width=0.32\textwidth]{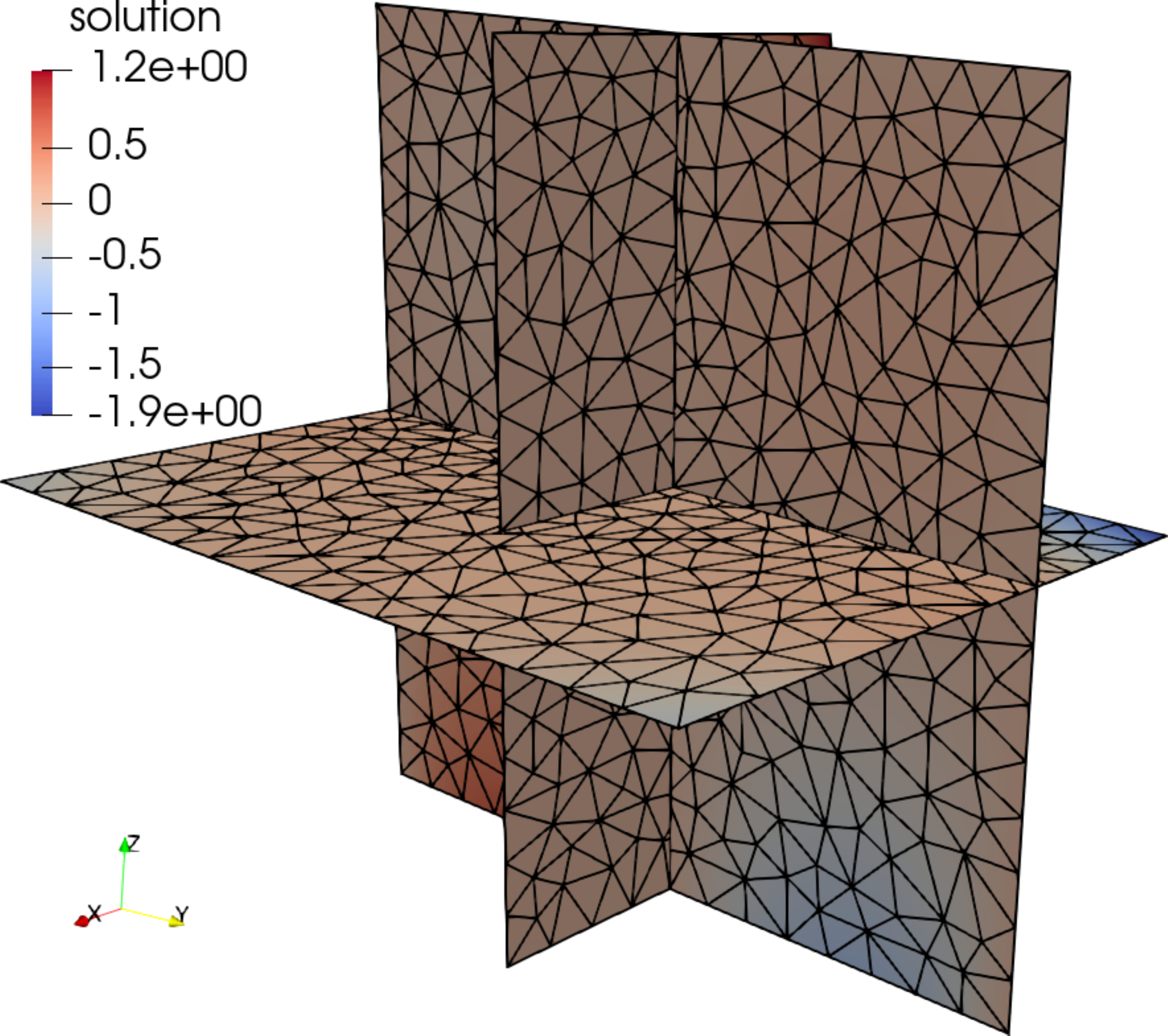}}\
    \subfloat[$\AgglomerationParameter = 0.25$\label{plot:RUN_T5_SOL_C1}]{\includegraphics[width=0.32\textwidth]{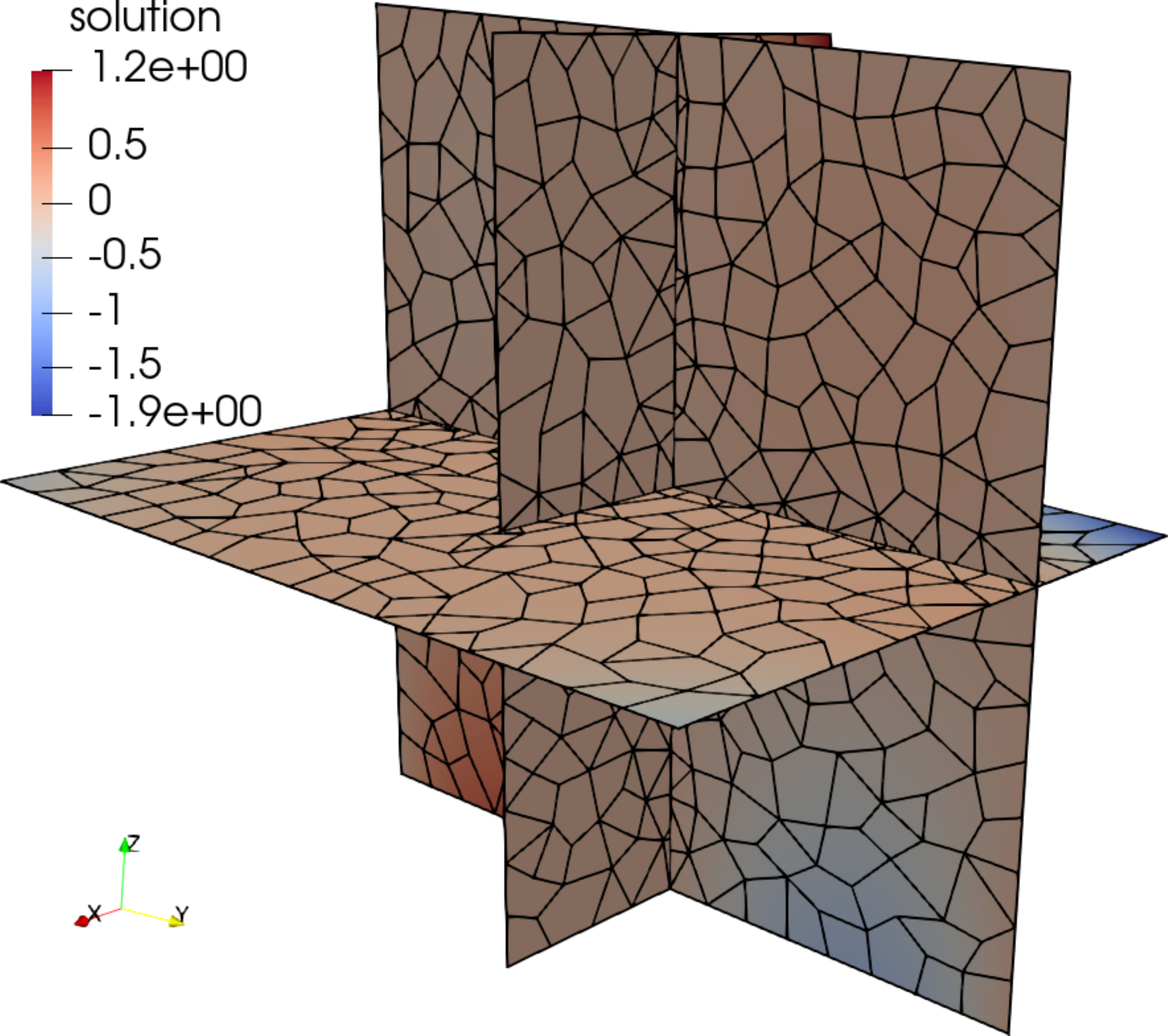}}\
    \subfloat[$\AgglomerationParameter = 1.0$\label{plot:RUN_T5_SOL_C2}]{\includegraphics[width=0.32\textwidth]{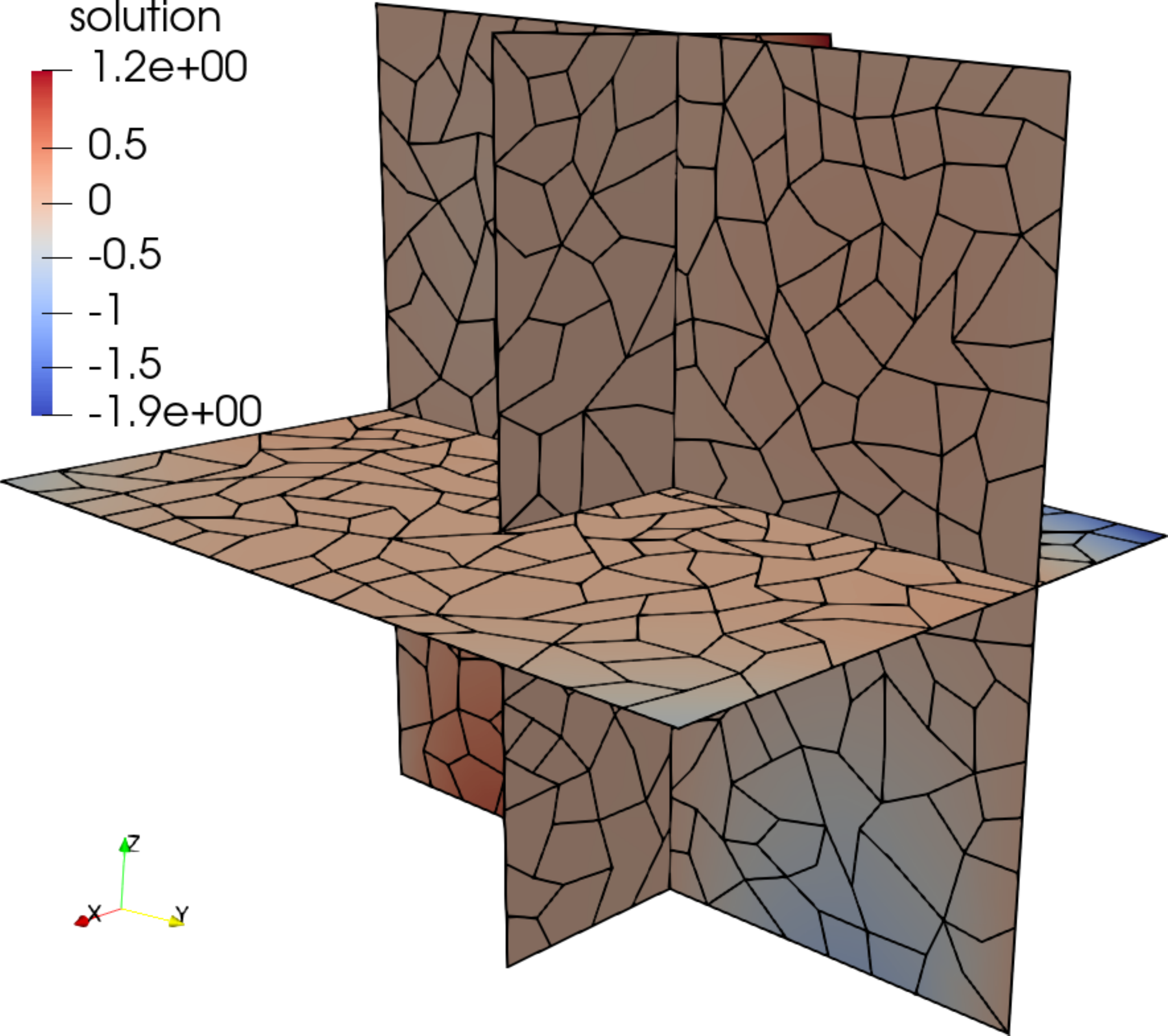}}
    \caption{Network 1, mesh M2; solution computed with $\VemOrder = 3$ for each $\AgglomerationParameter$ value.}
    \label{plot:RUN_T5_SOL}
\end{figure}

We set the numerical problem by imposing the Dirichlet boundary conditions in accordance with the exact solution
\begin{eqnarray*}
    \HydraulicHead_1(x,y) = -\frac{1}{10} (\frac{1}{2}+x) [x^3+8xy(x^2+y^2)\textrm{atan2}(y,x)] & \text{on } \Fracture_1,\\
    \HydraulicHead_2(x,z) = -\frac{1}{10} (\frac{1}{2}+x)x^3 + \pi\frac{4}{5}(\frac{1}{2}+x)x^3 \ABS{z} & \text{on } \Fracture_2,\\
    \HydraulicHead_3(y,z) = y(y-1)(y+1)z(z-1) & \text{on } \Fracture_3.
\end{eqnarray*}
In Figure~\ref{plot:RUN_T5_SOL} we show the solutions computed with the proposed method with $\VemOrder = 3$ on the mesh M2 for the three different levels of agglomeration.
We can appreciate that no qualitative differences can be observed in the three results, despite a significant reduction of the DOFs (up to $50\%$) in the agglomeration process, see Table~\ref{tab:RUN_T5_O1}.

\begin{figure}[htbp!]
    \setcounter{subfigure}{0}
    \centering
    \subfloat[$\VemOrder = 1$\label{plot:RUN_T5_CONV_O1}]{\includegraphics[width=0.3\textwidth]{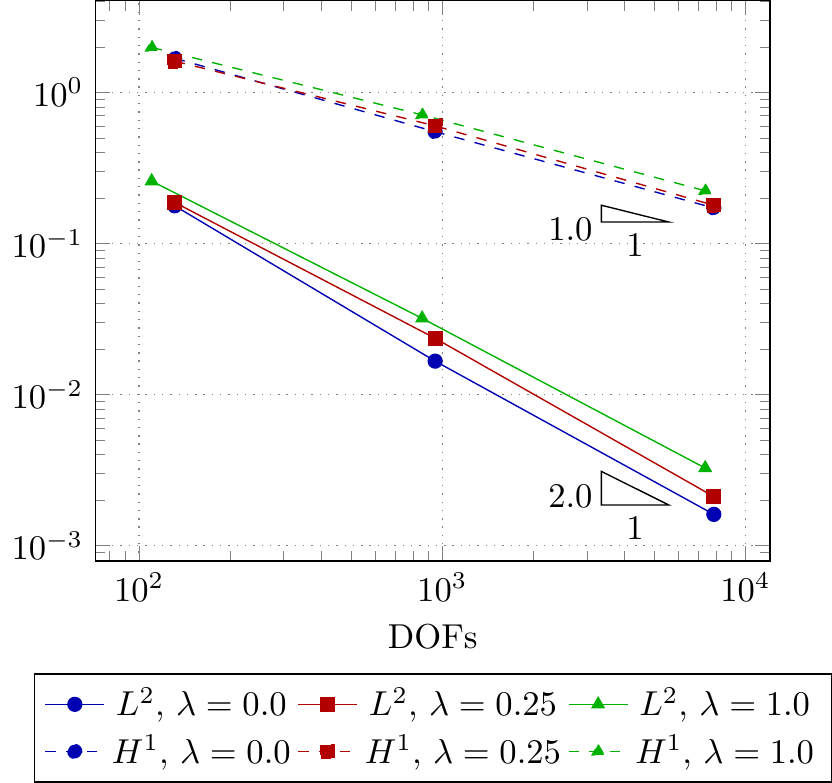}}
    \subfloat[$\VemOrder = 2$\label{plot:RUN_T5_CONV_O2}]{\includegraphics[width=0.3\textwidth]{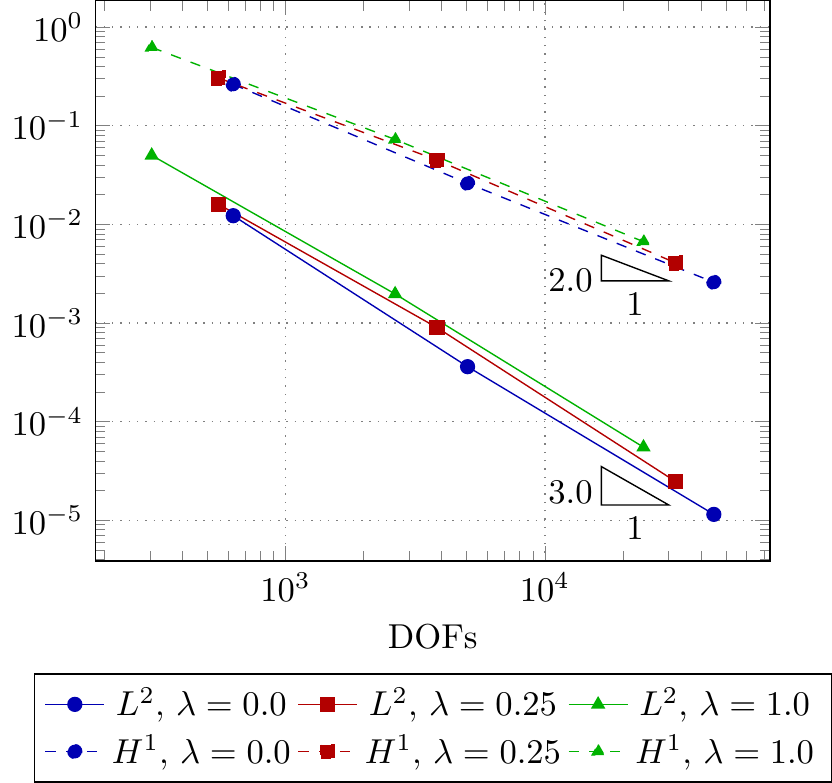}}
    \subfloat[$\VemOrder = 3$\label{plot:RUN_T5_CONV_O3}]{\includegraphics[width=0.3\textwidth]{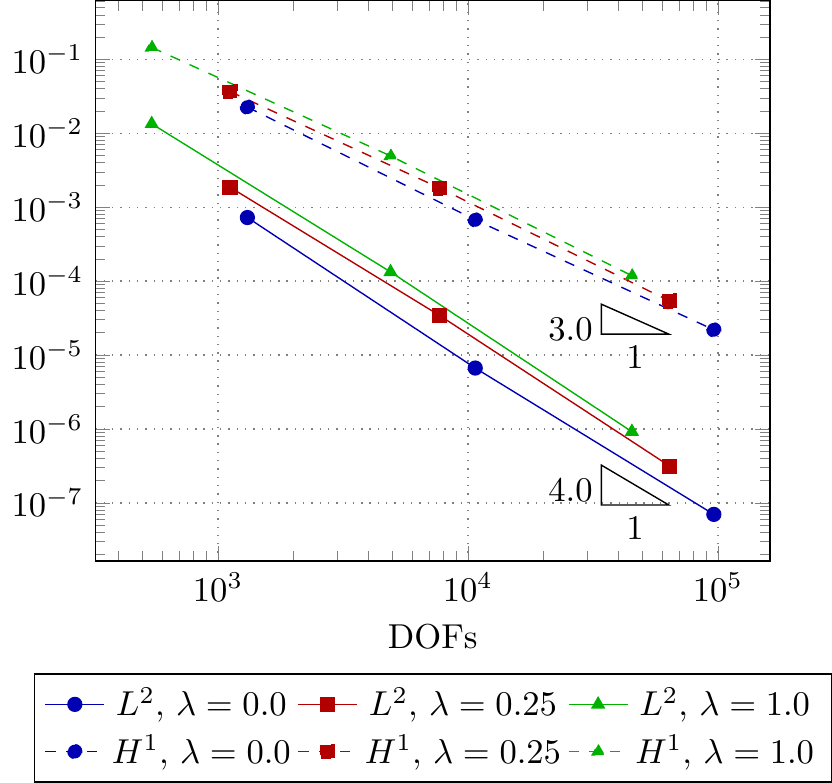}}
    \caption{Network 1; convergence curves of $L^2$ and $H^1$ errors for different $k$ values.}
    \label{plot:RUN_T5_CONV}
\end{figure}

\begin{figure}[htbp!]
    \setcounter{subfigure}{0}
    \centering
    \subfloat[$\VemOrder = 1$\label{plot:RUN_T5_EXP_O1}]{\includegraphics[width=0.3\textwidth]{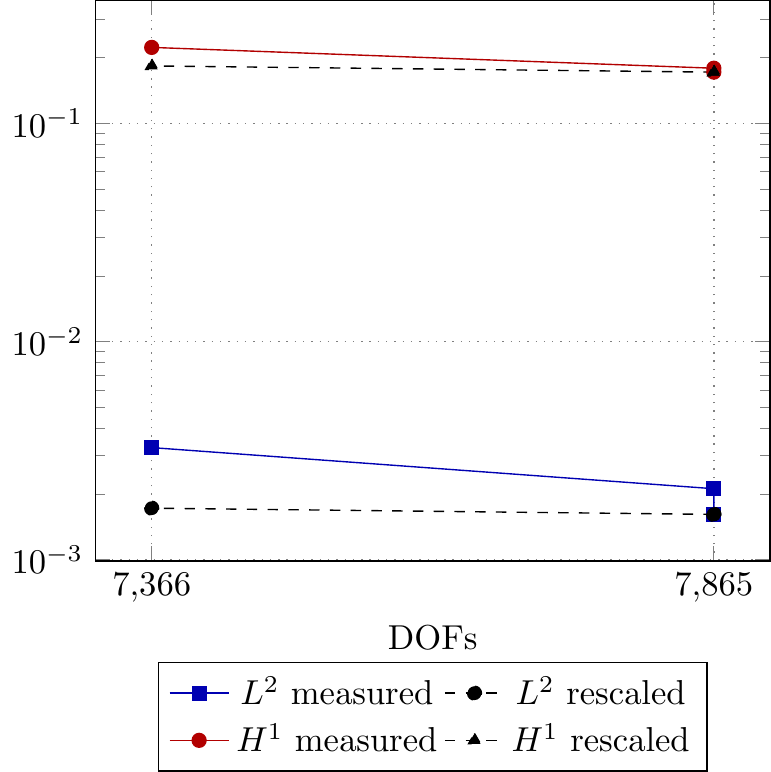}}
    \subfloat[$\VemOrder = 2$\label{plot:RUN_T5_EXP_O2}]{\includegraphics[width=0.3\textwidth]{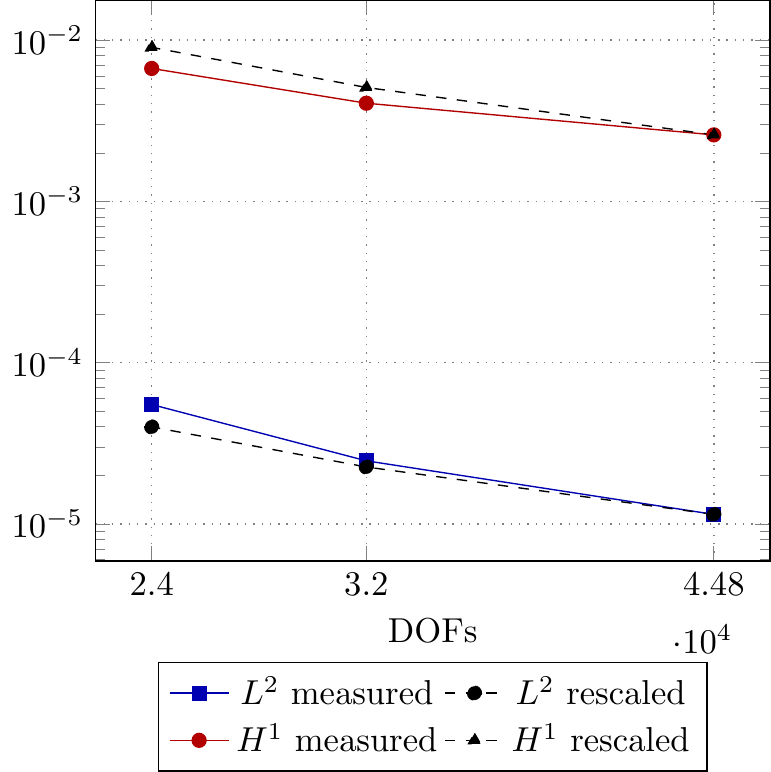}}
    \subfloat[$\VemOrder = 3$\label{plot:RUN_T5_EXP_O3}]{\includegraphics[width=0.3\textwidth]{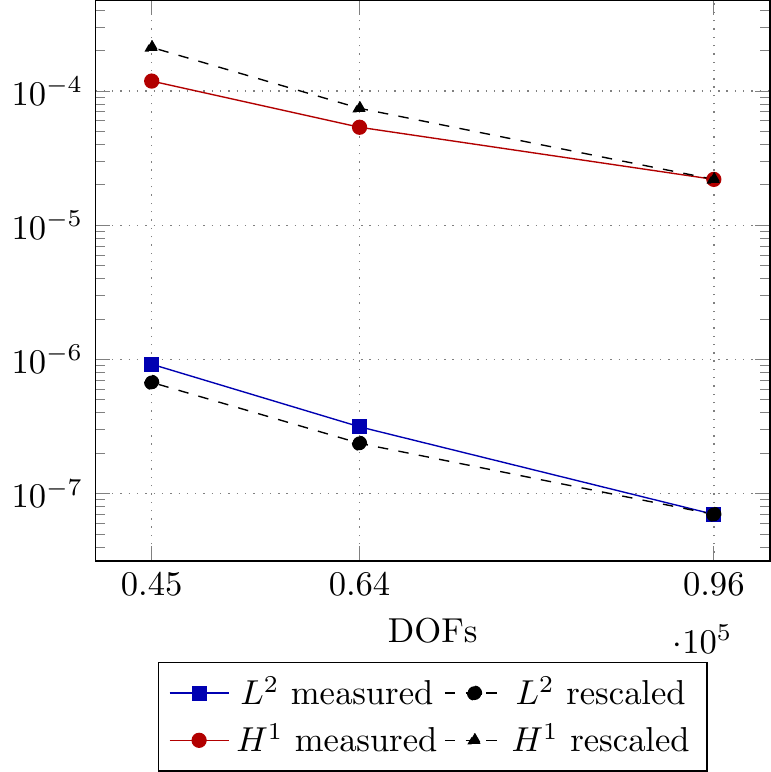}}
    \caption{Network 1, mesh M3; measured errors vs expected errors without mesh optimization.}
    \label{plot:RUN_T5_EXP}
\end{figure}

In the middle columns of Table~\ref{tab:RUN_T5} we show the global approximation errors
\begin{align*}
    \HilbertNorm[\SobolevL{2}]{\DiscreteError} := \HilbertNorm[\SobL{\Network}]{\HydraulicHead - \HydraulicHead_{\MeshParameter}}, \qquad
    \HilbertNorm[\SobolevH{1}]{\DiscreteError} := \seminorm[\SobH{1}{\Network}]{\HydraulicHead - \HydraulicHead_{\MeshParameter}}.
\end{align*}
Figure~\ref{plot:RUN_T5_CONV} shows the error curves versus the total number of degrees of freedom.
The slopes of these curves reflect the convergence rate $\VemConvergence$ for the polynomial degrees $k=1,2,3$.
Optimal convergence rates are evident in all plots.
The rate lines are similar in all the agglomerations.
These optimal rates are achieved thanks to the global conformity of the mesh, despite the low regularity of the solution around the extremity of $\Trace_3$ that falls inside $\Fracture_1$, see Figure~\ref{plot:RUN_T5_DMN}.
Note that the dots relative to agglomerated meshes are shifted leftwards, as they contain smaller numbers of DOFs.
To analyze pointwise the errors produced by an agglomerated mesh, we rescale the errors measured in the original mesh by the number of DOFs of the agglomerated mesh.
Let $\HilbertNorm[\star]{\DiscreteError[_1]}$ and $\HilbertNorm[\star]{\DiscreteError[_2]}$ for $\star\in\big\{\SobolevL{2},\SobolevH{1}\big\}$ the errors measured on two meshes with different mesh size $\MeshParameter_1$ and $\MeshParameter_2$.
We expect $\HilbertNorm[\star]{\DiscreteError[_1]} \MeshParameter_2^{\VemConvergence} \approx \HilbertNorm[\star]{\DiscreteError[_2]} \MeshParameter_1^{\VemConvergence}$.
Thus, in the plots of Figure~\ref{plot:RUN_T5_EXP}, we compare the errors $\HilbertNorm[\SobolevL{2}]{\DiscreteError}$ and $\HilbertNorm[\SobolevH{1}]{\DiscreteError}$ and the expected errors obtained by rescaling the non-optimized error norms $\HilbertNorm[\SobolevL{2}]{\DiscreteError}^{\AgglomerationParameter=0}$ and $\HilbertNorm[\SobolevH{1}]{\DiscreteError}^{\AgglomerationParameter=0}$ with the corresponding convergence rates; the order of magnitude of the two norms are comparable even when the reduction of DOFs is greater than $50 \%$, see Table~\ref{tab:RUN_T5_O2} as a reference.

\begin{figure}[htbp!]
    \subfloat[$\AgglomerationParameter = 0.0$, $\VemOrder = 1$\label{plot:RUN_T5_DMN_3_H1_O1_C0}]{\includegraphics[width=0.32\textwidth]{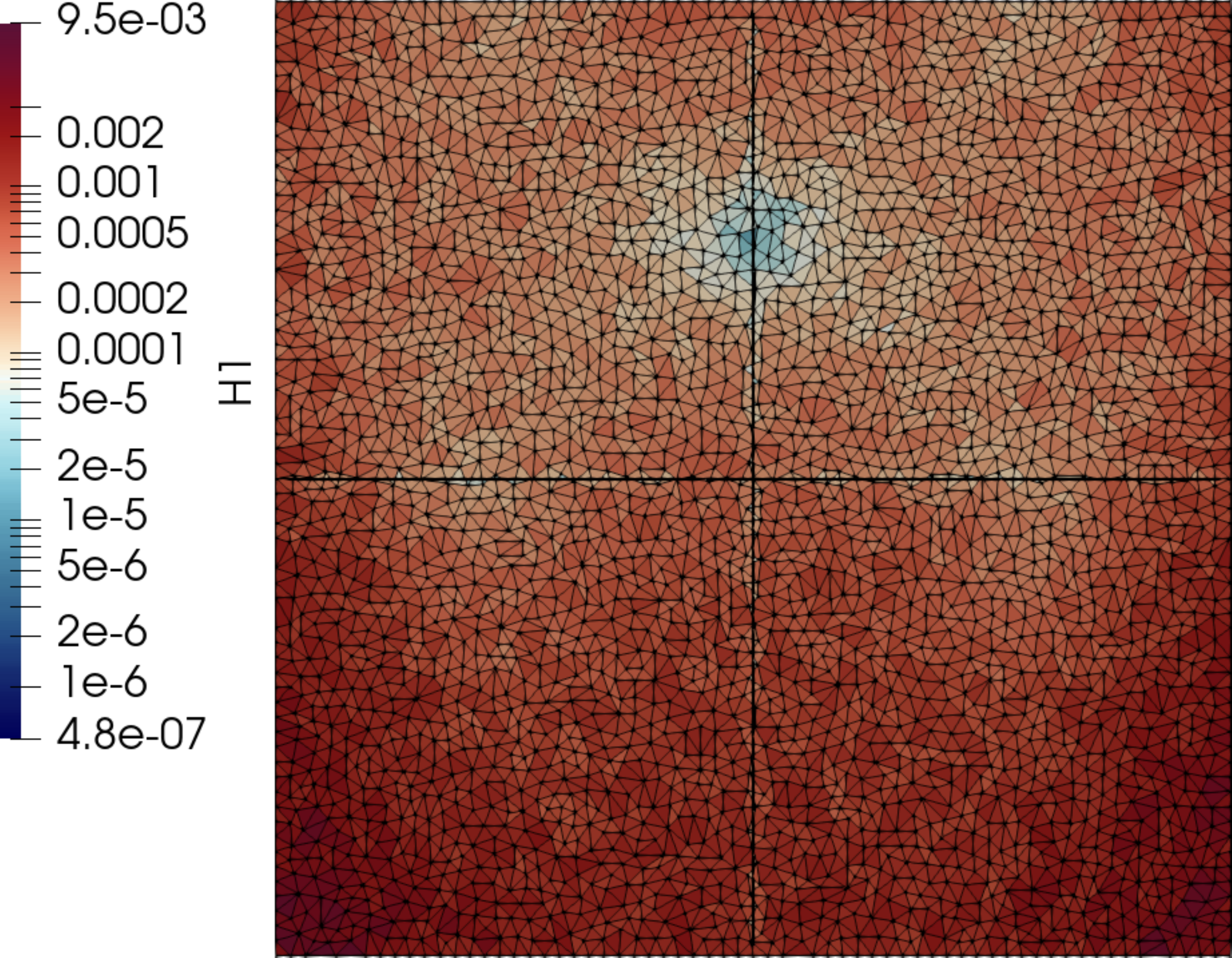}}\
    \subfloat[$\AgglomerationParameter = 0.25$, $\VemOrder = 1$\label{plot:RUN_T5_DMN_3_H1_O1_C1}]{\includegraphics[width=0.32\textwidth]{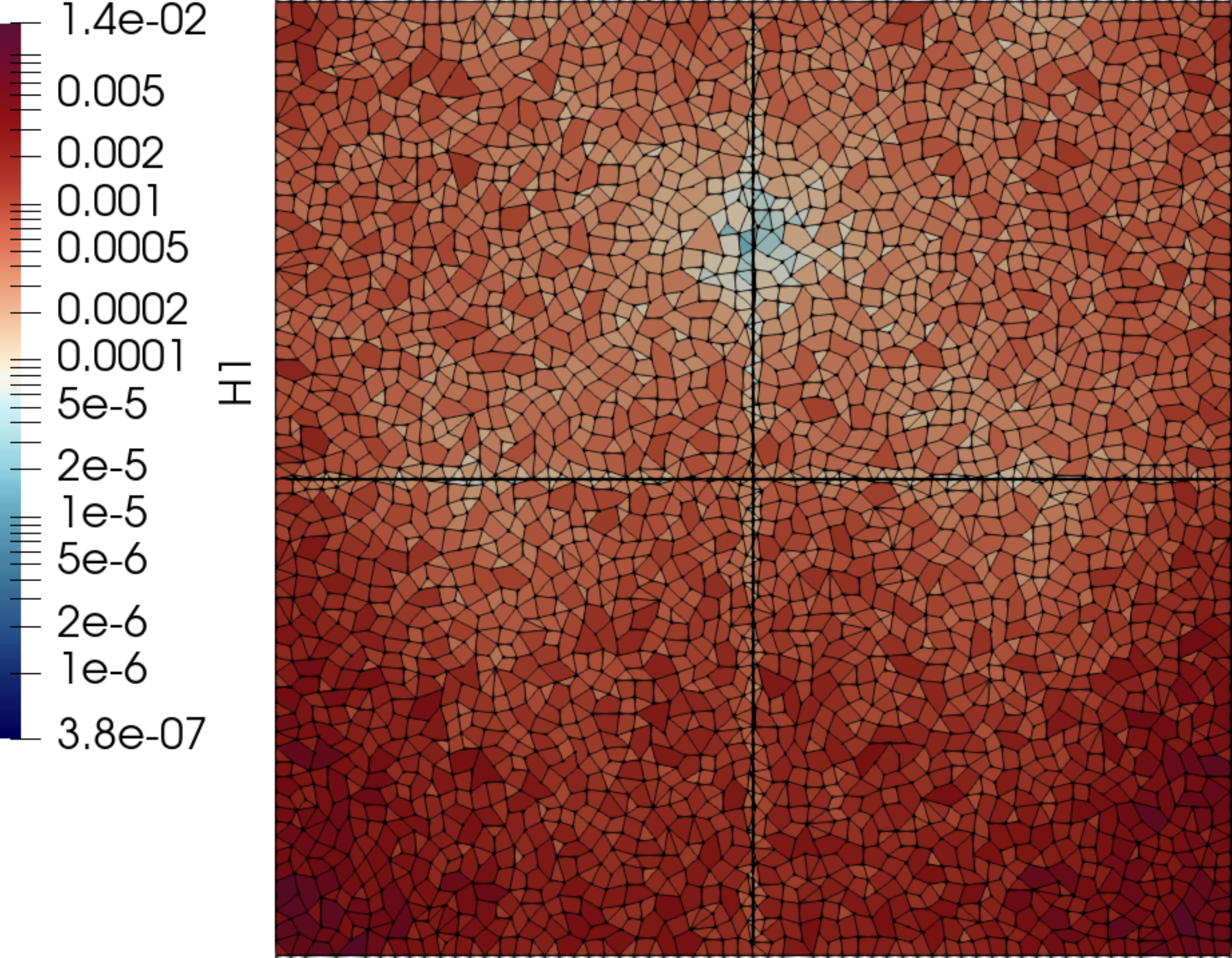}}\
    \subfloat[$\AgglomerationParameter = 1.0$, $\VemOrder = 1$\label{plot:RUN_T5_DMN_3_H1_O1_C2}]{\includegraphics[width=0.32\textwidth]{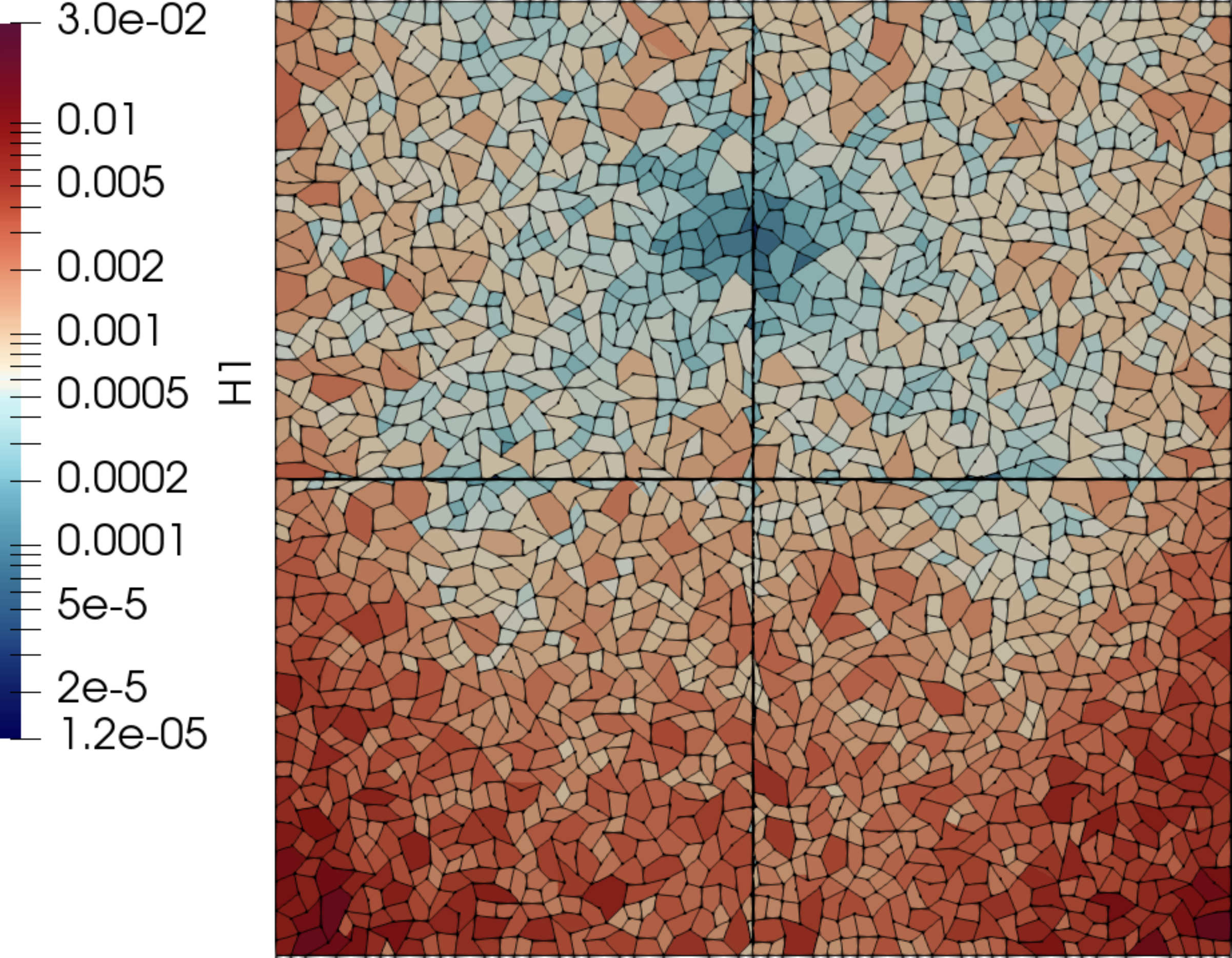}}\\
    \subfloat[$\AgglomerationParameter = 0.0$, $\VemOrder = 2$\label{plot:RUN_T5_DMN_3_H1_O2_C0}]{\includegraphics[width=0.32\textwidth]{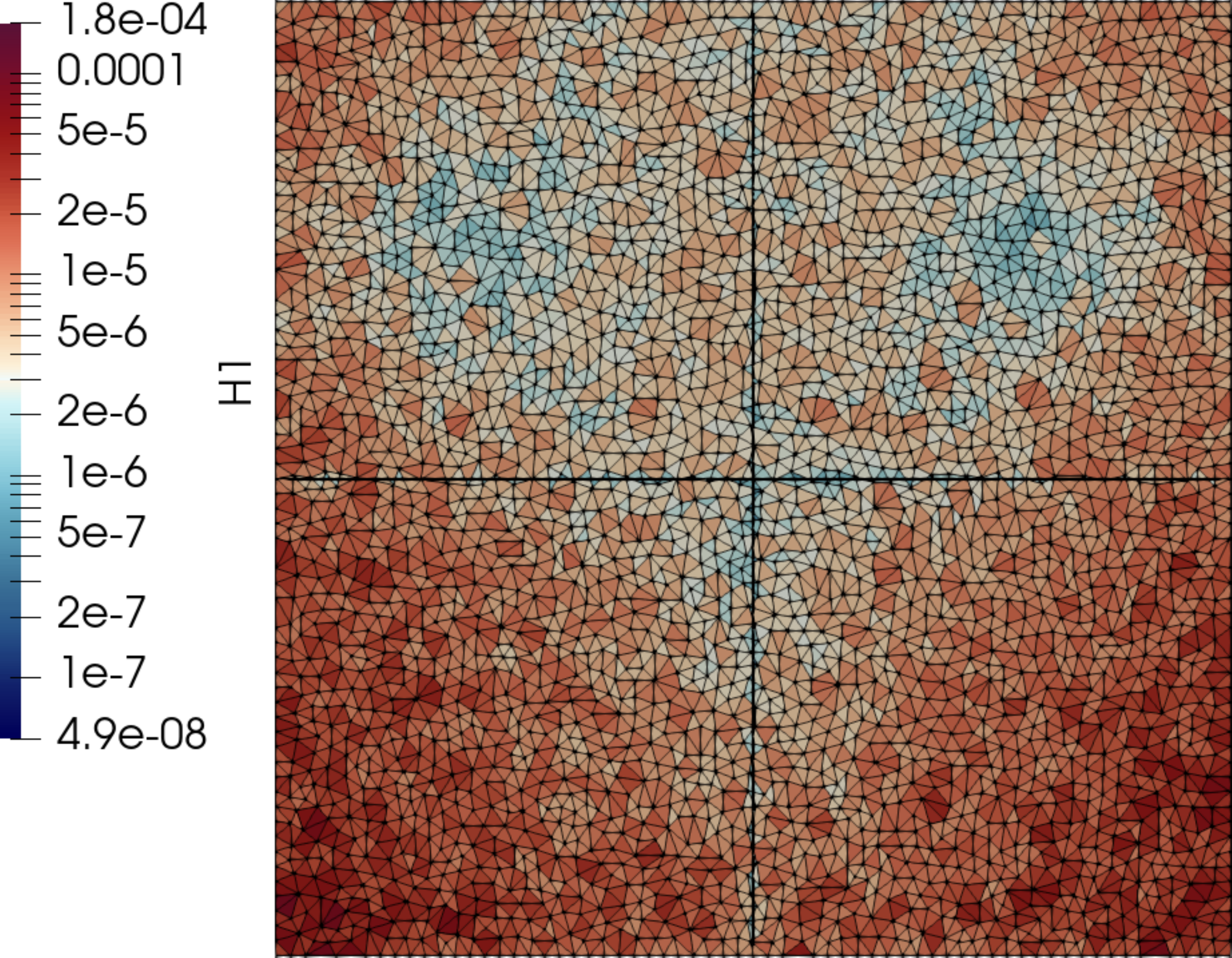}}\
    \subfloat[$\AgglomerationParameter = 0.25$, $\VemOrder = 2$\label{plot:RUN_T5_DMN_3_H1_O2_C1}]{\includegraphics[width=0.32\textwidth]{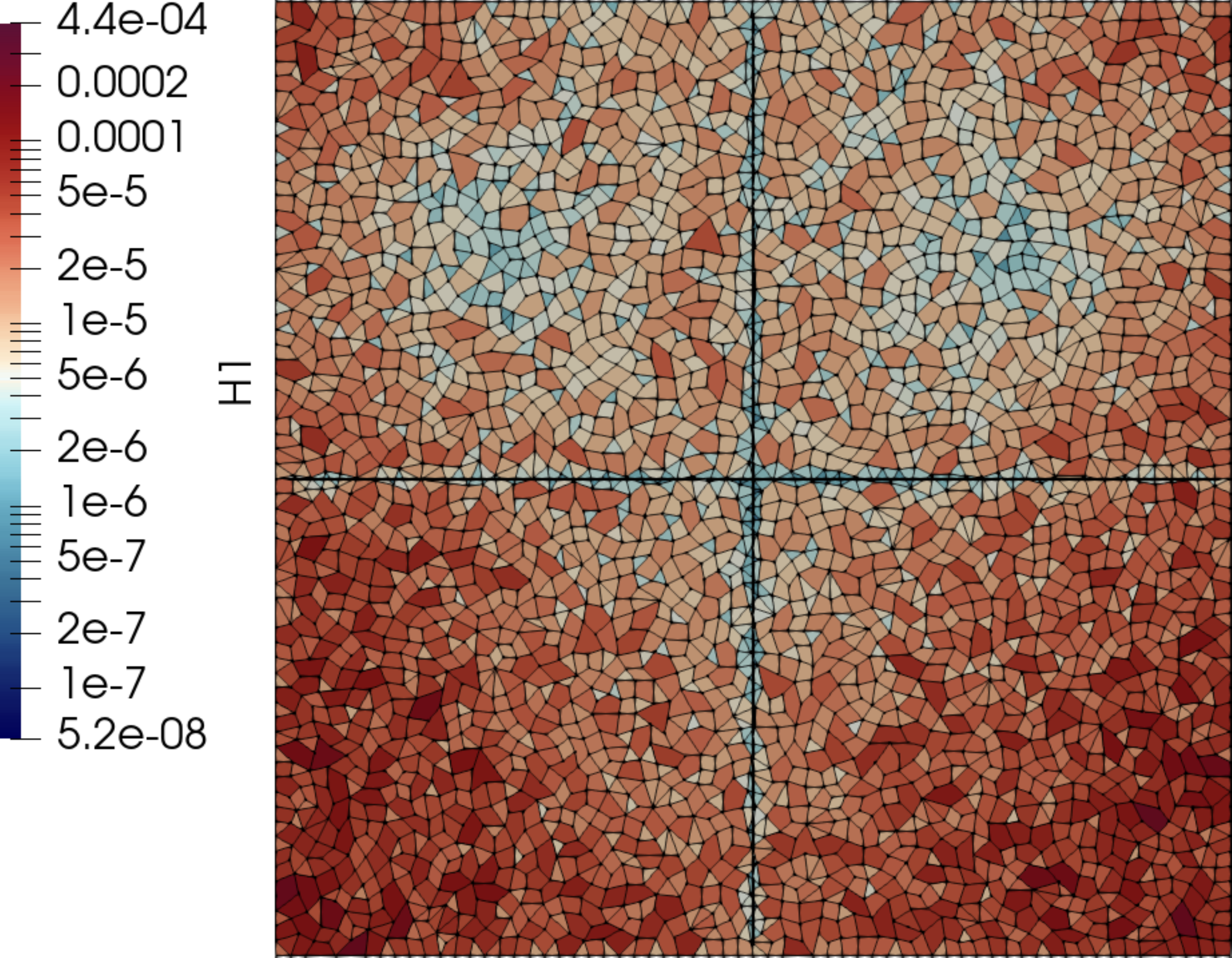}}\
    \subfloat[$\AgglomerationParameter = 1.0$, $\VemOrder = 2$\label{plot:RUN_T5_DMN_3_H1_O2_C2}]{\includegraphics[width=0.32\textwidth]{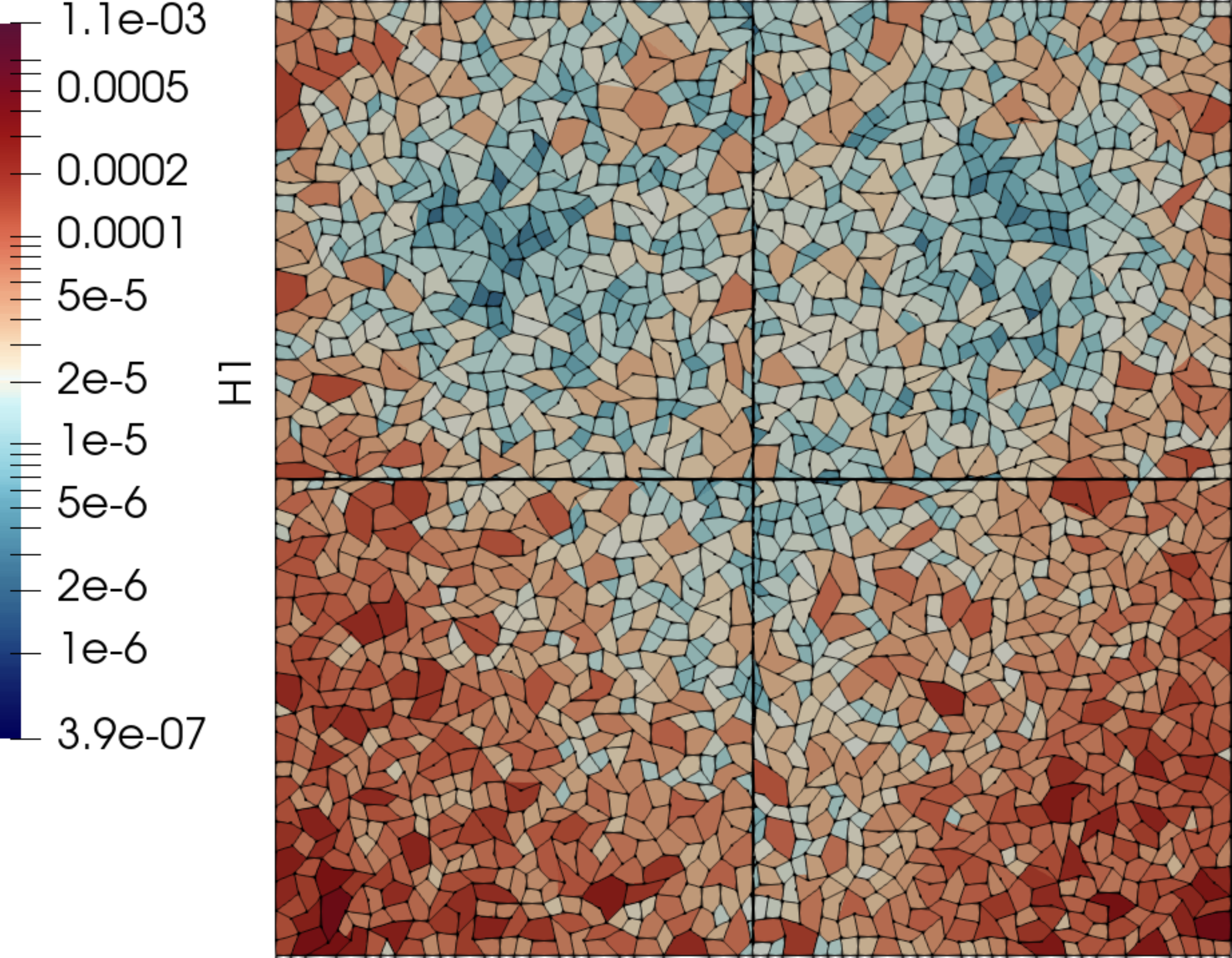}}\\
    \subfloat[$\AgglomerationParameter = 0.0$, $\VemOrder = 3$\label{plot:RUN_T5_DMN_3_H1_O3_C0}]{\includegraphics[width=0.32\textwidth]{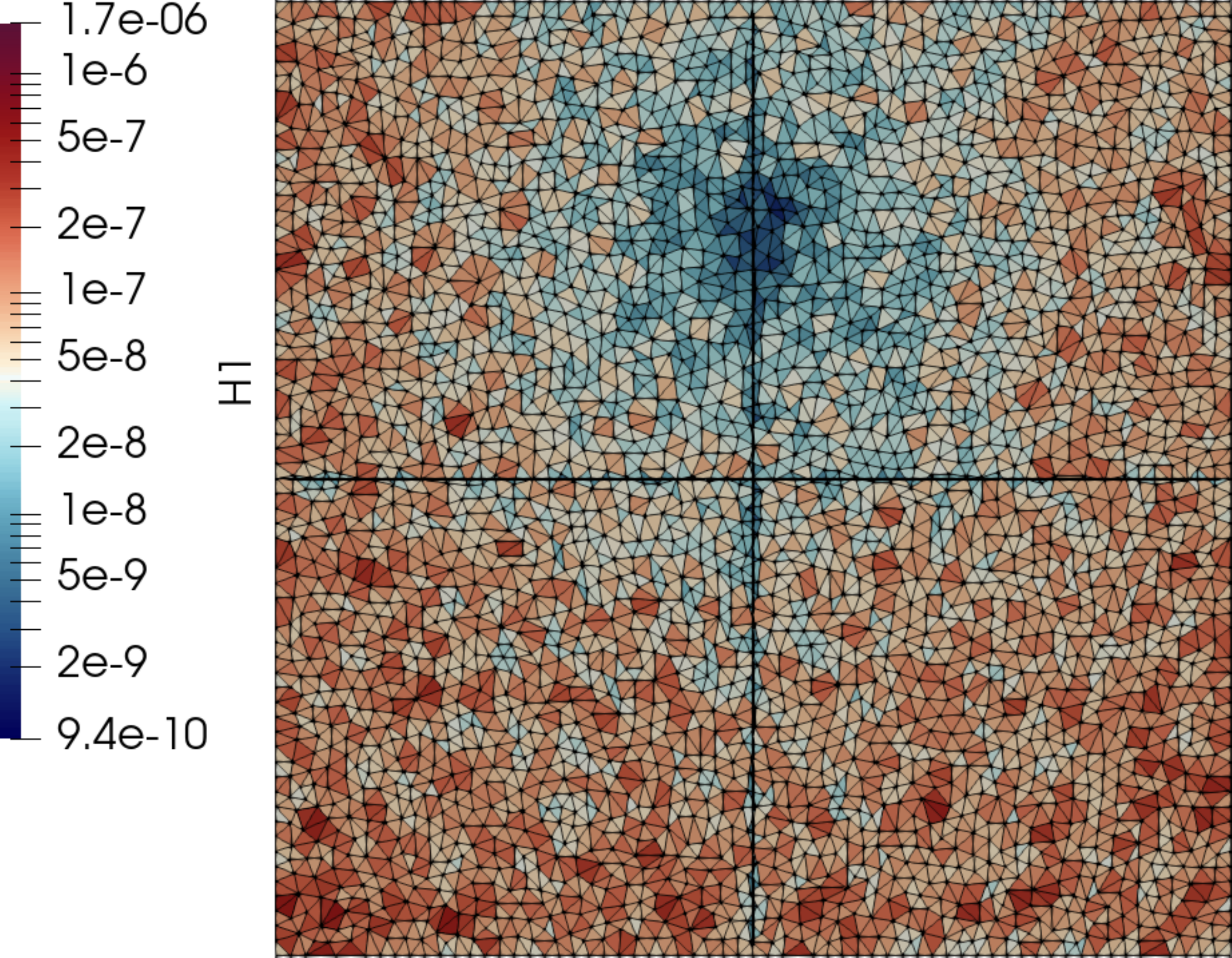}}\
    \subfloat[$\AgglomerationParameter = 0.25$, $\VemOrder = 3$\label{plot:RUN_T5_DMN_3_H1_O3_C1}]{\includegraphics[width=0.32\textwidth]{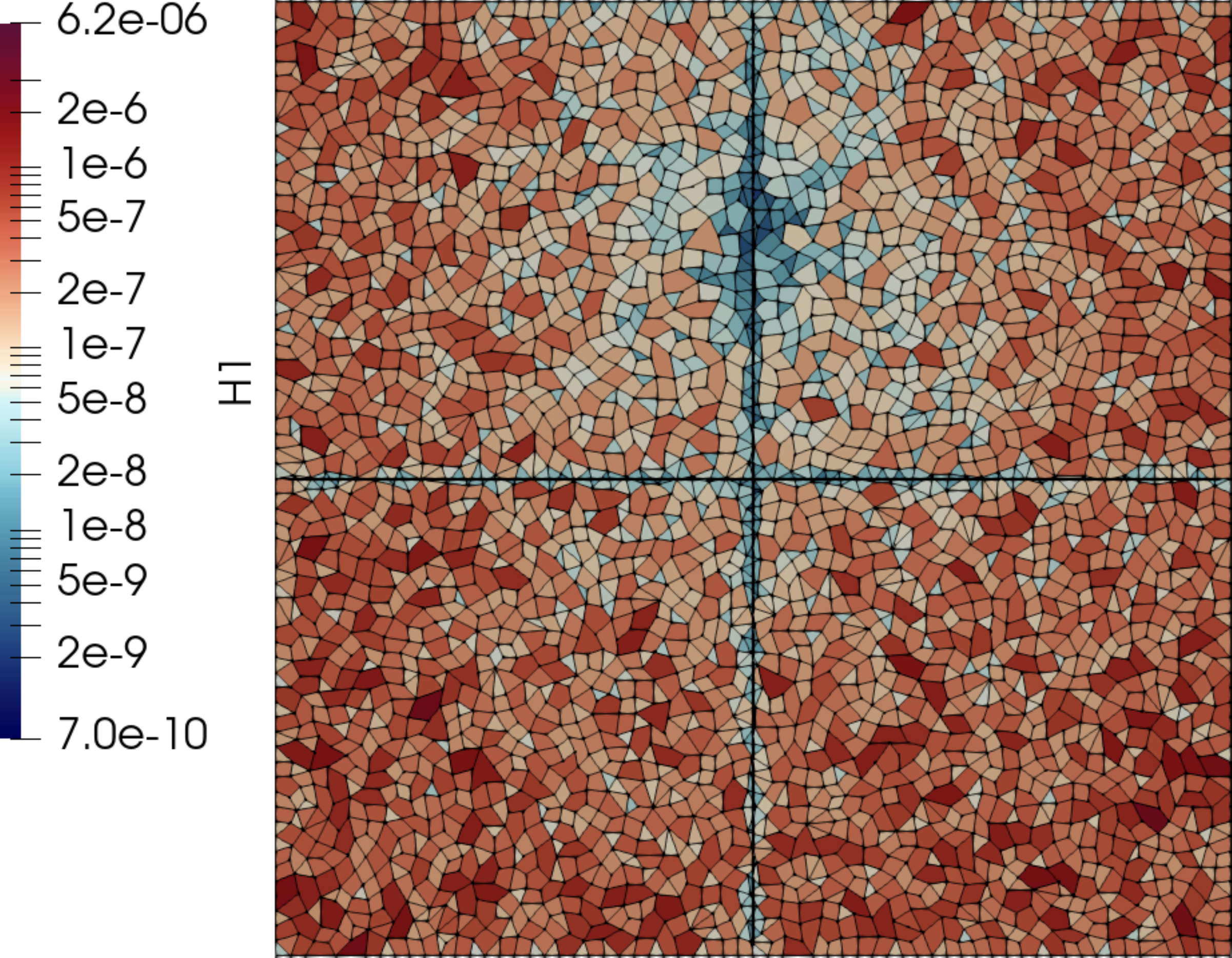}}\
    \subfloat[$\AgglomerationParameter = 1.0$, $\VemOrder = 3$\label{plot:RUN_T5_DMN_3_H1_O3_C2}]{\includegraphics[width=0.32\textwidth]{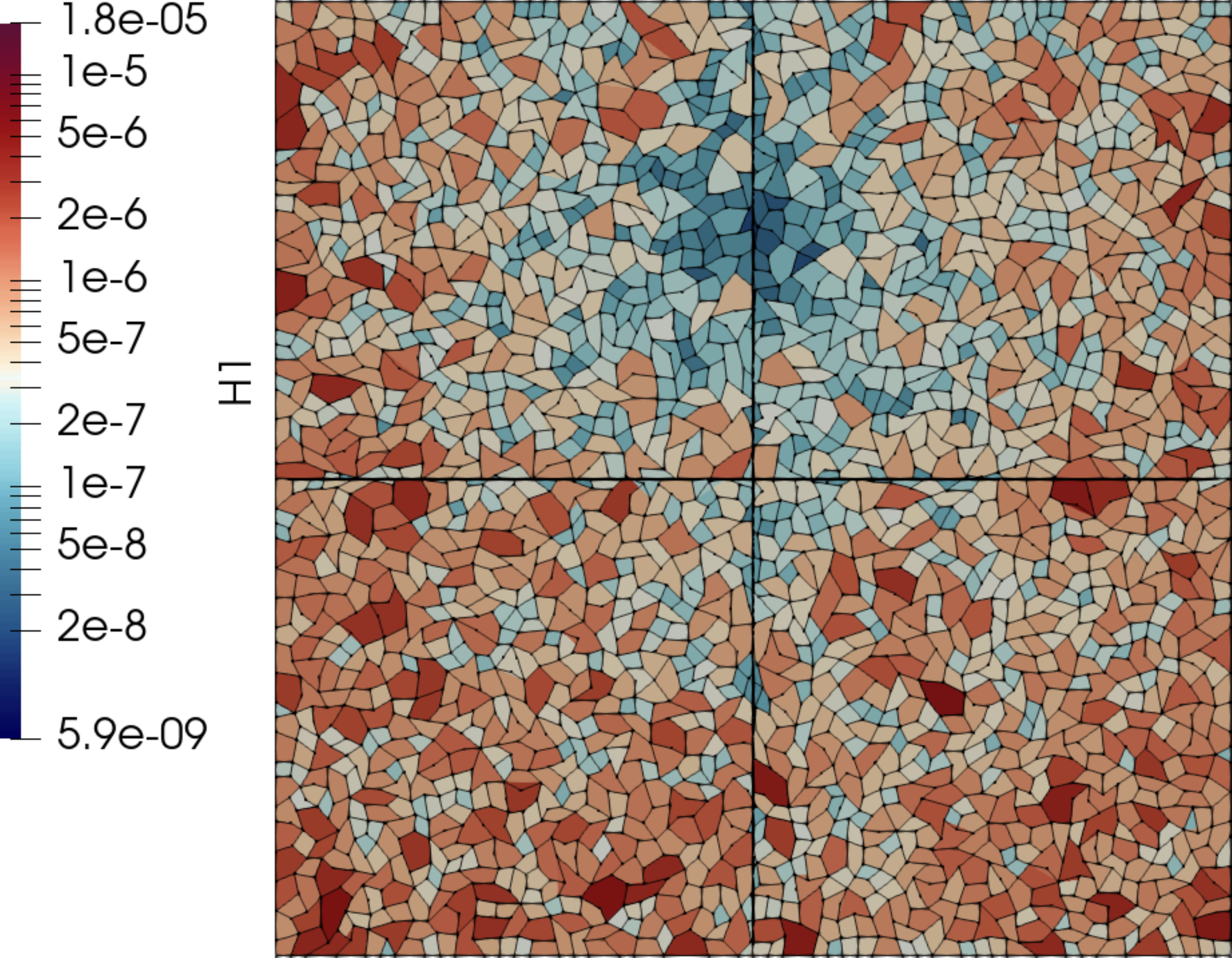}}
    \caption{Network 1, mesh M3; localization of the error $\HilbertNorm[\SobolevH{1}]{\DiscreteError}$ on each cell of $F_3$.}
    \label{plot:RUN_T5_DMN_3_H1}
\end{figure}

In particular, when the order of the VEM space increases (Figures~\ref{plot:RUN_T5_EXP_O2}-\ref{plot:RUN_T5_EXP_O3}) the measured errors in the $\SobolevH{1}$-seminorm become slightly lower than the expected ones.
This fact seems corroborated by Figure~\ref{plot:RUN_T5_DMN_3_H1}, in which we report the localization of $\HilbertNorm[\SobolevH{1}]{\DiscreteError}$ on the domain $\Fracture_3$ for the finer mesh M3.
Indeed, we can note from these images that the error roughly remains in the same order of magnitude across the agglomeration process (see the color bar at the left of each subplot), and this is more evident for the high VEM orders. 
Moreover, the blue zones, which correspond to areas with lower error values, seem to have a greater extension in the optimized meshes, meaning that the error is more evenly distributed in those cases.

To conclude the analysis, in the last columns of Table~\ref{tab:RUN_T5} we report the number NNZ of non-zero elements and the condition number cond$(\Matrixize{A})$ of the stiffness matrix associated with the discrete problem \eqref{eq:discrete}.
From the non-zero values measurements, we can assert that the optimization process does not excessively affect the sparsity pattern of matrix $\Matrixize{A}$, since the percentage values for the different $\AgglomerationParameter$ values are quite close.
On the other hand, we observe that the condition number of the stiffness matrix is under control for $\VemOrder=1,2$, and increases for $\VemOrder=3$; this growth is typical when using the basis \eqref{eq:monomial} as the polynomial base for the VEM, see \cite{ML, BERRONE201714}.
In support of this statement, in Table~\ref{tab:RUN_T5_QUALITY} we measure the $l^2$-norm quantities $\seminorm[2]{\Projection_{\VemOrder, *}^{\nabla} \mathcal{D} - \mathcal{I}}$ and $\seminorm[2]{\Projection_{\VemOrder, *}^{0} \mathcal{D} - \mathcal{I}}$ on each mesh and for each $\AgglomerationParameter$ value, as an estimate of the approximation error produced by the projectors involved in the computation of the discrete quantities.
As done in \cite{BBMR}, we denote by $\Projection_{\VemOrder, *}^{\nabla}$ and $\Projection_{\VemOrder, *}^{0}$ the matrices containing on each column the coefficients of the projection $\Projection_{\VemOrder}^{\nabla} \varphi_i$ and $\Projection_{\VemOrder}^{0} \varphi_i$ respectively of each basis function $\varphi_i$ of the VEM space $\SolutionSpace^{\VemOrder}_{\MeshParameter}$; moreover, we indicate with $\mathcal{D}$ the matrix formed by the elements $\mathcal{D}_{ij} = \text{dof}_i(m_j)$, with $m_j \in \Basis_{\VemOrder}(\MeshCellTwoD)$ and the operator $\text{dof}_i$ introduced in \eqref{eq:stabilization}.
We can observe how the error produced by the projectors increases with the VEM order as expected; this effect can be mitigated by making a different choice for the basis \eqref{eq:monomial}, such as the one proposed in \cite{ML, BERRONE201714}, but this is out of the context of this work.
In conclusion, the optimization process slightly reduces the approximation errors and the condition number of $\Matrixize{A}$; this effect is related to the agglomeration operation, which is able to remove most of the small edges and the bad-shaped ``long'' elements of the original discretization, as detailed in Figure~\ref{plot:RUN_T5_DMN_3_C2_ZOOM}.

\begin{table}[htbp!]
    \caption{Network 1; discrepancy of the projection matrices $\Projection_{\VemOrder, *}^{\nabla}$ and $\Projection_{\VemOrder, *}^{0}$; note that for $\VemOrder < 3$ we have $\Projection_{\VemOrder, *}^{\nabla} = \Projection_{\VemOrder, *}^{0}$.}
    \centering
    \includegraphics[width=0.85\textwidth]{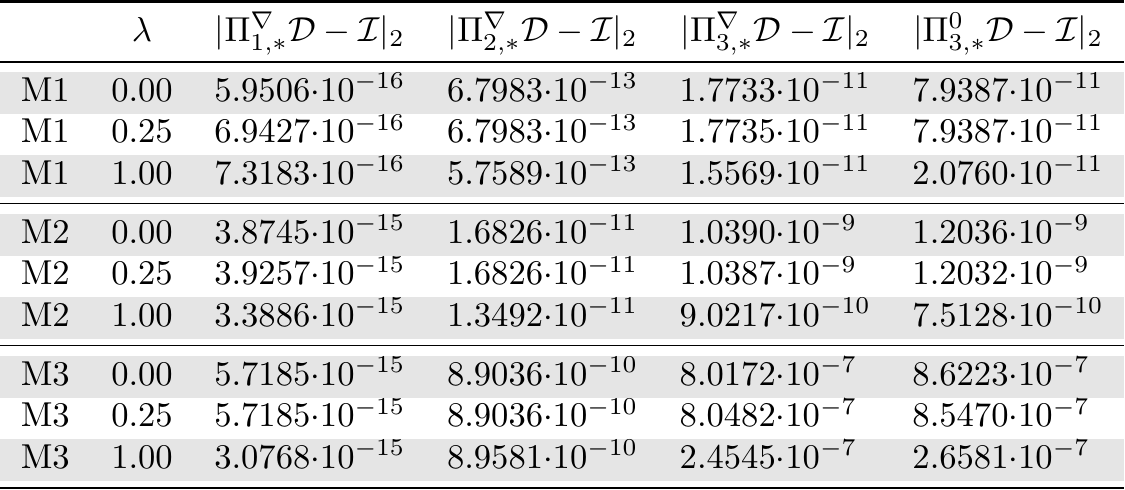}
    \label{tab:RUN_T5_QUALITY}
\end{table}

\subsection{Network 2 - Realistic DFN problem}
\label{sec:results:2}
\begin{figure}[htbp!]
    \centering
    \includegraphics[width=0.6\textwidth]{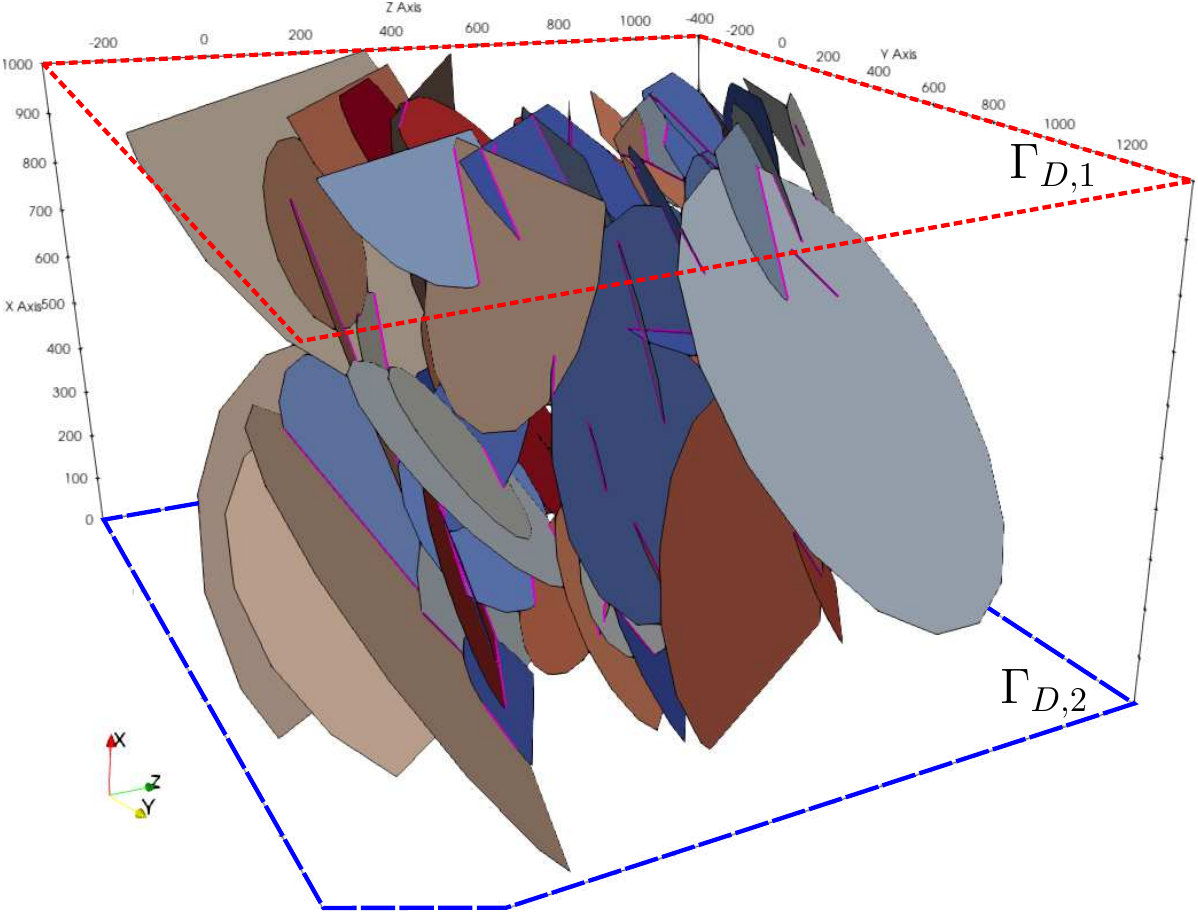}
    \caption{Network 2; domain description with interfaces $S_m$ highlighted in magenta.}
    \label{plot:RUN_T7_DMN}
\end{figure}

The second numerical test considers the computation of the hydraulic head distribution in a realistic DFN setting.
The network, presented in Figure~\ref{plot:RUN_T7_DMN}, is randomly generated inside the box $[0, 1000]\times[-400, 1400]\times[-350, 1200]$ following the strategy proposed in \cite{JOUR}, and it is composed by $86$ rounded fractures and $159$ interfaces.

\begin{figure}[htbp!]
    \subfloat[$\AgglomerationParameter = 0.0$\label{plot:RUN_T7_D3_C0}]{\includegraphics[width=0.48\textwidth]{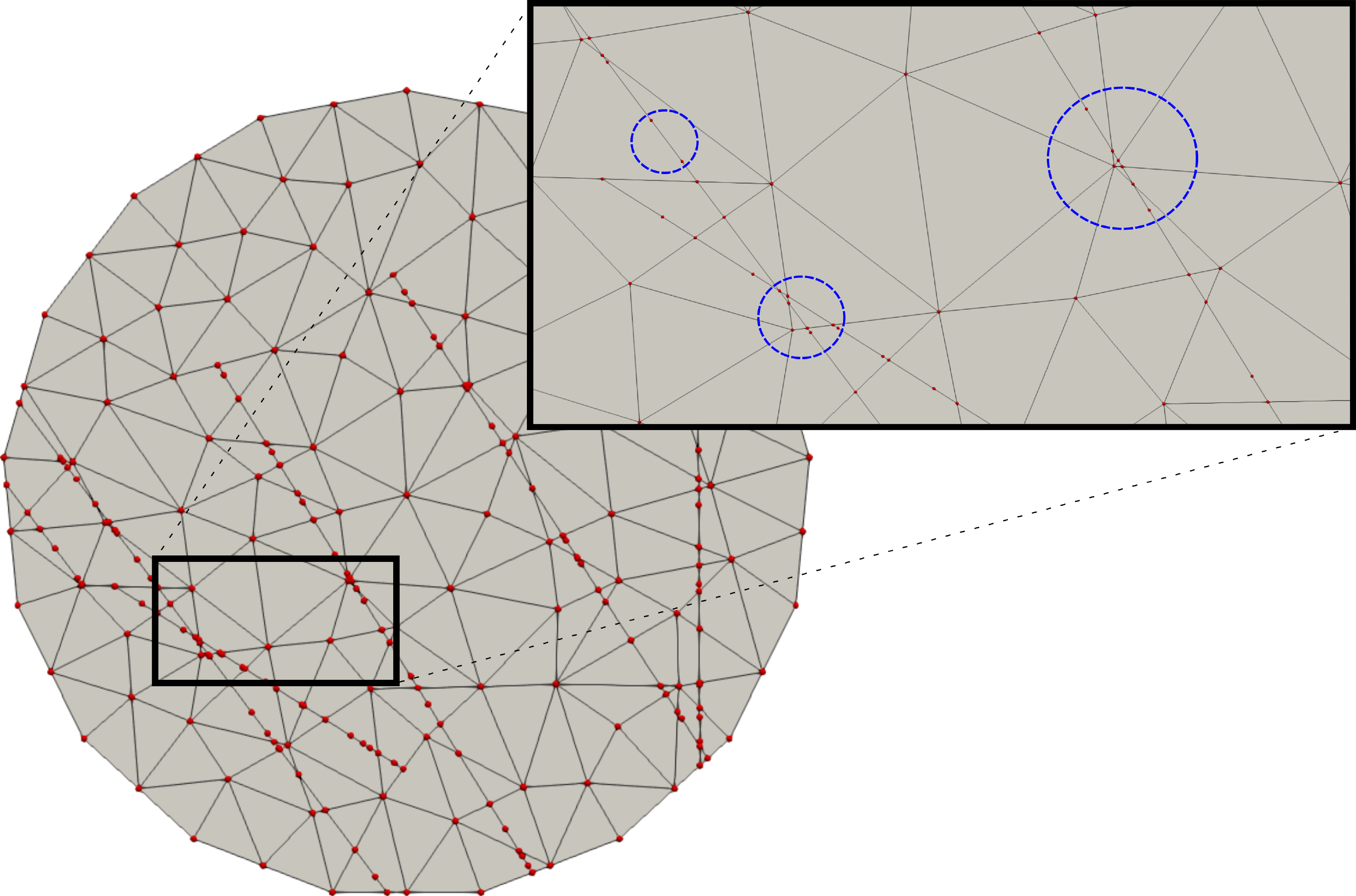}}\
    \subfloat[$\AgglomerationParameter = 1.0$\label{plot:RUN_T7_D3_C2}]{\includegraphics[width=0.48\textwidth]{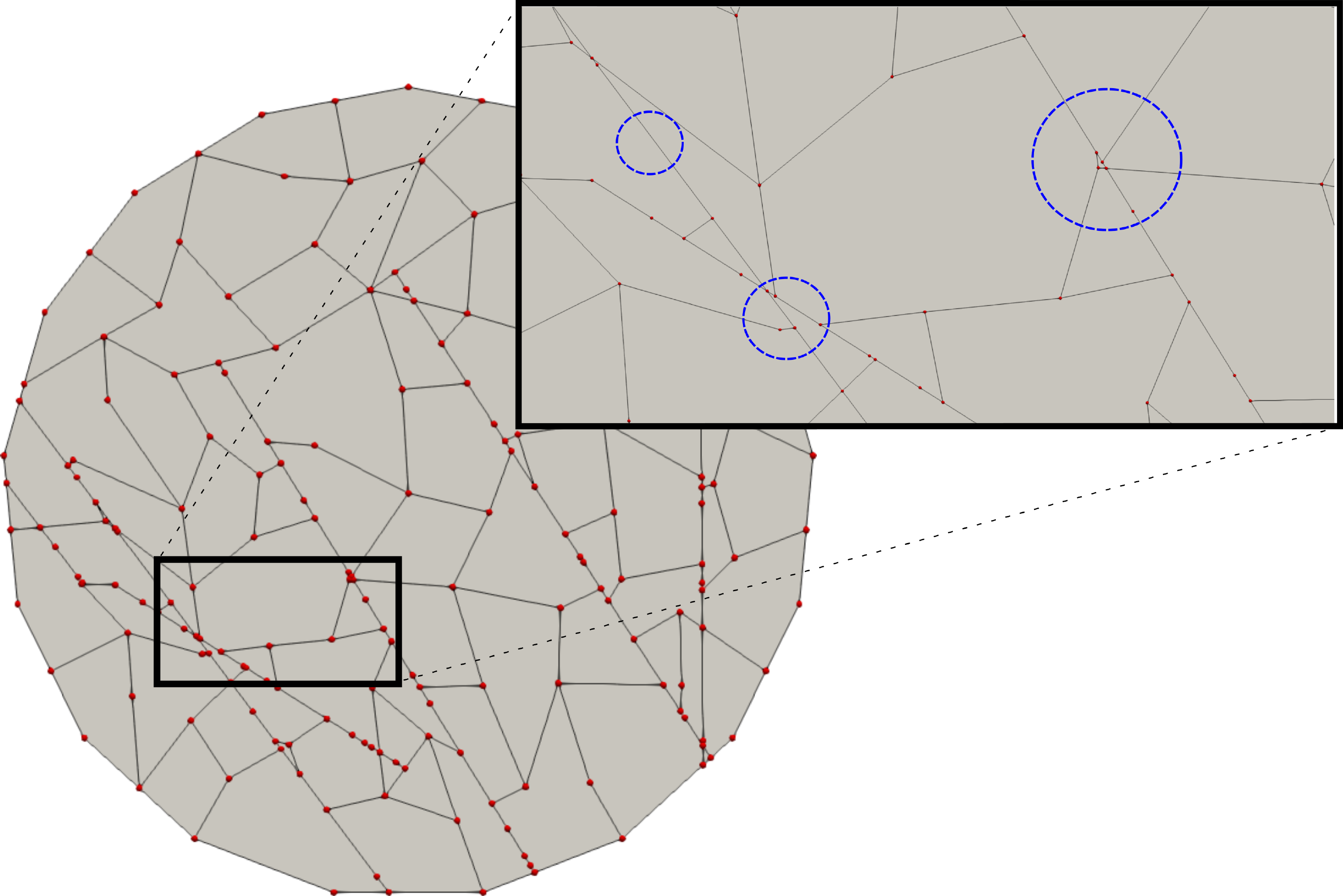}}
    \caption{Network 2, mesh M1; small edges and cells in the original fracture $\Fracture_3$ (a) are removed after the optimization process with $\AgglomerationParameter=1.0$ (b). Blue circles compare three example of mesh quality improvements.}
    \label{plot:RUN_T7_D3}
\end{figure}

We generate two meshes, labeled M1 and M2, by discretizing each fracture $\Fracture_i$ with a fixed number of triangular elements ($100$ and $200$, respectively).
As in Section~\ref{sec:results:1}, we optimize the meshes through the Quality Agglomeration algorithm of Section~\ref{sec:agglomeration}.
Figure~\ref{plot:RUN_T7_D3} shows the effects of the agglomeration over a fracture of mesh M1: a lot of small edges and cells are merged into coarser and more regular elements.

\begin{table}[htbp!]
    \caption{Network 2; analysis of the numerical errors for the different meshes and VEM orders $\VemOrder$. $\AgglomerationParameter$ is the agglomeration parameter, $\seminorm{\Mesh}$ the number of elements in the mesh, $\HydraulicHead$ the discrete solution, $\Matrixize{A}$ the stiffness matrix and NNZ its non-zero elements.}
    \centering
    \begin{subtable}{.8\textwidth}
        \centering
        \caption{$\VemOrder = 1$}
        \includegraphics[width=\textwidth]{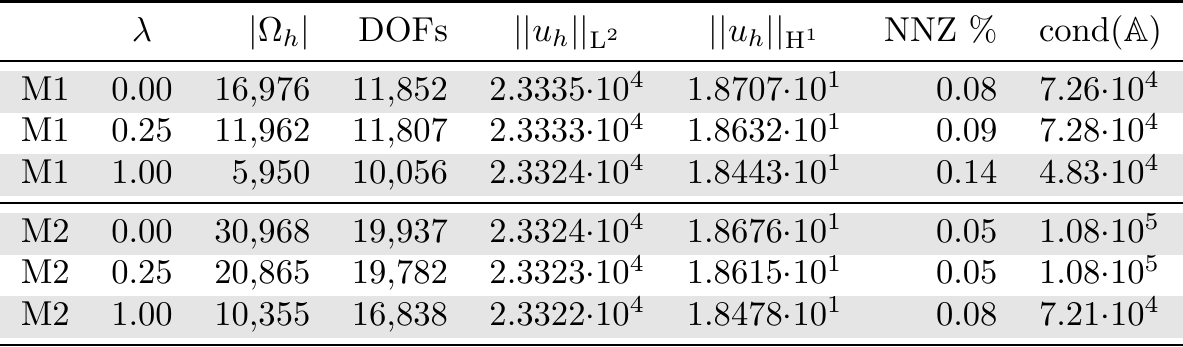}
        \label{tab:RUN_T7_O1}
    \end{subtable}
    \begin{subtable}{.8\textwidth}
        \centering
        \caption{$\VemOrder = 2$}
        \includegraphics[width=\textwidth]{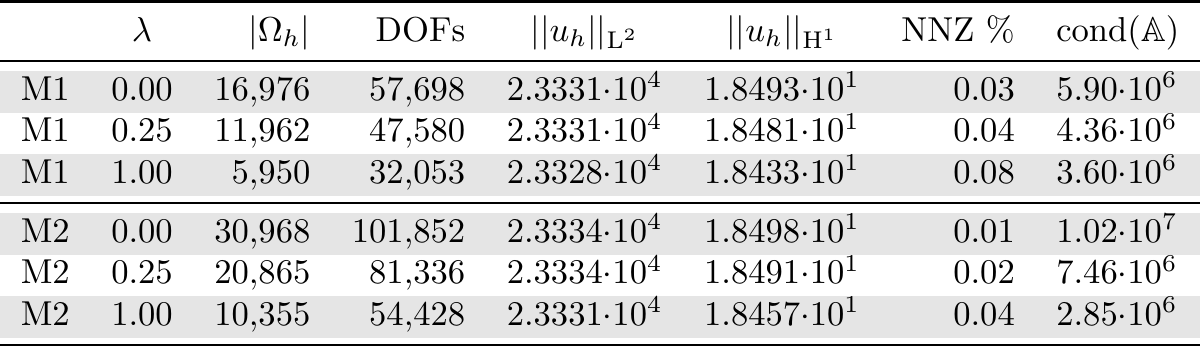}
        \label{tab:RUN_T7_O2}
    \end{subtable}
    \begin{subtable}{.8\textwidth}
        \centering
        \caption{$\VemOrder = 3$}
        \includegraphics[width=\textwidth]{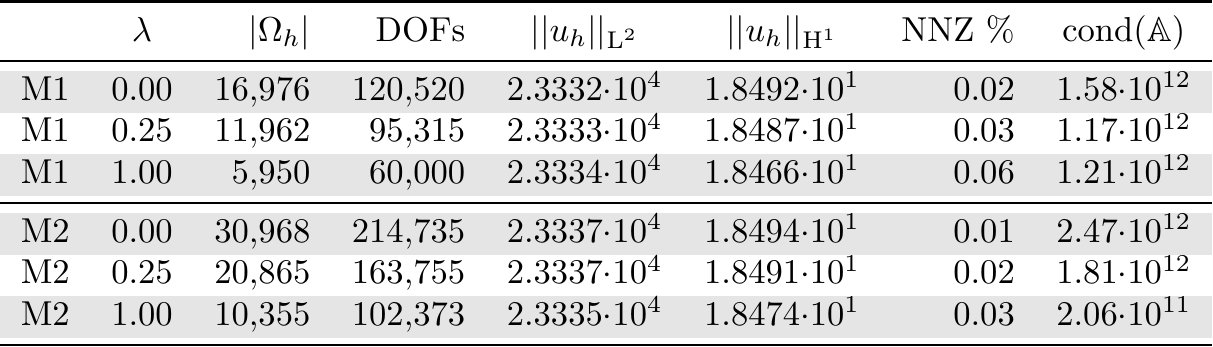}
        \label{tab:RUN_T7_O3}
    \end{subtable}
    \label{tab:RUN_T7}
\end{table}

Details on the number of elements and DOFs are given in Table~\ref{tab:RUN_T7} and are quite similar to those of Table~\ref{tab:RUN_T5}.
The number of elements $\seminorm{\Mesh}$ decreases by $30\%$ with $\AgglomerationParameter=0.25$ and by $70\%$ with $\AgglomerationParameter=1.0$; again, the number of DOFs does not change sognificantly in the case $k=1$, while for $k>1$ it decreases by $20\%$ and $50\%$ with $\AgglomerationParameter=0.25$ and 1.0, respectively.

In Table~\ref{tab:RUN_T7_ENERGY} we report the values of the energy functional $\Energy$ obtained for a generic fracture of the network ($\Fracture_3$, presented in Figure~\ref{plot:RUN_T7_D3}).
The results achieved by the optimization algorithm in terms of number of iterations, saved energy, and different energy contributions, are comparable to the ones obtained in Table~\ref{tab:RUN_T5_ENERGY}.
The measurements suggest that the optimization algorithm is robust with respect to the geometric complexity of the fracture, since the shape and the interfaces of $\Fracture_3$ are highly more complex, compared to those of the fractures in Network 1.
\begin{table}[htbp!]
    \caption{Network 2; energy functional over $\Fracture_3$ before ($\calE_1$) and after ($\calE_2$) the optimization, with the detail of the data cost ($\calE_{dc}$) and smoothness cost ($\calE_{sc}$) contribution.
    Column \textit{It} shows the number of iterations needed to converge.}
    \centering
    \includegraphics[width=0.65\textwidth]{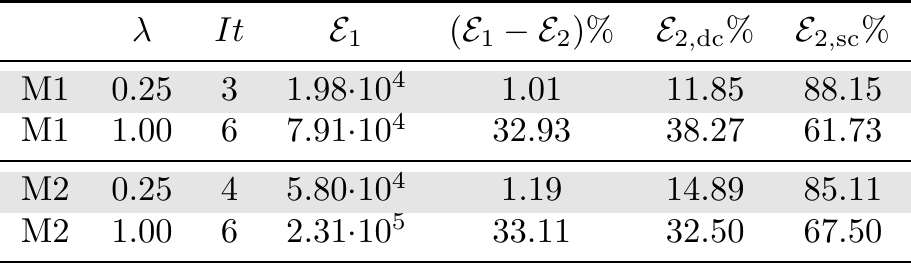}
    \label{tab:RUN_T7_ENERGY}
\end{table}

We set the numerical problem by imposing two Dirichlet boundary conditions on the mesh edges which intersect the top and the bottom planes of the box:
\begin{align*}
\begin{cases}
    \HydraulicHead = 0.0
    & \mbox{ on } \Gamma_{D,1} := \{(x,y,z) \in \RSet[3] : x=1000\} \\
    \HydraulicHead = 10.0 
    & \mbox{ on } \Gamma_{D,2} := \{(x,y,z) \in \RSet[3] : x=0\}
\end{cases}
\end{align*}
We recall that homogeneous Neumann boundary conditions are imposed on the other borders.

\begin{figure}[htbp!]
    \subfloat[\label{plot:RUN_T7_SOL_1}]{\includegraphics[width=0.52\textwidth]{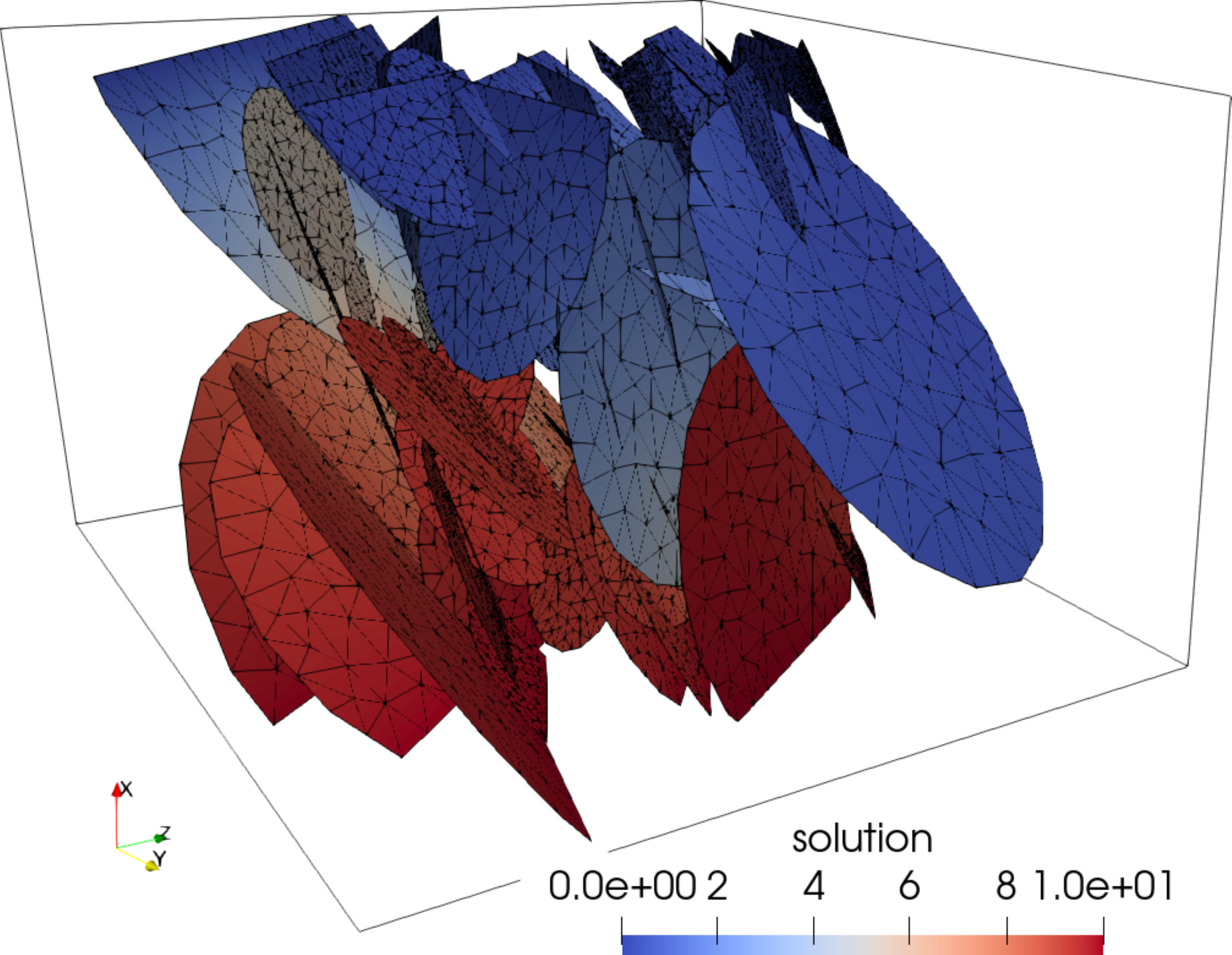}}
    \subfloat[\label{plot:RUN_T7_SOL_2}]{\includegraphics[width=0.46\textwidth]{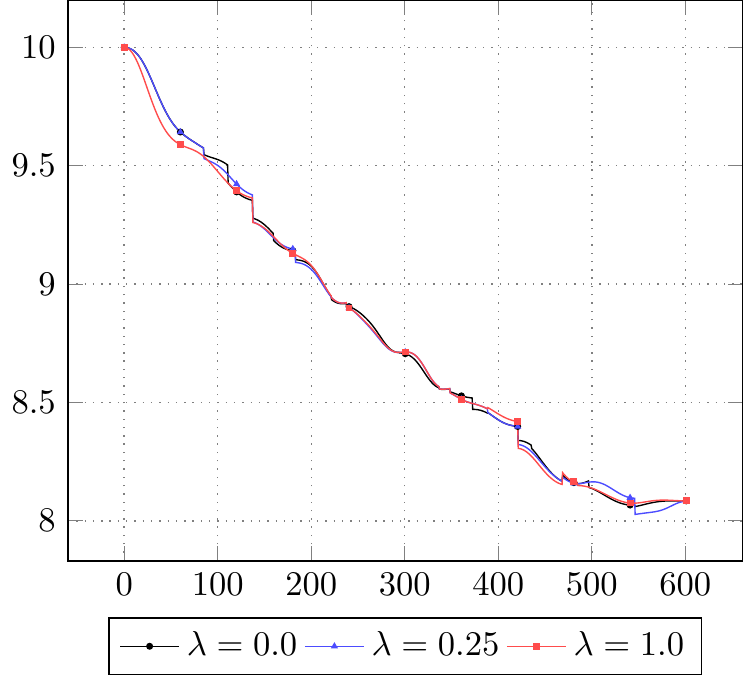}}
    \caption{Network 2; (a) discrete solution $\HydraulicHead_{\MeshParameter}$ for $\VemOrder = 1$ on mesh M1; (b) $\HydraulicHead_{\MeshParameter}$ for $\VemOrder = 3$ over a line on fracture $F_{49}$ of mesh M2.}
    \label{plot:RUN_T7_SOL}
\end{figure}

The numerical solution for $\VemOrder = 1$ on the coarser mesh M1 is shown in Figure~\ref{plot:RUN_T7_SOL_1}.
As no exact solution is available for this network, we measure the $\SobolevL{2}$-norm and the $\SobolevH{1}$-seminorm of the discrete solution $\HydraulicHead_{\MeshParameter}$ obtained on the whole network as a quantitative benchmark indicator.
Values of $\HilbertNorm[\SobL{\Network}]{\HydraulicHead_{\MeshParameter}}$ and $\HilbertNorm[\SobH{1}{\Network}]{\HydraulicHead_{\MeshParameter}}$ are reported in Table~\ref{tab:RUN_T7}.

In Figure~\ref{plot:RUN_T7_SOL_2} we perform a qualitative comparison, in which we report the value of $\HydraulicHead_{\MeshParameter}$ (computed for $\VemOrder = 3$) on a generic line in the middle of the random fracture $\Fracture_{49}$ for each optimization value $\AgglomerationParameter$.
No relevant differences can be identified between the solution on the non-optimized mesh and the one computed on the agglomerated one, despite the huge DOFs reduction (up to $50 \%$) measured in the DOFs column of Table~\ref{tab:RUN_T7}.

The last analysis is devoted to the stiffness matrix $\Matrixize{A}$ of the discrete problem \eqref{eq:discrete}, measuring its non-zero elements and its condition number in Table~\ref{tab:RUN_T7}.
We can assert that the sparsity pattern of $\Matrixize{A}$ is even less influenced by the optimization process than in the test of Section~\ref{sec:results:1}, due to the presence of a higher number of DOFs.
We immediately notice how the condition numbers rapidly increase for higher orders of the method, even faster than in the tests with Network $1$.
This effect is ascribed to the high complexity of the fracture intersections caused by the randomness of Network 2, which leads to a conforming mesh with very small edges and elongated cells, as highlighted in Figure~\ref{plot:RUN_T7_D3_C0}.
Such pathological elements are partially removed during the agglomeration process, as shown in Figure~\ref{plot:RUN_T7_D3_C2}, and we can observe a remarkable reduction of cond$(\Matrixize{A})$ after the optimization, see for instance mesh M2 in Table~\ref{tab:RUN_T7_O3}.

\begin{table}[htbp!]
    \caption{Network 2; discrepancy of the projection matrices $\Projection_{\VemOrder, *}^{\nabla}$ and $\Projection_{\VemOrder, *}^{0}$; note that for $\VemOrder < 3$ we have $\Projection_{\VemOrder, *}^{\nabla} = \Projection_{\VemOrder, *}^{0}$.}
    \centering
    \includegraphics[width=0.85\textwidth]{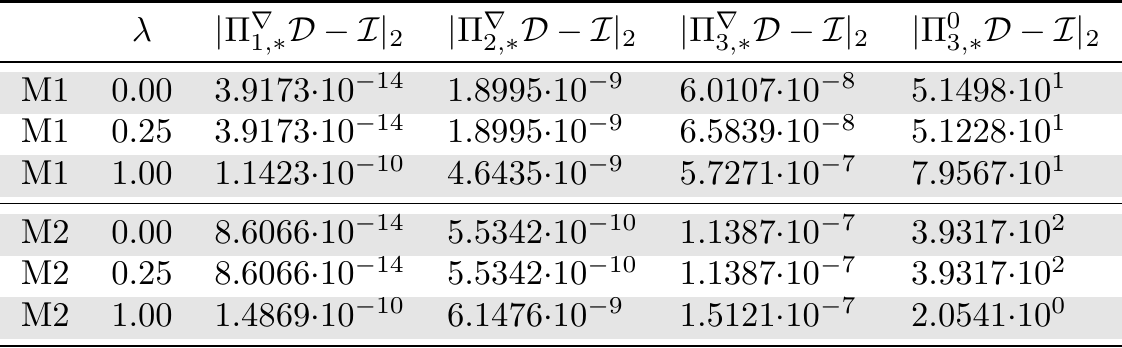}
    \label{tab:RUN_T7_QUALITY}
\end{table}

\begin{figure}[htbp!]
    \subfloat[\label{plot:RUN_T7_VEM_PI0}]{\includegraphics[width=0.5\textwidth]{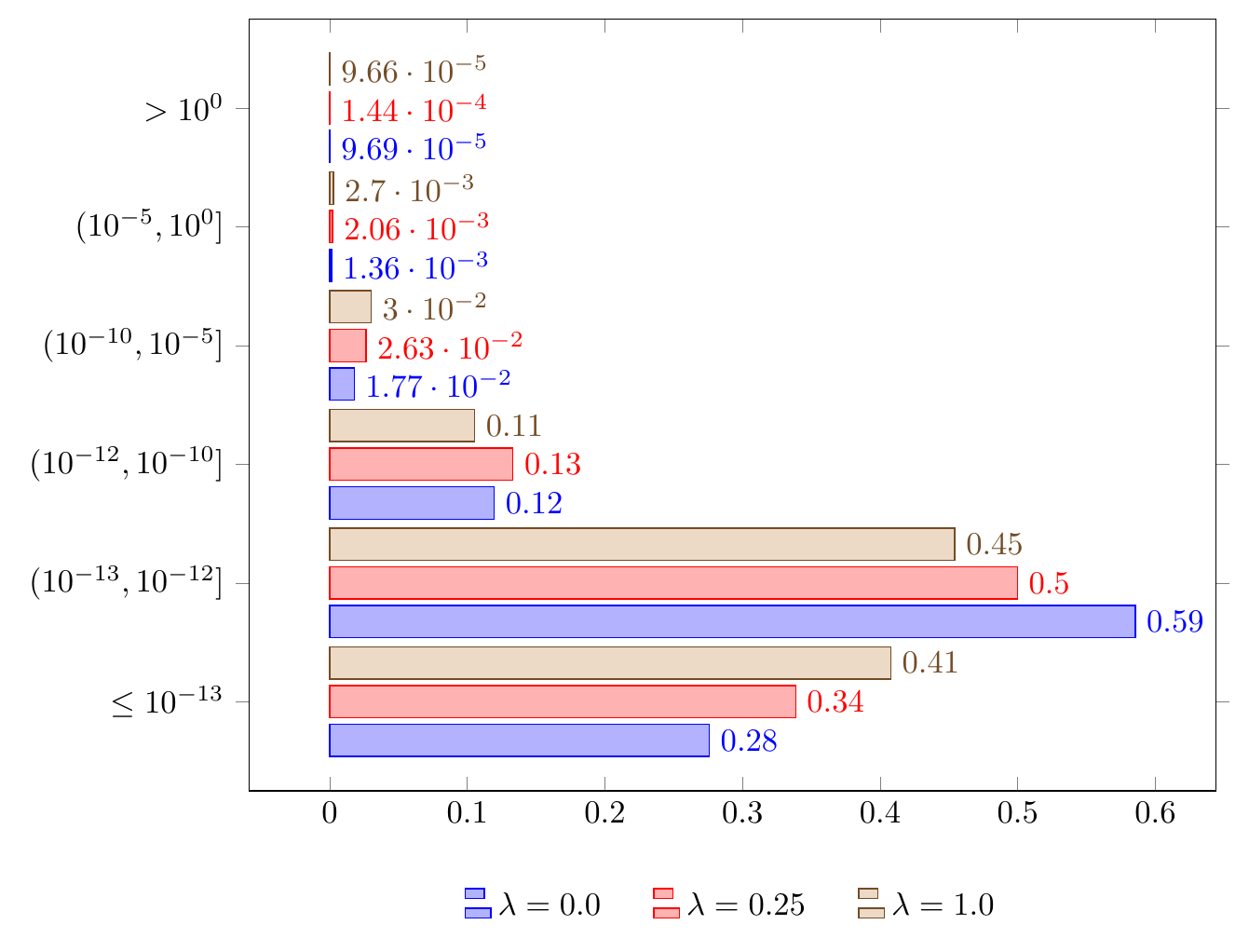}}
    \subfloat[\label{plot:RUN_T7_VEM_PIN}]{\includegraphics[width=0.5\textwidth]{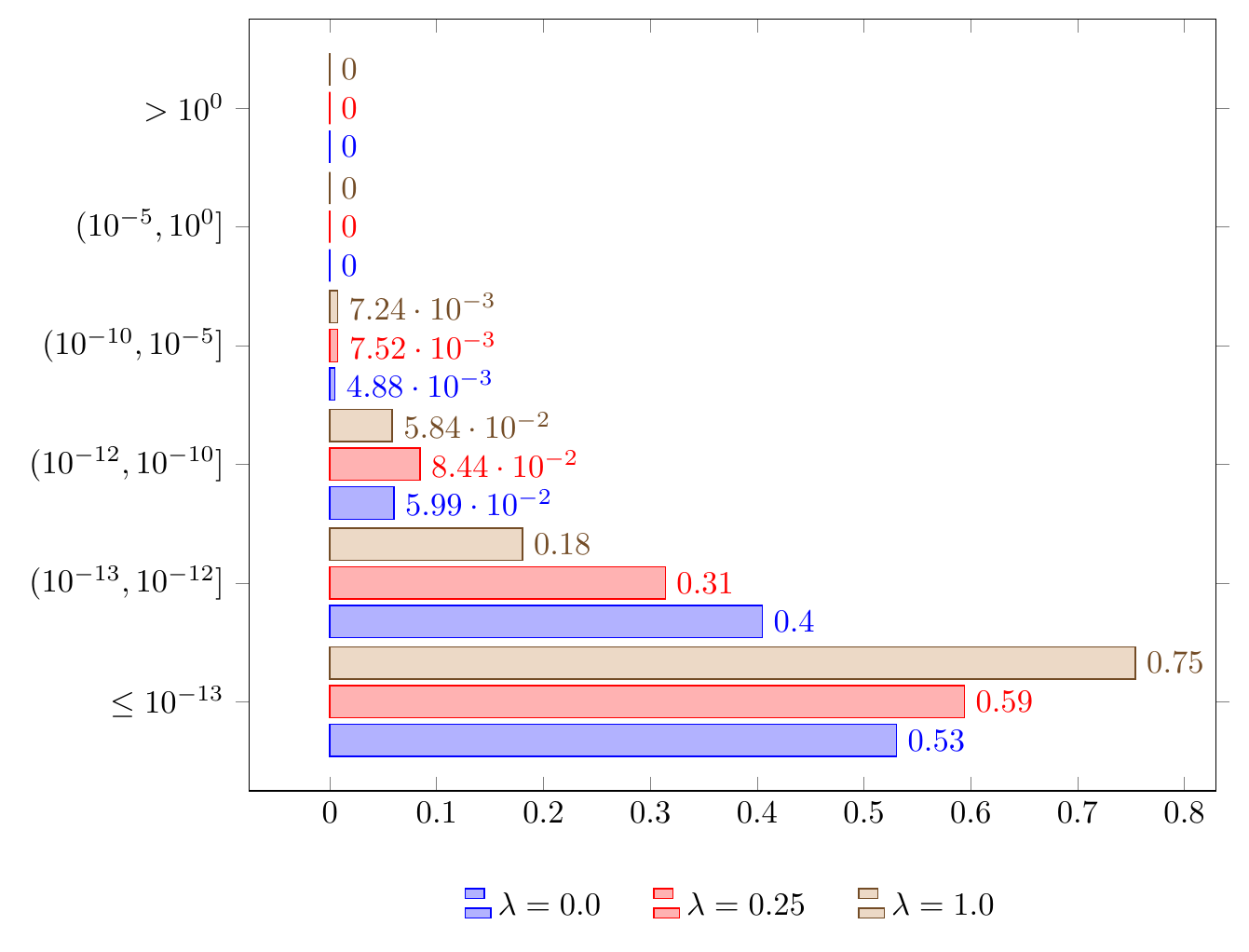}}
    \caption{Network 2, mesh M2; statistic distribution of $\seminorm[2]{\Projection_{3, *}^{0} \mathcal{D} - \mathcal{I}}$ (a) and $\seminorm[2]{\Projection_{3, *}^{\nabla} \mathcal{D} - \mathcal{I}}$  (b), with number of mesh elements reported at the end of each bar.}
    \label{plot:RUN_T7_VEM}
\end{figure}

To conclude, we show in Table~\ref{tab:RUN_T7_QUALITY} the performance of VEM projectors $\Projection_{\VemOrder}^{\nabla}$ and $\Projection_{\VemOrder}^{0}$, measuring the identities which estimate the approximation errors produced by the projection operation.
While for VEM order $\VemOrder = 1$ and $\VemOrder = 2$ the errors are acceptable, the case $\VemOrder = 3$ requires further analysis.
Thus, in Figure~\ref{plot:RUN_T7_VEM} we compare the distributions of the indicators $\seminorm[2]{\Projection_{3, *}^{\nabla} \mathcal{D} - \mathcal{I}}$ and $\seminorm[2]{\Projection_{3, *}^{0} \mathcal{D} - \mathcal{I}}$ on the elements of mesh M2, for the different values of $\AgglomerationParameter$.
From the reported data we can see how the projection errors remain under control in the vast majority of the mesh, with very few elements (less than $0.1 \%$) with projection error greater than $10^{-5}$.
As for the previous test, we are confident that those elements are curable by making a different choice for the basis \eqref{eq:monomial}.
Finally, the plots clearly highlight that the local approximation error of the projectors decreases thanks to the optimization process, leading to an improvement of the global condition number and to a more reliable solution.

\section{Conclusions}
\label{sec:conclusions}
The results obtained with the proposed mesh optimization approach in conjunction with the primal VEM-based formulation applied on DFNs are encouraging.
We remark that the Quality Agglomeration strategy presented is performed independently on each fracture of the network, thus its application in larger and more complex DFN configurations is perfectly feasible.

The Quality Agglomeration strategy has proved to be equally effective in both the simple and the complex DFN setting, being able to reduce the total number of mesh elements (up to $65\%$) and the number of DOFs (up to $50\%$), while preserving the VEM optimal convergence rates.
The former effect leads to a remarkable reduction in the computational effort for the discrete system assembly process; the latter effect translates into a significant gain in terms of computational cost to obtain the numerical solution.
Moreover, a slight improvement in the quality of the discretization was observed, both locally on the VEM projector approximation errors measured on each mesh element and globally on the condition number of the stiffness matrix.

A possible limitation of the agglomeration algorithm is that it only compares pairs of elements, i.e. given an element $\P$ with neighboring elements $\P'$ and $\P''$, the algorithm computes separately the potential quality of $\P\cup\P'$ and that of $\P\cup\P''$ not considering the total quality $\P\cup\P'\cup\P''$.
Therefore, some small mesh elements still persist in the optimized final mesh, particularly in the neighborhood of fracture intersection segments.
This issue could be partially controlled at a higher computational cost by merging the elements after each graph-cut minimization iteration.
However, the numerical results presented in this work show that this action is likely not necessary.  

In conclusion, we advise the application of the proposed optimization process in the computation of the solution for time-dependent/parametric problems, in which the computational cost reduction would reflect in each time iteration/parameter value resolution.

\section{Acknowledgements}
  Fabio Vicini and Stefano Berrone have been financially supported by
  the MIUR project ``Dipartimenti di Eccellenza 2018-2022'' (CUP
  E11G18000350001), by PRIN project ``Advanced polyhedral
  discretisations of heterogeneous PDEs for multiphysics problems''
  20204LN5N5\_003 and by INdAM-GNCS.
  Tommaso Sorgente, Silvia Biasotti, Gianmarco Manzini and Michela
  Spagnuolo have been financially supported by the ERC Project CHANGE,
  which has received funding from the European Research Council (ERC)
  under the European Union's Horizon 2020 research and innovation
  program (grant agreement no.~694515).
 Fabio Vicini, Stefano Berrone and Gianmarco Manzini are members of INdAM
  research group GNCSn.



\begin{thebibliography}{10}

\bibitem{Adams-Fournier:2003}
R~A Adams and J~J~F Fournier.
\newblock {\em Sobolev spaces}.
\newblock Pure and Applied Mathematics. Academic Press, Amsterdam, 2 edition,
  2003.

\bibitem{AHMED201749}
E~Ahmed, J~Jaffr\'e, and JE~Roberts.
\newblock A reduced fracture model for two-phase flow with different rock
  types.
\newblock {\em Mathematics and Computers in Simulation}, 137:49--70, 2017.
\newblock {MAMERN} {VI}-2015: 6th International Conference on Approximation
  Methods and Numerical Modeling in Environment and Natural Resources.

\bibitem{alliez2008recent}
P~Alliez, G~Ucelli, C~Gotsman, and M~Attene.
\newblock Recent advances in remeshing of surfaces.
\newblock {\em Shape Analysis and Structuring}, pages 53--82, 2008.

\bibitem{antonietti2016mimetic}
PF~Antonietti, L~Formaggia, A~Scotti, M~Verani, and N~Verzott.
\newblock Mimetic finite difference approximation of flows in fractured porous
  media.
\newblock {\em ESAIM: Mathematical Modelling and Numerical Analysis},
  50(3):809--832, 2016.

\bibitem{BPVEM}
L~{Beir\~{a}o da Veiga}, F~Brezzi, A~Cangiani, G~Manzini, LD~Marini, and
  A~Russo.
\newblock Basic principles of virtual element methods.
\newblock {\em Mathematical Models and Methods in Applied Sciences},
  23(01):199--214, 2013.

\bibitem{BBMR}
L~{Beir\~{a}o~da~Veiga}, F~Brezzi, LD~Marini, and A~Russo.
\newblock The hitchhiker's guide to the virtual element method.
\newblock {\em Mathematical Models and Methods in Applied Sciences},
  24(8):1541--1573, 2014.

\bibitem{VEMGSO}
L~{Beir\~{a}o~da~Veiga}, F~Brezzi, LD~Marini, and A~Russo.
\newblock Virtual element method for general second-order elliptic problems on
  polygonal meshes.
\newblock {\em Mathematical Models and Methods in Applied Sciences},
  26(4):729--750, 2016.

\bibitem{BeiraodaVeiga-Lipnikov-Manzini:2014}
L~{Beir\~ao~da~Veiga}, K~Lipnikov, and G~Manzini.
\newblock {\em The mimetic finite difference method}, volume~11 of {\em
  Modeling, Simulations and Applications}.
\newblock Springer, {I} edition, 2014.

\bibitem{BLR}
L~{Beir\~{a}o~da~Veiga}, C~Lovadina, and A~Russo.
\newblock Stability analysis for the virtual element method.
\newblock {\em Mathematical Models and Methods in Applied Sciences},
  27(13):2557--2594, 2017.

\bibitem{beirao2022sharper}
L~{Beir{\~a}o~da~Veiga} and G~Vacca.
\newblock Sharper error estimates for virtual elements and a bubble-enriched
  version.
\newblock {\em SIAM Journal on Numerical Analysis}, 60(4):1853--1878, 2022.

\bibitem{BENEDETTO2014135}
MF~Benedetto, S~Berrone, S~Pieraccini, and S~Scial{\`o}.
\newblock The virtual element method for discrete fracture network simulations.
\newblock {\em Computer Methods in Applied Mechanics and Engineering},
  280:135--156, 2014.

\bibitem{benedetto2016globally}
MF~Benedetto, S~Berrone, and S~Scial{\`o}.
\newblock A globally conforming method for solving flow in discrete fracture
  networks using the virtual element method.
\newblock {\em Finite Elements in Analysis and Design}, 109:23--36, 2016.

\bibitem{BERRONE201714}
S~Berrone and A~Borio.
\newblock Orthogonal polynomials in badly shaped polygonal elements for the
  virtual element method.
\newblock {\em Finite Elements in Analysis and Design}, 129:14--31, 2017.

\bibitem{BERRONE2021103502}
S~Berrone, A~Borio, and A~D'Auria.
\newblock Refinement strategies for polygonal meshes applied to adaptive {VEM}
  discretization.
\newblock {\em Finite Elements in Analysis and Design}, 186:103502, 2021.

\bibitem{BERRONE2019904}
S~Berrone, A~Borio, and F~Vicini.
\newblock Reliable a posteriori mesh adaptivity in discrete fracture network
  flow simulations.
\newblock {\em Computer Methods in Applied Mechanics and Engineering},
  354:904--931, 2019.

\bibitem{BERRONE2022103770}
S~Berrone and A~D'Auria.
\newblock A new quality preserving polygonal mesh refinement algorithm for
  polygonal element methods.
\newblock {\em Finite Elements in Analysis and Design}, 207:103770, 2022.

\bibitem{berrone2021three}
S~Berrone, D~Grappein, S~Pieraccini, and S~Scial{\'o}.
\newblock A three-field based optimization formulation for flow simulations in
  networks of fractures on nonconforming meshes.
\newblock {\em SIAM Journal on Scientific Computing}, 43(2):B381--B404, 2021.

\bibitem{BPST}
S~Berrone, S~Pieraccini, and S~Scial\`{o}.
\newblock A {PDE}-constrained optimization formulation for discrete fracture
  network flows.
\newblock {\em SIAM Journal on Scientific Computing}, 35(2):B487--B510, 2013.

\bibitem{BSV}
S~Berrone, S~Scial\`{o}, and F~Vicini.
\newblock Parallel meshing, discretization, and computation of flow in massive
  discrete fracture networks.
\newblock {\em SIAM Journal on Scientific Computing}, 41(4):C317--C338, 2019.

\bibitem{boykov2004experimental}
Y~Boykov and V~Kolmogorov.
\newblock An experimental comparison of min-cut/max-flow algorithms for energy
  minimization in vision.
\newblock {\em IEEE Transactions on Pattern Analysis and Machine Intelligence},
  26(9):1124--1137, 2004.

\bibitem{boykov2001fast}
Y~Boykov, O~Veksler, and R~Zabih.
\newblock Fast approximate energy minimization via graph cuts.
\newblock {\em IEEE Transactions on Pattern Analysis and Machine Intelligence},
  23(11):1222--1239, 2001.

\bibitem{BURMAN2019726}
E~Burman, P~Hansbo, MG~Larson, and K~Larsson.
\newblock Cut finite elements for convection in fractured domains.
\newblock {\em Computers \& Fluids}, 179:726--734, 2019.

\bibitem{CMSO}
A~Cangiani, G~Manzini, and OJ~Sutton.
\newblock {Conforming and nonconforming virtual element methods for elliptic
  problems}.
\newblock {\em IMA Journal of Numerical Analysis}, 37(3):1317--1354, 08 2016.

\bibitem{chalmeta2013measuring}
R~Chalmeta, F~Hurtado, V~Sacrist{\'a}n, and M~Saumell.
\newblock Measuring regularity of convex polygons.
\newblock {\em Computer-Aided Design}, 45(2):93--104, 2013.

\bibitem{chave2018hybrid}
F~Chave, DA~Di~Pietro, and L~Formaggia.
\newblock A hybrid high-order method for {D}arcy flows in fractured porous
  media.
\newblock {\em SIAM Journal on Scientific Computing}, 40(2):A1063--A1094, 2018.

\bibitem{DF}
WS~Dershowitz and C~Fidelibus.
\newblock Derivation of equivalent pipe network analogues for three-dimensional
  discrete fracture networks by the boundary element method.
\newblock {\em Water Resources Research}, 35(9):2685--2691, 1999.

\bibitem{erten2009mesh}
H~Erten, A~{\"U}ng{\"o}r, and C~Zhao.
\newblock Mesh smoothing algorithms for complex geometric domains.
\newblock In {\em Proceedings of the 18th International Meshing Roundtable},
  pages 175--193. Springer, 2009.

\bibitem{fumagalli2019numerical}
A~Fumagalli and I~Berre.
\newblock {\em Numerical methods for processes in fractured porous media}.
\newblock Birkh\"auser, Cham, 2019.

\bibitem{fumagalli2018dual}
A~Fumagalli and E~Keilegavlen.
\newblock Dual virtual element method for discrete fractures networks.
\newblock {\em SIAM Journal on Scientific Computing}, 40(1):B228--B258, 2018.

\bibitem{goldberg1988new}
AV~Goldberg and RE~Tarjan.
\newblock A new approach to the maximum-flow problem.
\newblock {\em Journal of the ACM}, 35(4):921--940, 1988.

\bibitem{eigenweb}
G~Guennebaud and B~Jacob.
\newblock Eigen v3.
\newblock http://eigen.tuxfamily.org, 2010.

\bibitem{PFA}
F~H{\'e}din, G~Pichot, and A~Ern.
\newblock A hybrid high-order method for flow simulations in discrete fracture
  networks.
\newblock In Fred~J. Vermolen and Cornelis Vuik, editors, {\em Numerical
  Mathematics and Advanced Applications ENUMATH 2019}, pages 521--529, Cham,
  2021. Springer International Publishing.

\bibitem{huang2020anisotropic}
W~Huang and Y~Wang.
\newblock Anisotropic mesh quality measures and adaptation for polygonal
  meshes.
\newblock {\em Journal of Computational Physics}, 410:109368, 2020.

\bibitem{knupp2001algebraic}
P~Knupp.
\newblock Algebraic mesh quality metrics.
\newblock {\em SIAM Journal on Scientific Computing}, 23(1):193--218, 2001.

\bibitem{knupp2012introducing}
P~Knupp.
\newblock Introducing the target-matrix paradigm for mesh optimization via
  node-movement.
\newblock {\em Engineering with Computers}, 28(4):419--429, 2012.

\bibitem{kolmogorov2004energy}
V~Kolmogorov and R~Zabin.
\newblock What energy functions can be minimized via graph cuts?
\newblock {\em IEEE Transactions on Pattern Analysis and Machine Intelligence},
  26(2):147--159, 2004.

\bibitem{KMJR}
M~K\"oppel, V~Martin, J~Jaffr\'e, and JE~Roberts.
\newblock A {L}agrange multiplier method for a discrete fracture model for flow
  in porous media.
\newblock {\em Computers \& Fluids}, 23:239--253, 2019.

\bibitem{Lipnikov-Manzini-Shashkov:2014}
K~Lipnikov, G~Manzini, and MJ~Shashkov.
\newblock Mimetic finite difference method.
\newblock {\em Journal of Computational Physics}, 257:1163--1227, 2014.

\bibitem{lo2014finite}
DSH Lo.
\newblock {\em Finite element mesh generation}.
\newblock CRC Press, 2014.

\bibitem{martin2005modeling}
V~Martin, J~Jaffr{\'e}, and JE~Roberts.
\newblock Modeling fractures and barriers as interfaces for flow in porous
  media.
\newblock {\em SIAM Journal on Scientific Computing}, 26(5):1667--1691, 2005.

\bibitem{ML}
L~Mascotto.
\newblock Ill-conditioning in the virtual element method: stabilizations and
  bases.
\newblock {\em Numerical Methods for Partial Differential Equations},
  34(4):1258--1281, 2018.

\bibitem{misztal2009tetrahedral}
MK~Misztal, JA~B{\ae}rentzen, F~Anton, and K~Erleben.
\newblock Tetrahedral mesh improvement using multi-face retriangulation.
\newblock In {\em Proceedings of the 18th International Meshing Roundtable},
  pages 539--555. Springer, 2009.

\bibitem{Brenner-Scott:2008}
L~Ridgway~Scott and SC~Brenner.
\newblock {\em The mathematical theory of finite element methods}.
\newblock Texts in Applied Mathematics 15. Springer-Verlag, New York, 3
  edition, 2008.

\bibitem{sorgente2021role}
T~Sorgente, S~Biasotti, G~Manzini, and M~Spagnuolo.
\newblock The role of mesh quality and mesh quality indicators in the virtual
  element method.
\newblock {\em Advances in Computational Mathematics}, 48(1):3, 2021.

\bibitem{sorgente2021indicator}
T~Sorgente, S~Biasotti, G~Manzini, and M~Spagnuolo.
\newblock Polyhedral mesh quality indicator for the virtual element method.
\newblock {\em Computers \& Mathematics with Applications}, 114:151--160, 2022.

\bibitem{sorgente2021vem}
T~Sorgente, D~Prada, D~Cabiddu, S~Biasotti, G~Patane, M~Pennacchio,
  S~Bertoluzza, G~Manzini, and M~Spagnuolo.
\newblock {\em VEM and the mesh}, volume~31 of {\em SEMA SIMAI Springer
  series}, chapter~1, pages 1--54.
\newblock Springer, Nature Switzerland AG, 2021.
\newblock ISBN: 978-3-030-95318-8.

\bibitem{JOUR}
S~Srinivasan, J~Hyman, S~Karra, D~O'Malley, H~Viswanathan, and G~Srinivasan.
\newblock Robust system size reduction of discrete fracture networks: a
  multi-fidelity method that preserves transport characteristics.
\newblock {\em Computational Geosciences}, 22(6):1515--1526, 2018.

\bibitem{stimpson2007verdict}
CJ~Stimpson, CD~Ernst, P~Knupp, PP~P{\'e}bay, and D~Thompson.
\newblock The {V}erdict library reference manual.
\newblock {\em Sandia National Laboratories Technical Report}, 9(6), 2007.

\bibitem{vartziotis2008mesh}
D~Vartziotis, T~Athanasiadis, I~Goudas, and J~Wipper.
\newblock Mesh smoothing using the geometric element transformation method.
\newblock {\em Computer Methods in Applied Mechanics and Engineering},
  197(45-48):3760--3767, 2008.

\bibitem{gco}
O~Veksler and A~Delong.
\newblock Multi-label optimization.
\newblock ~\\ https://vision.cs.uwaterloo.ca/code/, 2010.

\bibitem{zunic2004new}
J~Zunic and PL~Rosin.
\newblock A new convexity measure for polygons.
\newblock {\em IEEE Transactions on Pattern Analysis and Machine Intelligence},
  26(7):923--934, 2004.

\end{thebibliography}
\end{document}